\documentclass[bj-nolines,authoryear]{imsart}
\usepackage{bm}
\usepackage{bbm}
\usepackage{appendix}
\usepackage{calrsfs}
\usepackage{amsfonts}
\usepackage[mathscr]{eucal}
\usepackage{dsfont}
\usepackage{amsmath}
\usepackage{natbib}

\startlocaldefs
\theoremstyle{plain}

\newtheorem{theorem}{Theorem}[section]
\newtheorem{lemma}[theorem]{Lemma}
\newtheorem{corollary}[theorem]{Corollary}
\theoremstyle{remark}

\newcommand{\Mle}{\widehat\bTheta}
\newcommand{\mle}{\widehat\btheta}

\newcommand{\truth}{\btheta^\star}
\newcommand{\mI}{\mathcal{I}}

\newcommand{\bDelta}{\bm{\Delta}}
\newcommand{\Aavg}{A_{\textup{avg}}}

\newcommand{\avg}{\textup{avg}}

\newcommand{\interior}{\textnormal{int}}

\newcommand{\mA}{\mathscr{A}}
\newcommand{\mB}{\mathscr{B}}

\newcommand{\mD}{\mathscr{D}}

\newcommand{\mG}{\mathscr{G}}

\newcommand{\mN}{\mathscr{N}}

\newcommand{\mcC}{\mathcal{C}}

\newcommand{\mcG}{\mathcal{G}}
\newcommand{\mcH}{\mathcal{H}}
\newcommand{\mcI}{I}

\newcommand{\mcR}{\mathcal{R}}

\newcommand{\mcU}{\mathcal{U}}
\newcommand{\mcV}{\mathcal{V}}

\newcommand{\mbE}{\mathds{E}}

\newcommand{\mbM}{\mathds{M}}

\newcommand{\mbP}{\mathds{P}}

\newcommand{\mbR}{\mathds{R}}

\newcommand{\mbX}{\mathds{X}}
\newcommand{\mbY}{\mathds{Y}}

\newcommand{\bu}{\bm{u}}
\newcommand{\bv}{\bm{v}}
\newcommand{\bw}{\bm{w}}
\newcommand{\bx}{\bm{x}}
\newcommand{\by}{\bm{y}}
\newcommand{\bz}{\bm{z}}

\newcommand{\bA}{\bm{A}}

\newcommand{\bC}{\bm{C}}

\newcommand{\bG}{\bm{G}}
\newcommand{\bH}{\bm{H}}
\newcommand{\bI}{\bm{I}}

\newcommand{\bR}{\bm{R}}
\newcommand{\bS}{\bm{S}}

\newcommand{\bW}{\bm{W}}
\newcommand{\bX}{\bm{X}}
\newcommand{\bY}{\bm{Y}}
\newcommand{\bZ}{\bm{Z}}

\newcommand{\bmu}{\mbox{\boldmath$\mu$}}

\newcommand{\bTheta}{\bm\Theta}

\newcommand{\btheta}{\boldsymbol{\theta}}

\newcommand{\one}{\mathbbm{1}}

\DeclareMathOperator*{\argmax}{arg\,max}

\renewcommand{\bar}{\overline}
\newcommand{\dsum}{\displaystyle\sum\limits}

\newcommand{\dprod}{\displaystyle\prod\limits}
\newcommand{\var}{\mathds{V}}
\newcommand{\cov}{\mathds{C}}
\newcommand{\norm}[1]{|\!|#1|\!|} 
\newcommand{\mnorm}[1]{|\!|\!|#1|\!|\!|}

\newcommand{\maxw}{\widetilde{\lambda}^{\star}_{\max,W}}
\newcommand{\minw}{\widetilde{\lambda}^{\epsilon}_{\min,W}}
\newcommand{\maxb}{\widetilde{\lambda}^{\star}_{\max,B}}
\newcommand{\minb}{\widetilde{\lambda}^{\epsilon}_{\min,B}}
\newcommand{\maxwe}{\widetilde{\lambda}^{\epsilon}_{\max,W}}
\newcommand{\maxbe}{\widetilde{\lambda}^{\epsilon}_{\max,B}}
\newcommand{\minwt}{\widetilde{\lambda}^{\star}_{\min,W}}
\newcommand{\minbt}{\widetilde{\lambda}^{\star}_{\min,B}}

\newtheorem{remark}{Remark}

\newcommand{\be}{\begin{equation}\arraycolsep=4pt\begin{array}{lllllllllllllllll}}
\newcommand{\beno}{\begin{equation}\arraycolsep=4pt\begin{array}{lllllllllllll}\nonumber}
\newcommand{\ee}{\end{array}\end{equation}}

\newcommand{\bi}{\begin{itemize}}
\newcommand{\ei}{\end{itemize}}

\newcommand{\ben}{\begin{enumerate}}
\newcommand{\een}{\end{enumerate}}

\renewcommand{\=}{&=&}

\newcommand{\hide}[1]{}
\newcommand{\ghost}[1]{}

\newcommand{\s}{\vspace{0.25cm}}

\newcounter{com}

\newcounter{assumption}
\newenvironment{assumption}[1][]{\refstepcounter{assumption}\par\noindent%
\textbf{Assumption~\theassumption #1}. \rmfamily}{\medskip}

\newcommand{\nat}{\btheta}

\newcounter{ttproof}
\newcounter{llproof}
\newcommand{\llproof}{%
\addtocounter{llproof}{1}%
{}\textsc{Proof of Lemma}{}
}

\endlocaldefs

\begin{document}

\begin{frontmatter}
\title{Rates of convergence and normal approximations for estimators of local dependence random graph models}
\runtitle{Rates of convergence for local dependence random graph models}

\begin{aug}
%%%%%%%%%%%%%%%%%%%%%%%%%%%%%%%%%%%%%%%%%%%%%%%
%% ORCID can be inserted by command:         %%
%% \orcid{0000-0000-0000-0000}               %%
%%%%%%%%%%%%%%%%%%%%%%%%%%%%%%%%%%%%%%%%%%%%%%%
\author[A]{\inits{F.}\fnms{Jonathan R.}~\snm{Stewart}\ead[label=e1]{jrstewart@fsu.edu}}
%%%%%%%%%%%%%%%%%%%%%%%%%%%%%%%%%%%%%%%%%%%%%%
%% Addresses                                %%
%%%%%%%%%%%%%%%%%%%%%%%%%%%%%%%%%%%%%%%%%%%%%%
\address[A]{Department of Statistics,
Florida State University,
Tallahassee, FL, USA\printead[presep={,\ }]{e1}}
\end{aug}

\begin{abstract}
Local dependence random graph models are a class of block models for network data which allow for dependence
among edges under a local dependence assumption defined around the block structure of the network.
Since being introduced by \citet{Schweinberger2015},
research in the statistical network analysis and network science
literatures have demonstrated the potential and utility of this class of models.
In this work,
we provide the first theory for 
estimation and inference 
which ensures consistent and valid inference
of parameter vectors of local dependence random graph models. 
This is accomplished by deriving convergence rates of estimation and inference procedures
for local dependence random graph models based on a single observation of the graph,
allowing both the number of model parameters and the sizes of blocks to tend to infinity.
First, we derive non-asymptotic bounds on
the $\ell_2$-error of maximum likelihood estimators with convergence rates,
outlining conditions under which
these rates are minimax optimal.
Second,
and more importantly,
we derive non-asymptotic bounds on the
error of the multivariate normal approximation.
These theoretical results are the first
to achieve both optimal rates of convergence and non-asymptotic bounds on the error of the
multivariate normal approximation for parameter vectors of local dependence random graph models.

\end{abstract}


\end{frontmatter}

\section{Introduction} 
\label{sec:intro}

Local dependence random graph models,
introduced by \citet{Schweinberger2015}, 
are a class of statistical models for network data built around block structure,
where a population of nodes $\mN$,
which we take without loss to be $\mN \coloneqq \{1, \ldots, N\}$ ($N \geq 3$),
is partitioned 
into $K \in \{1, 2, \ldots\}$ subsets $\mA_1, \ldots, \mA_K$ called blocks 
(also referred to as communities or subpopulations within the literature). 
The class owes its name to the fundamental assumption that dependence among edges is constrained 
to block-based subgraphs. 
We formally review local dependence random graph models in Section \ref{sec:loc_dep_rg}. 

There are two key aspects to local dependence random graph models which help to explain the research interest received 
in both the statistical network analysis and network science literatures
\citep[][]{Stewart2019,SchweinbergerStewart2020,Babkin2020,Whetsell2021,Mele2022,Dahbura2021,Agneessens2022,Dahbura2023,Tolochko2024}.
First, 
block structure (or community structure) is a well-established structural phenomena relevant to many applications and 
networks encountered in our world 
\citep[e.g.,][]{Holland1983, Newman2004, Stewart2019}.
Second, 
local dependence random graph models 
possess 
desirable  properties and behavior 
that circumvent early difficulties in constructing models of edge dependence,
which include producing non-degenerate models of edge dependence (including transitivity)
and consistency results for estimators
\citep[][]{Schweinberger2015, SchweinbergerStewart2020}.

Utilization of local dependence random graph models requires knowledge or estimates 
of both the block memberships of the nodes in the network, 
as well as the parameters of interest which determine the amount of probability mass placed on different 
configurations of the network.
In practice, 
the parameter vectors of local dependence random graph models must always be estimated, 
whereas the block memberships of nodes can either be observed as part of the observation process 
\citep{Stewart2019,SchweinbergerStewart2020},
or can be estimated \citep{Babkin2020,Schweinberger2020-Bernoulli}. 
We will focus on the problem of estimating parameter vectors under the assumption that 
the block memberships of nodes have either been observed as part of the observation process or have been estimated.

In this work, 
we advance the literature on local dependence random graph models by 
providing the first statistical theory 
%disclaimers 
which elaborates conditions under which 
estimation and inference methodology 
%procedures 
%for parameter vectors 
based on a single observation of the graph 
can be expected to produce consistent and valid inference  
of parameter vectors of  local dependence random graph models. 
The main results are non-asymptotic and cover settings where the number of model parameters  
and the sizes of the blocks tend to infinity. 
The main contributions of this work include: 
\ben
\item Establishing the first non-asymptotic bounds on the $\ell_2$-error of maximum likelihood 
estimators of parameters vectors of local dependence random graph models which hold with high probability, \s
\item Outlining conditions under which the rates of convergence implied by the bounds on the $\ell_2$-error 
of maximum likelihood estimators are minimax optimal, and \s  
\item Deriving the first non-asymptotic bound on the error of the multivariate normal approximation of
a standardization of maximum likelihood estimators. 
\een
%Leveraging these contributions, 
%we are able to provide convergence rates for the statistical error of estimators
%and for the error of the multivariate normal approximation,
%enhancing existing methodological and theoretical work in the literature on these models
%\citep[e.g.,][]{Schweinberger2015,
%Stewart2019,
%Schweinberger2020-Bernoulli,
%Babkin2020,
%SchweinbergerStewart2020}.
All results are stated in terms of interpretable quantities,
allowing us to quantify the effect of key aspects of the statistical model and network structure
upon convergence rates of the aforementioned errors.
In so doing, 
we introduce the first principled approach to estimation and inference for 
local dependence random graph models
by developing theoretical results which achieve both optimal rates of convergence
and  non-asymptotic bounds on the error of the multivariate normal approximation of 
maximum likelihood estimators.

\subsection{Local dependence random graph models} 
\label{sec:loc_dep_rg}

We consider simple, undirected random graphs $\bX \in \mbX \coloneqq \{0, 1\}^{\tbinom{N}{2}}$ 
which are defined on the set of nodes 
$\mN \coloneqq \{1, \ldots, N\}$ ($N \geq 3$). 
Edge variables between pairs of nodes $\{i,j\} \subset \mN$ are given by 
\beno
X_{i,j} 
\= \begin{cases}
1 & \text{Nodes } i \text{ and } j \text{ are connected in the graph} \\ 
0 & \text{Otherwise}, 
\end{cases}
\ee
assuming throughout that $X_{i,j} = X_{j,i}$ ($\{i,j\} \subset \mN$) and $X_{i,i} = 0$ ($i \in \mN$). 
%adopting the usual convention that $X_{i,j} = X_{j,i}$ (for all $\{i,j\} \subset \mN$)
%and $X_{i,i} = 0$ (for all $i \in \mN$) with probability $1$.
%Extensions to directed random graphs are straightforward. 

\begin{figure}[t]
\centering 
\includegraphics[width = .325 \linewidth, keepaspectratio]{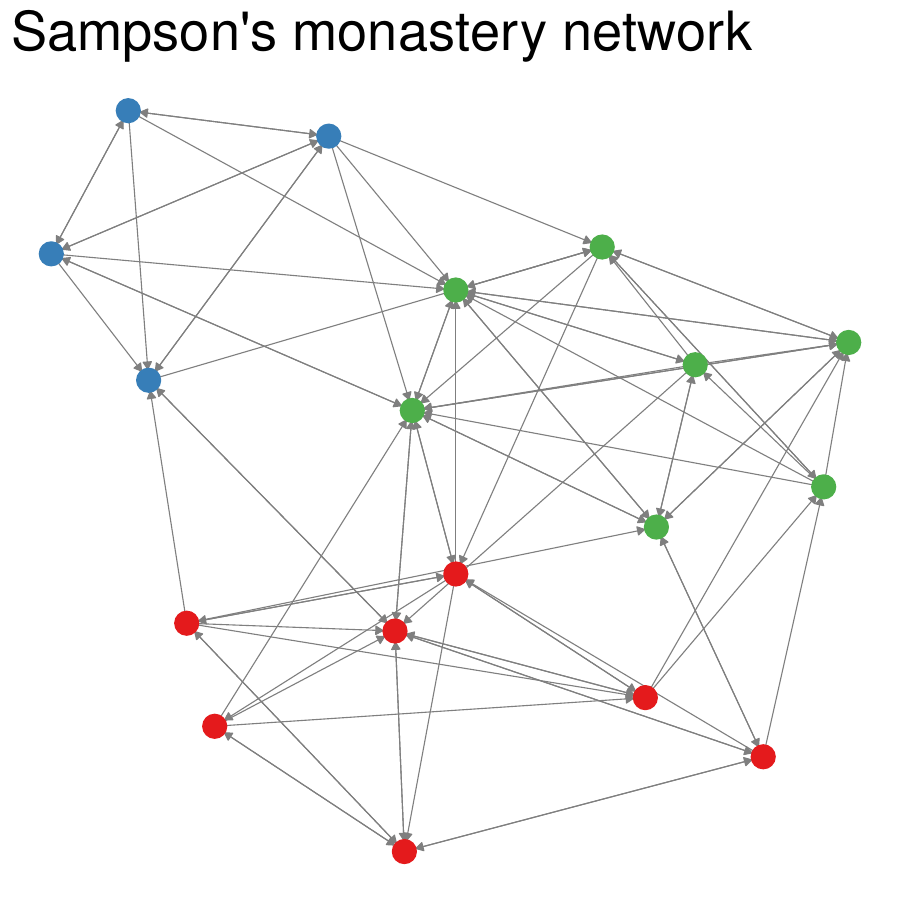} % 
\includegraphics[width = .325 \linewidth, keepaspectratio]{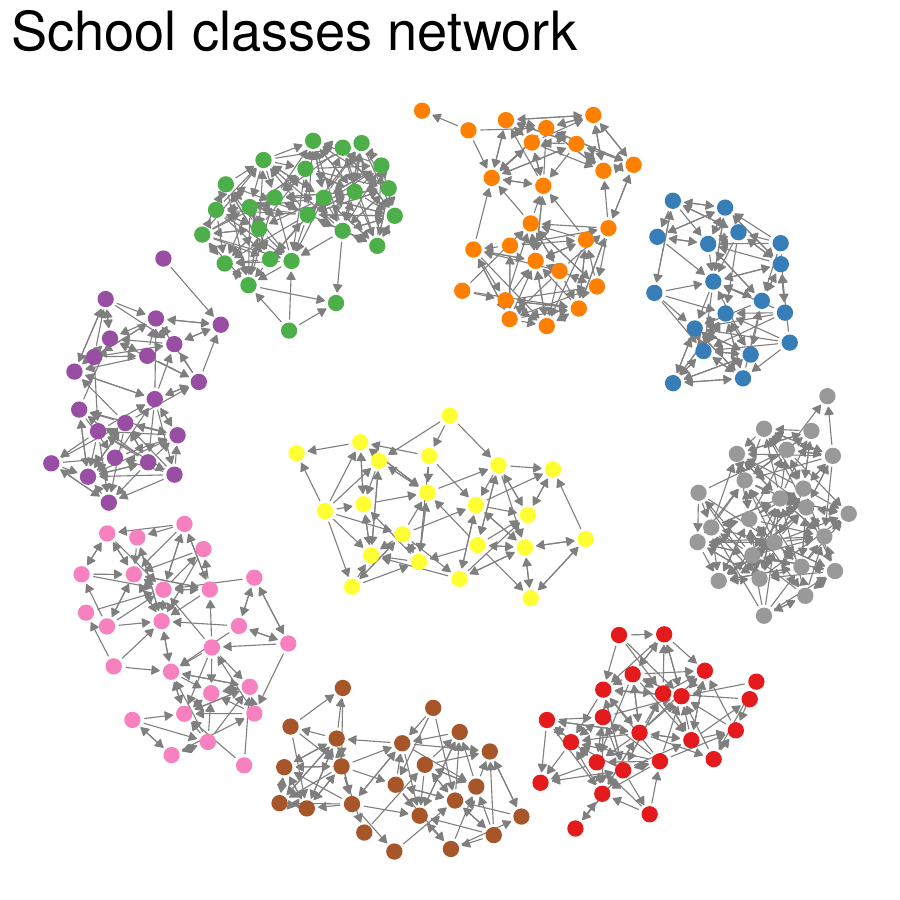} %
\includegraphics[width = .325 \linewidth, keepaspectratio]{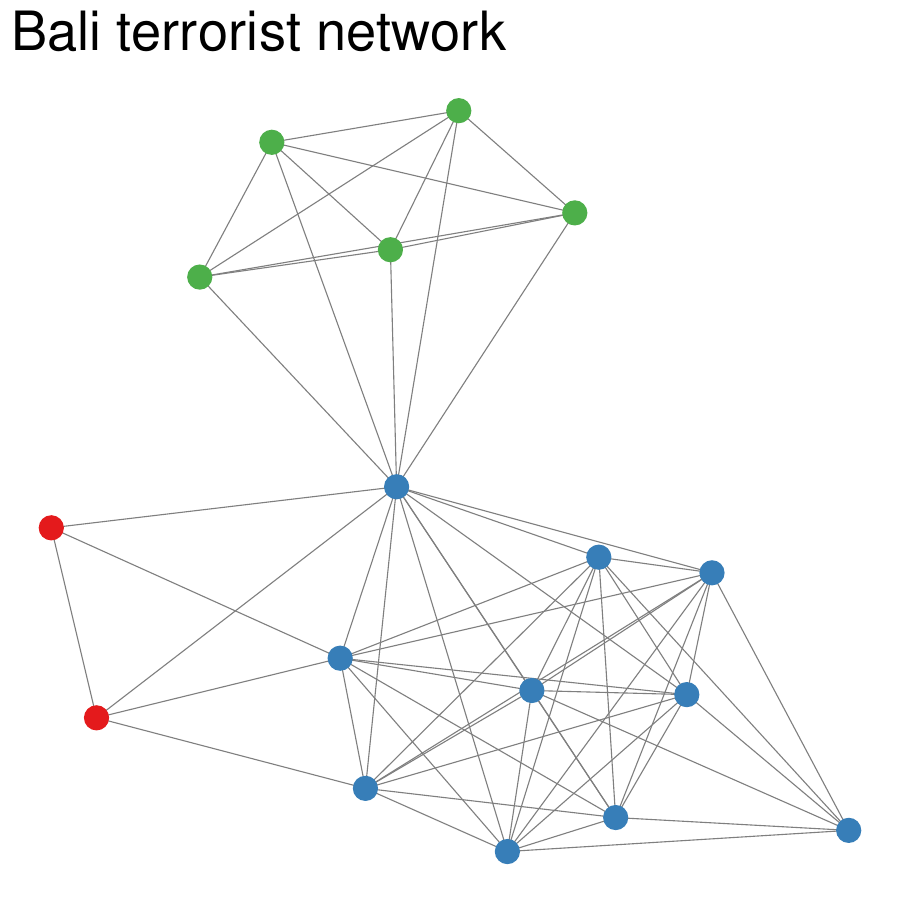}
\caption{\label{fig:nets} Three real data examples of networks 
for which local dependence random graph models would be applicable, 
including Sampson's monastery network, 
the school classes data set from \citet{Stewart2019}, 
and the Bali terrorist network studied in \citet{Schweinberger2015}.
Node colors correspond to block memberships.}
\end{figure}

A {\it local dependence random graph} 
\citep{Schweinberger2015}
is a random graph $\bX$ where the set of nodes $\mN$ 
is partitioned into $K$ {\it blocks} $\mA_1, \ldots, \mA_K$ with probability distributions $\mbP$ 
of the form 
\be
\label{eq:loc_dep}
\mbP(\bX = \bx)
\= \dprod_{1 \leq k \leq l \leq K} \, \mbP_{k,l}(\bX_{k,l} = \bx_{k, l}), 
&& \bx \in \mbX, 
\ee
where the subgraphs $\bX_{k,l}$ ($1 \leq k \leq l \leq K$) are defined as follows: 
\hide{ 
\bi
\item {\it Within-block subgraphs:} 
\beno 
\bX_{k,l} 
&\coloneqq& (X_{i,j} \,:\, i < j, \; i \in \mA_k, j \in \mA_k) 
&\in& \mbX_{k,k} 
&\coloneqq& \{0,1\}^{\tbinom{|\mA_k|}{2}},
&& \mbox{for all } 1 \leq k = l \leq K. 
\ee
\item {\it Between-block subgraphs:} 
\beno
\bX_{k,l} 
&\coloneqq& (X_{i,j} \,:\, i \in \mA_k, \; j \in \mA_l)
&\in& \mbX_{k,k} 
&\coloneqq& \{0,1\}^{|\mA_k| \, |\mA_l|}
&& \mbox{for all } 1 \leq k < l \leq K.
\ee
\ei
}
\beno
\bX_{k,l}
&\coloneqq& \begin{cases}
(X_{i,j})_{\{(i,j) \,:\, i < j, \; i \in \mA_k, j \in \mA_k\}}
\hspace{.1cm}\in\hspace{.1cm} \mbX_{k,k} \hspace{.06cm}\coloneqq\; \{0,1\}^{\tbinom{|\mA_k|}{2}}
& \mbox{if } k = l\s\s\\
(X_{i,j})_{\{(i,j) \,:\, i \in \mA_k, \; j \in \mA_l\}}
\hspace{.67cm}\in\hspace{.1cm} \mbX_{k,l} \hspace{.13cm}\coloneqq\; \{0,1\}^{|\mA_k|\, |\mA_l|} & \mbox{if } k\neq l.\s
\end{cases}
\ee
We refer to the subgraphs $\bX_{k,k}$ ($1 \leq k \leq K$)  
as the {\it within-block subgraphs}
and to the subgraphs $\bX_{k,l}$ ($1 \leq k < l \leq K$) as the {\it between-block subgraphs}.
The probability distribution $\mbP_{k,l}$ is the marginal probability distribution of 
the subgraph 
$\bX_{k,l}$
($ 1 \leq k \leq l \leq K$). 
A {\it local dependence random graph model} is any probability distribution $\mbP$ for $\bX$
of the form \eqref{eq:loc_dep}. 
Figure \ref{fig:nets}
visualizes three networks  which can be studied 
using local dependence random graph models.
%A key aspect of \eqref{eq:loc_dep} inherent to local dependence random graph models 
%is that the block-based subgraphs $\bX_{k,l}$ ($1 \leq k \leq l \leq K$) are independent.
While the collection of block-based subgraphs 
$\bX_{k,l}$ ($1 \leq k \leq l \leq K$) are independent, 
edges within the same block-based subgraph can be dependent. 
The joint distribution $\mbP$ 
can be specified by specifying the marginal probability distributions $\mbP_{k,l}$ 
for the block-based subgraphs 
$\bX_{k,l}$ ($1 \leq k \leq l \leq K$). 

It is worth noting that the block memberships 
are known in both Sampson's monastery network and the school classes network visualized in  
Figure \ref{fig:nets}, 
whereas the block memberships of the Bali terrorist network were estimated as in \citet{Schweinberger2015}. 
When the block memberships correspond to tangible and observable quantities 
(e.g., school class memberships of students), 
data on the block memberships can be collected as part of the observation process. 
When this is not possible, 
the block memberships must be estimated, 
for example by 
using the two-step estimation methodology of \citet{Babkin2020}, 
which estimates both the   
block memberships of nodes and the parameter vectors of local dependence random graph models.

Exponential families account for the most prevalent specifications of local dependence random graph models 
\citep[e.g.,][]{Schweinberger2015, Stewart2019, Dahbura2021, Schubert2022, Tolochko2024},
indeed having been the statistical foundations for the class in the seminal work by \citet{Schweinberger2015}.  
Moreover, 
exponential families provide a flexible statistical platform for constructing models of edge dependence 
in network data applications 
\citep{Lusher2012,Schweinberger2020}, 
and have been shown to be possess desirable statistical properties 
in local dependence random graph models,
including the consistency of maximum likelihood estimators of canonical and curved exponential families 
\citep{SchweinbergerStewart2020}. 
An {\it exponential-family local dependence random graph model} 
can be specified via the marginal probability distributions
of the block-based subgraphs: 
\be
\label{eq:exp_fam_sub}
\mbP_{k,l,\btheta_{k,l}}(\bX_{k,l} = \bx_{k,l})
\= h_{k,l}(\bx_{k,l}) \;
\exp\left( \langle \btheta_{k,l}, \, s_{k,l}(\bx_{k,l}) \rangle - \psi_{k,l}(\btheta_{k,l})\right),
%&& \bv \in \mbX_{k,l}, 
\ee
defined for each $\bx_{k,l} \in \mbX_{k,l}$, 
where 
\bi
\setlength\itemsep{.5em}
\item $s_{k,l} : \mbX_{k,l} \mapsto \mbR^{p_{k,l}}$ is a vector of sufficient statistics; 
\item $\nat_{k,l} \in \mbR^{p_{k,l}}$ is the natural parameter vector; 
\item $h_{k,l} : \mbX_{k,l} \mapsto [0, \infty)$ is the reference function of the exponential family;
and 
\item 
$\psi_{k,l}(\btheta_{k,l}) 
= \log \sum_{\bv \in \mbX_{k,l}} h_{k,l}(\bv) \, \exp(\langle \btheta_{k,l}, \, s_{k,l}(\bv)\rangle)$
is the log-normalizing constant. 
\ei
It is straightforward to show that exponential family specifications 
of the marginal probability distributions of the within-block and between-block subgraphs
will lead to a joint distribution which is also an exponential family.  

A diverse range of models with the local dependence property in \eqref{eq:loc_dep}
can be constructed 
through different specifications of the sufficient statistics and reference functions.
To allow for a general scope of well-structured models, 
we assume that the joint distributions of $\bX$ take the form 
\be
\label{eq:exp_fam}
\mbP_{\nat}(\bX = \bx)
\= \dprod_{1 \leq k \leq l \leq K} \, \mbP_{k,l,\btheta_{k,l}}(\bX_{k,l} = \bx_{k,l})
\= h(\bx) \, \exp\left( \langle \btheta, \, s(\bx) \rangle - \psi(\btheta)\right),
\ee
where $\btheta = (\btheta_W, \btheta_B) \in \mbR^{p+q}$ and  
$s(\bx) = (s_W(\bx_W), s_B(\bx_B)) \in \mbR^{p+q}$, 
with the definitions 
\beno
& \bx_W \coloneqq (\bx_{1,1}, \ldots, \bx_{K,K}),
\quad\quad
\bx_B \coloneqq (\bx_{1,2}, \ldots, \bx_{1,K}, \bx_{2,3}, \bx_{2,4}, \ldots, \bx_{K-1,K}), \s\s\\  
%& s_{W}(\bx) 
%\;\;\coloneqq\;\; \dsum_{k=1}^{K} \, s_{k,k}(\bx_{k,k}), \quad\quad
%s_{B}(\bx)
%\;\coloneqq\; \dsum_{1 \leq k < l \leq K} \, s_{k,l}(\bx_{k,l}), \s\\
& h(\bx)
\;\coloneqq\; \dprod_{1 \leq k \leq l \leq K} h_{k,l}(\bx_{k,l}), \quad\quad\mbox{and}\quad\quad 
\psi(\nat)
\;\;\coloneqq\;\; \dsum_{k=1}^{K} \psi_{k,k}(\btheta_W) \;+ \dsum_{1 \leq k < l \leq K} \psi_{k,l}(\btheta_B).
\ee
Throughout, 
we will assume that $p = \dim(\btheta_W)$ and $q = \dim(\btheta_B)$. 
The exponential family is then the set of probability distributions $\{\mbP_{\nat} : \nat \in \mbR^{p+q}\}$,
where we note that the natural parameter space is equal to $\mbR^{p+q}$,
a fact which follows trivially due to the fact that the support $\mbX$ of $\bX$ is a finite set. 
We additionally assume throughout this work that the exponential family implied by \eqref{eq:exp_fam} is minimal.
While the product of the block-based subgraph distributions in \eqref{eq:exp_fam_sub}
will form an exponential family, 
it may not be minimal,
in which case 
we assume that the representation in \eqref{eq:exp_fam} 
is the minimal representation of the exponential family
obtained through reduction by sufficiency, reparameterization, and proper choice of reference measure;
see Proposition 1.5 of \citet{Brown1986}.
The assumption that an exponential family
is minimal is not restrictive,
as any non-minimal exponential family can be reduced to a minimal exponential family
\citep[Proposition 1.5,][]{Brown1986}.

We next provide examples of exponential-family local dependence random graph models 
in order to motivate the broad scope of this class of models, 
as well as to demonstrate how to construct local dependence random graph models. 
%models in this class. 
As the scope of possible models that can be constructed is large, 
we are unable to present a complete primer on the topic, 
and refer to works by \citet{Schweinberger2015}, \citet{Stewart2019}, and \citet{SchweinbergerStewart2020},
for further information on and concrete examples of exponential-family local dependence random graph models.

\subsection{Examples of exponential-family local dependence random graph models}
\label{sec:examples}

%We present examples of local dependence random graph models in the exponential family framework,
%with select examples being utilized in the simulation studies in Section \ref{sec:simulation}.  

\subsubsection{Example 1: The stochastic block model} 

As a first example, 
we review the stochastic block model \citep{Holland1983},
which is a special case of a local dependence random graph model. 
The joint distribution for $\bX$ is given by  
\be
\label{sb}
\mbP_{\nat}(\bX = \bx)
&\propto& \left[ 
\dprod_{k =1}^{K} \; 
\dprod_{i < j \,:\, i,j \in \mA_k} \, 
\exp(\theta_{k,k} \, x_{i,j})
\right]
\left[ \dprod_{1 \leq k < l \leq K} \, 
\dprod_{(i,j) \in \mA_k \times \mA_l} \, 
\exp(\theta_{k,l} \, x_{i,j})
\right] \s\\
&\propto& \exp\left( \dsum_{k=1}^{K} \, \theta_{k,k} \; \dsum_{i < j \,:\, i,j \in \mA_k} x_{i,j} 
+ \dsum_{1 \leq k < l \leq K} \, \theta_{k,l} \, \dsum_{(i,j) \in \mA_k \times \mA_l} x_{i,j} \right),
\ee
where $\theta_{k,l} \in \mbR$ ($1 \leq k \leq l \leq K$). 
The second line of \eqref{sb} implies the minimal exponential family, 
where each block-based subgraph is a collection of independent and identically distributed Bernoulli random variables 
whose edge probability depends on the subgraph index $(k,l)$ and the value of $\theta_{k,l} \in \mbR$. 

\subsubsection{Example 2: Transitivity in local dependence random graphs} 

The second example we present captures stochastic tendencies towards edge transitivity in networks,
by including a sufficient statistic which models the stochastic tendency for an edge in the network 
to belong to a triangle. 
For this example,
we consider joint distributions $\{\mbP_{\nat} : \nat \in \mbR^3\}$ for $\bX$ of the form  
\beno
\mbP_{\nat}(\bX = \bx) 
&\propto& \exp\big( \theta_1 \, s_1(\bx) + \theta_2 \, s_2(\bx) + \theta_3 \, s_3(\bx) \big),
\ee
with natural parameters $(\theta_1, \theta_2, \theta_3) \in \mbR^3$, 
and where the sufficient statistics are given by 
\beno 
s_1(\bx) \= \dsum_{k=1}^{K} \; \dsum_{i < j \,:\, i,j \in \mA_k} x_{i,j} \s\s\\
s_2(\bx) \= \dsum_{k=1}^{K} \; \dsum_{i < j \,:\, i,j \in \mA_k} \,
x_{i,j} \; \one\left( \dsum_{h \in \mA_k \setminus \{i,j\}} x_{i,h} \, x_{j,h} \,\geq\, 1 \right) \s\s\\
s_3(\bx) \= \dsum_{1 \leq k < l \leq K} \, \dsum_{(i,j) \in \mA_k \times \mA_l} \, x_{i,j}. 
\ee
In words, 
$s_1(\bx)$ counts the number of edges in each of the within-block subgraphs
$\bx_{k,k}$ ($1 \leq k \leq K$), 
whereas  
$s_3(\bx)$ counts the number of edges in each of the between-block subgraphs 
$\bx_{k,l}$ ($1 \leq k < l \leq K$). 
Neither of these statistics induce dependence, 
as when $\theta_2 = 0$,
the joint distribution will factorize with respect to the edge variables
in the graph,
implying edges are independent.  

The second sufficient statistic induces dependence among edge variables contained in the same 
within-block subgraph,
noting that the form of the statistic in $s_2(\bx)$ 
ensures that distributions will not factorize with respect to the edge variables within the graph 
when $\theta_2 \neq 0$. 
The second statistic $s_2(\bx)$ counts the number of edges $x_{i,j}$ 
between pairs of nodes $i$ and $j$ belonging to a common block $\mA_k$,
which are also mutually connected to at least one other node $h \in \mA_k$ also belonging to the same common block,
i.e.,
it counts the number of within-block edges which form at least one triangle within the respective block-based subgraph. 
We call such edges {\it transitive edges}.

Motivation for taking this approach to constructing models of edge dependence 
lies in foundational properties of exponential families. 
%Recall that 
The mean-value parameter map of the exponential family is given by 
$\bmu(\nat) \coloneqq \mbE_{\nat} \, (s_1(\bX), s_2(\bX), s_3(\bX))$ \citep[p. 73--74,][]{Brown1986},
mapping the natural parameter space $\mbR^3$ to the interior of the convex hull of the image of 
$\mbX$ under the vector of sufficient statistics $s : \mbX \mapsto \mbR^3$:  
\beno
\bmu(\nat) 
&\in& \mbM 
&\coloneqq& \interior\left(  \text{ConHull}\left(\left\{ s(\bx) \in \mbR^3 \;:\; \bx \in \mbX \right\}\right)   \right),
\ee
where $\text{ConHull}(\mathscr{S})$ 
represents the convex hull of the set $\mathscr{S}$.  
Moreover, 
for a minimal exponential family (of which this example is), 
the map $\bmu : \mbR^3 \mapsto \mbM$ defines a homeomorphism between $\mbR^3$ and $\mbM$
\citep[Theorem 3.6,][]{Brown1986}. 
This last point emphasizes a key modeling aspect of exponential-family local dependence random graph models, 
as for any point $\bu \in \mbM$ parameterizing the expected values (mean values) of the sufficient statistics 
$(s_1(\bX), s_2(\bX), s_3(\bX))$,
we are guaranteed to be able to find a natural parameter vector $\nat \in \mbR^3$ 
for which $\mbE_{\nat} \, (s_1(\bX), s_2(\bX), s_3(\bX)) = \bu$,
allowing specified models to flexibly capture average tendencies of networks, 
including density, 
transitivity, 
and much more.

\subsubsection{Example 3: Incorporating node and block heterogeneity into models}

The third example shows how we are able to incorporate heterogeneous parameterizations for blocks,
as well as for the stochastic propensities of different nodes to form edges, 
demonstrating ways in which the dimension of parameter vectors can grow in applications.  
For ease of presentation, 
we will build on Example 2 by extending the sufficient statistics which
were specified in that example.  

First, 
we will demonstrate how heterogeneity in node degrees can be incorporated into models. 
Suppose that nodes are divided into $M$ non-overlapping groups or categories $\{1, \ldots, M\}$
which we represent as sets $\mcG_1, \ldots, \mcG_M$. 
Note that these are distinct from the blocks $\mA_1, \ldots, \mA_K$. 
In applications, 
these groups might be ranks in a department, 
gender, 
race, 
or any other categorical covariate which can be observed and treated as fixed. 
As such, 
each block may be comprised of different amount of nodes from each of the groups $\mcG_1, \ldots, \mcG_M$,
an example of which is the school classes data set studied in \citet{Stewart2019},
where each school class was comprised of different numbers of male and female students.  
 
We replace $s_1(\bx)$ in Example 2 by multiple sufficient statistics:  
\beno
s_m(\bx) 
\= \dsum_{k=1}^{K} \; \dsum_{i \in \mcG_m \cap \mA_k} \; \dsum_{j \in \mA_k \setminus \{i\}} \, x_{i,j},
&& m \in \{1, \ldots, M\}, 
\ee
with 
natural parameters $(\theta_1, \ldots, \theta_M) \in \mbR^M$. 
A version of this statistic is implemented in the {\tt R} package {\tt ergm}
under the name {\tt nodefactor} \citep{Krivitsky2023}.
In words, 
the model includes a sufficient statistic that, 
based on the value of the corresponding natural parameter, 
adjusts the baseline propensity for within-block edge formation involving nodes in that group. 
With no other sufficient statistics in the model,  
the log-odds of an edge would be given by 
\beno
\log \, \dfrac{\mbP_{\nat}(X_{i,j} = 1)}{\mbP_{\nat}(X_{i,j} = 0)}
\= \theta_m + \theta_n,
&& i \in \mcG_m\cap \mA_k, \; j \in \mcG_n \cap \mA_k,
&& k \in \{1, \ldots, K\}.   
\ee 
This is reminiscent of the $p1$ model \citep{Holland1981} and the $\beta$-model \citep{Chatterjee2011},
in which 
each node is given its own distinct class.
%i.e., 
%in the $\beta$-model, 
%$M = N$ and $\mcG_i = \{i\}$ ($i \in \{1, \ldots, N\}$). 
 
We now show how heterogeneity can arise in the block-based subgraphs,
by allowing different parameterizations for different blocks. 
The statistic $s_2(\bx)$ in Example 2 counts the number of transitive edges 
in each within-block subgraph $\bX_{k,k}$ ($1 \leq k \leq K$), 
using the same parameter for each within-block subgraph. 
It may be that different blocks display different tendencies towards transitivity. 
To make this concrete, 
suppose that the individual blocks $\{1, \ldots, K\}$ 
are partitioned into $L$ groups or categories 
$\mcH_1, \ldots, \mcH_L$.
We then replace the sufficient statistic $s_2(\bx)$ in Example 2 
by multiple statistics:  
\beno
s_{M+l}(\bx)
\= \dsum_{k \in \mcH_l} \;
\dsum_{i < j \,:\, i,j \in \mA_k} \;
x_{i,j} \; \one\left( \dsum_{h \in \mA_k \setminus \{i,j\}} x_{i,h} \, x_{j,h} \,\geq\, 1 \right), 
&& l \in \{1, \ldots, L\},  
\ee
with natural parameters $(\theta_{M+1}, \ldots, \theta_{M+L}) \in \mbR^L$. 
The complete model is given by  
\beno
\mbP_{\nat}(\bX = \bx) 
&\propto& \exp\left( 
\dsum_{t=1}^{M+L+1} \, \theta_t \, s_t(\bx) 
%+ \dsum_{h=M+1}^{M+L+1} \, \theta_h \, s_h(\bx)
%+ \theta_{M+L+1} \, s_{M+L+1}(\bx) 
\right),
\ee
with natural parameter space $\mbR^{M+L+1}$,
where the last sufficient statistic is equal to  
\beno
s_{M+L+1}(\bx)
\= 
\dsum_{1 \leq k < l \leq K} \, \dsum_{(i,j) \in \mA_k \times \mA_l} x_{i,j}. 
\ee 
Example 3 helps to demonstrate how the number of model parameters can grow quickly in applications when significant
generality,
heterogeneity,
or adaptability is needed to capture important aspects of the application.
A version of Example 3 will be used in the simulation studies conducted in Section \ref{sec:simulation}.

\section{Theoretical guarantees} 
\label{sec:theory}

Our main theoretical results are presented in this section.  
We first review 
exponential family theory for local dependence random graph models   
in Section \ref{sec:exp_fam_prelims}.
Our consistency theory is then presented in Section \ref{sec:main_res_estimators}.
Section \ref{sec:L2} derives rates of convergence 
in the $\ell_2$-norm
of maximum likelihood estimators, 
whereas 
Section \ref{sec:minimax} presents bounds on the minimax risk in the $\ell_2$-norm
which help to establish the minimax optimality (under mild conditions) of the upper bounds presented in Section \ref{sec:L2}. 
Lastly, 
but importantly,  
rates of convergence of the error of the multivariate normal approximation are obtained in  Section \ref{sec:main_res_norm},
providing both non-asymptotic and asymptotic theory for multivariate normal approximations  
of maximum likelihood estimators of local dependence random graph models.  

Due to space restrictions, all proofs are presented in the supplement \citep{Stewart2024-sup}.

\subsection{Preliminaries for exponential families} 
\label{sec:exp_fam_prelims}

The log-likelihood of an exponential-family local dependence random graph model is 
\beno
\ell(\btheta, \bx)
&\coloneqq& \log \, \mbP_{\nat}(\bX = \bx) 
\= \dsum_{k=1}^{K} \ell_{k,k}(\btheta_W, \bx_{k,k}) + 
\dsum_{1 \leq k < l \leq K} \ell_{k,l}(\btheta_B, \bx_{k,l}), 
\ee
where 
\beno
\ell_{k,k}(\btheta_W, \bx_{k,k}) 
&\coloneqq&  \langle \btheta_W, \, s_{k,k}(\bx_{k,k}) \rangle - \psi_{k,k}(\btheta_W) + \log h_{k,k}(\bx_{k,k}) \s \\ 

\ell_{k,l}(\btheta_B, \bx_{k,l})
&\coloneqq& \langle \btheta_B, \, s_{k,l}(\bx_{k,l}) \rangle - \psi_{k,l}(\btheta_B) + \log h_{k,l}(\bx_{k,l}).
\ee
The gradient 
$\nabla_{\btheta} \, \ell(\btheta, \bx) = (\nabla_{\btheta_W} \, \ell(\btheta, \bx), \nabla_{\btheta_B} \, \ell(\btheta, \bx))$
is given by 
\beno
\nabla_{\btheta_W} \, \ell(\btheta, \bx)
\=  \dsum_{k=1}^{K} \, \left[ s_{k,k}(\bx_{k,k}) - \mbE_{k,k,\btheta_W}  \, s_{k,k}(\bX_{k,k}) \right] \s\\
\nabla_{\btheta_B} \, \ell(\btheta, \bx)
\= \dsum_{1 \leq k < l \leq K} \,  \left[ s_{k,l}(\bx_{k,l}) - \mbE_{k,l,\btheta_B} \, s_{k,l}(\bX_{k,l}) \right], 
\ee 
where $\mbE_{k,k,\btheta_W}$ and $\mbE_{k,l,\btheta_B}$ 
are the expectation operators with respect to the marginal probability distributions 
$\mbP_{k,k,\btheta_W}$ of $\bX_{k,k}$
and 
$\mbP_{k,l,\btheta_B}$ of $\bX_{k,l}$,
respectively
\citep[Lemma 6.1,][]{Stewart2024-sup}.  
We denote the set of maximum likelihood estimators for a given observation $\bx \in \mbX$ by 
\beno
\widehat\bTheta
&\equiv& \widehat\bTheta(\bx) 
&\coloneqq& \left\{ \btheta^\prime \in \mbR^{p+q} \,:\, \ell(\btheta^\prime, \bx) = \sup\limits_{\btheta \in \mbR^{p+q}} \, \ell(\btheta, \bx) \right\}.  
\ee
For minimal and regular exponential families, 
the maximum likelihood estimator exists uniquely when it is exists, 
i.e.,
$|\widehat\bTheta| \in \{0, 1\}$ 
\citep[Proposition 3.13,][]{Sundberg2019}. 
Regarding existence, 
the maximum likelihood estimator of natural parameter vectors of minimal exponential families 
exists when the sufficient statistic vector falls within the interior of the mean-value parameter space 
\citep[Theorem 5.5., p. 148,][]{Brown1986},
in which case there exists a parameter vector $\widehat\btheta \in \mbR^{p+q}$ 
for which $\bmu(\widehat\btheta) = s(\bx)$ 
for a given observation $\bx \in \mbX$ of the random graph $\bX$,
defining $\bmu(\nat) \coloneqq \mbE_{\nat} \, s(\bX)$ to be the mean-value parameter map \citep[p. 73--74,][]{Brown1986}.    

In practice, 
computing maximum likelihood estimators is not straightforward, 
as the log-normalizing constants are generally computationally intractable
unless the marginal probability distributions of the block-based subgraphs $\bX_{k,l}$ ($1 \leq k \leq l \leq K$)
are assumed to factorize further to reduce the computational burden of computing the normalizing constants,  
because $\psi_{k,k}(\nat_W)$ involves the summation of $\tbinom{|\mA_k|}{2}$ terms 
and $\psi_{k,l}(\nat_B)$ involves the summation of $|\mA_k| \, |\mA_l|$ terms. 
It becomes infeasible to compute these summations in practice even for modest block sizes. 
The prevailing method for estimating exponential families of random graph models is 
Monte-Carlo maximum likelihood estimation (MCMLE) \citep{Hunter2006}. 
%In fact, 
%this is the default method for estimation in the {\tt R} package {\tt ergm} \citep{Krivitsky2023}. 
The algorithm  outlined in \citet{Hunter2006} applies directly to 
exponential-family local dependence random graph models \citep{mlergm, Stewart2019, SchweinbergerStewart2020}, 
and is used in the simulation studies conducted in Section \ref{sec:simulation}
through the implementation in  
the {\tt R} package {\tt mlergm} \citep{mlergm}.

We summarize the key aspects of MCMLE with exponential families of random graph models outlined in \citet{Hunter2006}. 
%in order to keep the manuscript mostly self-contained. 
The essential idea of MCMLE is to approximate intractable likelihood functions with stochastic approximations 
utilizing Markov Chain Monte Carlo (MCMC) methods. 
The crux of the methodology rests on a simple approximation of normalizing constants via importance sampling. 
To introduce the idea, 
let $\nat_0 \in \mbR^{p+q}$ be a fixed parameter vector in the natural parameter space of an 
exponential-family local dependence random graph model.
We can equivalently find maximum likelihood estimators $\mle$ of $\truth$ by 
\beno
\mle 
\= \argmax\limits_{\nat \in \mbR^{p+q}} \, 
\left[ \ell(\nat, \bx) - \ell(\nat_0,\bx) \right]  
\=  \argmax\limits_{\nat \in \mbR^{p+q}} \,
\left[ \langle \nat - \nat_0, \, s(\bx) \rangle 
- \log\left( \exp\left( \psi(\nat) - \psi(\nat_0)\right)\right) \right].
\ee
In order to solve the above optimization problem, 
we need to be able to approximate the gradient corresponding to the above objective function,
which is given by
\be
\label{gradient}
\nabla_{\nat} \, \left[ \ell(\nat, \bx) - \ell(\nat_0,\bx) \right]
\= s(\bx) - \nabla_{\nat} \, \log\left( \exp\left( \psi(\nat) - \psi(\nat_0)\right)\right). 
\ee
The intractability of the normalizing constants in \eqref{gradient} makes direct computation infeasible, 
as discussed. 

We approximate the term $\exp\left( \psi(\nat) - \psi(\nat_0)\right)$ via 
a change of measure argument: 
\beno
\exp\left( \psi(\nat) - \psi(\nat_0)\right)
\= \exp(-\psi(\nat_0)) \, 
\dsum_{\bx \in \mbX} \, h(\bx) \, \exp(\langle \nat, s(\bx) \rangle) \s\\
\= \exp(-\psi(\nat_0)) \,
\dsum_{\bx \in \mbX} \, h(\bx) \, \exp(\langle \nat, s(\bx) \rangle) \, 
\dfrac{\exp(\langle \nat_0, s(\bx)\rangle)}{\exp(\langle\nat_0, s(\bx)\rangle)} \s\\
\= \mbE_{\nat_0} \exp(\langle \nat - \nat_0, s(\bX) \rangle), 
\ee
where $\mbE_{\nat_0}$ is the expectation operator corresponding to $\mbP_{\nat_0}$.
As a result, 
if we can approximate the expectation $\mbE_{\nat_0} \exp(\langle \nat - \nat_0, s(\bX) \rangle)$ via Monte Carlo methods,
then we can approximate the ratio of normalizing constants 
$\exp\left( \psi(\nat) - \psi(\nat_0)\right)$.
A key advantage of this approach lies in the fact that the expectation is taken with respect to 
a fixed distribution $\mbP_{\nat_0}$. 
In general, 
we will not be able to sample directly from the distributions
and will need 
to rely on MCMC sampling methods \citep[see, e.g., ][]{Snijders2002, Krivitsky2023}.  
Let $\widetilde\bX_1, \ldots, \widetilde\bX_n$ be an MCMC sample from $\mbP_{\nat_0}$. 
Then, returning to \eqref{gradient}, 
we have the following approximation:  
\beno
\nabla_{\nat} \, \left[ \ell(\nat, \bx) - \ell(\nat_0,\bx) \right]
&\approx& s(\bx) - \nabla_{\nat} \, \log\left(\dfrac{1}{n} \, \dsum_{i=1}^{n} \, 
\exp\left(\langle \nat - \nat_0, \, s(\widetilde\bx_i) \rangle\right) \right)  \s\\
\= s(\bx) - \dsum_{i=1}^{n} \, \left( \dfrac{\exp\left(\langle \nat - \nat_0, \, s(\widetilde\bx_i) \rangle\right)}{
\sum_{j=1}^{n} \, \exp\left(\langle \nat - \nat_0, \, s(\widetilde\bx_j) \rangle\right)} \right) \, s(\widetilde\bx_i). 
\ee
Using this approximation,  
root finding algorithms---such as stochastic gradient descent or Fisher scoring algorithms---can 
be utilized to find the MCMLE approximation to the MLE; 
\citet{Hunter2006} outlines a stochastic Fisher scoring algorithm.  
The convergence of the MCMLE to the MLE depends on the convergence 
of the exact log-likelihood to the stochastic approximation 
\citep[see, e.g., discussions in][]{Geyer1992}, 
which will depend upon properties of the Markov chain utilized to generate sample networks. 
In the usual implementations \citep[e.g.,][]{Krivitsky2023},
these chains will be geometrically mixing toward the target sampling distribution 
and will provide good approximations provided sufficient computational resources have been expended.
As a final point on the computational complexity, 
different model specifications, implementations of MCMC methodology,
and block structures will have different mixing times and thus will require differing amounts of computational resources. 
With regards to scalability, 
access to parallel computing presents a significant opportunity to improve computation times 
by exploiting the independence of the block-based subgraphs to parallelize simulation; 
see discussions in \citet{Babkin2020}, 
which analyzed networks with over 10,000 nodes 
utilizing parallel computing and an implementation of the stochastic Fisher scoring algorithm of 
\citet{Hunter2006} in the {\tt R} package  {\tt mlergm} \citep{mlergm}. 

\subsection{Convergence rates of maximum likelihood estimators} 
\label{sec:main_res_estimators}

We derive non-asymptotic bounds on the $\ell_2$-error of maximum likelihood estimators  
which hold with high probability.
Our results extend those of \citet{SchweinbergerStewart2020}, 
who derived consistency results for maximum likelihood estimators of canonical and curved exponential-family 
local dependence random graph models, 
but did not report rates of convergence. 
Additionally, 
\citet{SchweinbergerStewart2020} focused on estimation of only the within-block parameter vectors $\nat_W$. 
In contrast, 
we establish consistency theory with rates of convergence 
for entire parameter vectors $(\nat_W, \nat_B)$ of exponential-family local dependence random graph models,
covering settings  
where the number of model parameters 
%$p+q$ 
and sizes of blocks 
%$|\mA_1|, \ldots, |\mA_K|$
may tend to infinity, 
at appropriate rates. 
The consistency theory in this work is related to---but distinct from---the results in 
\citet{Stewart2020},
who prove a general theorem for establishing consistency 
and rates of convergence 
of maximum likelihood and pseudolikelihood-based estimators
of random graph models with dependent edges with respect to the $\ell_{\infty}$-norm
under a more general weak dependence assumption.  
First, 
we focus specifically on local dependence random graph models 
and quantify rates of convergence 
in the $\ell_2$-norm for this class of models
and in terms of interpretable quantities related to local dependence random graphs,
namely  
properties of the block structure,
graph, 
and model.  
Second, 
our method of proof is fundamentally different from that of both \citet{SchweinbergerStewart2020} and \citet{Stewart2020}, 
and the consistency theory in this work cannot be proved as a corollary to an existing result.
\hide{
As a final point of contrast, 
we 
%leverage properties of local dependence random graph models 
establish 
the first non-asymptotic bound on the error of the multivariate normal approximation 
of a standardization of maximum likelihood estimators,
which in turn provides the first rigorous derivation of a statistical inference procedure 
for estimators of local dependence random graph models. 
Our normality theory is presented in Section \ref{sec:main_res_norm}.  
}

We outline some notational definitions and regularity assumptions for our theorems to follow, 
subsequently discussing each in turn.  
Let $\mB_2(\bv, r) \coloneqq \{\bv^\prime \in \mbR^{\dim(\bv)} : \norm{\bv^\prime - \bv}_2 < r\}$
be the open $\ell_2$-ball 
with center $\bv$ and radius $r > 0$
and denote by $\lambda_{\min}(\bA)$ 
and $\lambda_{\max}(\bA)$ the smallest and largest eigenvalues,
respectively, of the matrix $\bA \in \mbR^{d \times d}$.
We write $a_N = O(b_N)$ when there exists
a constant $C > 0$ and integer $N_0 \geq 1$ 
such that $a_N \leq C \, b_N$ for all $N \geq N_0$,
and write $a_N = o(b_N)$ when there exists, 
for all $\delta > 0$, 
an integer $N_0(\delta) \geq 1$ 
such that $a_N \leq \delta \, b_N$ for all $N \geq N_0(\delta)$. 

\s

\begin{assumption}
\label{a1}
Assume there exist 
$C_W > 0$ and $C_B > 0$,
independent of $N$, $p$, and $q$,
such that
\beno
\sup\limits_{\bx_{k,k} \in \mbX_{k,k}} \,
&\norm{s_{k,k}(\bx_{k,k})}_{\infty}
&\leq& C_W \; \dbinom{|\mA_k|}{2},
&& k \in \{1, \ldots, K\}, \s\s\\

\sup\limits_{\bx_{k,l} \in \mbX_{k,l}} \,
&\norm{s_{k,l}(\bx_{k,l})}_{\infty}
&\leq& C_B \; |\mA_k| \, |\mA_l|,
&& \{k,l\} \subseteq \{1,\ldots,K\}.
\ee
\end{assumption}

\vspace{-.6cm}

\begin{assumption}
\label{a2}
Assume there exists $\epsilon > 0$,
independent of $N$, $p$, and $q$,
such that
\beno
\minw &\coloneqq& \inf\limits_{\nat \in \mB_2(\truth, \epsilon)} \,
\dfrac{\lambda_{\min}\left(- \mbE \, \nabla_{\nat_W}^2 \, \ell(\btheta, \bX)\right)}{K}
&>& 0 \s\\
\minb &\coloneqq&
\inf\limits_{\nat \in \mB_2(\btheta^\star, \epsilon)} \,
\dfrac{\lambda_{\min}\left(-\mbE \, \nabla_{\nat_B}^2 \, \ell(\btheta, \bX)\right)}{\tbinom{K}{2}}
&>& 0. \\
\ee
\end{assumption}

\vspace{-.6cm}

\begin{assumption}
\label{a3}
Define 
$\Aavg \coloneqq K^{-1} \, \sum_{k=1}^{K} \, |\mA_k|$ to be the average block size and
\beno
\maxw 
&\coloneqq& \dfrac{\lambda_{\max}\left(- \mbE \, \nabla_{\nat_W}^2 \, \ell(\truth, \bX)\right)}{K}
&&\mbox{and}&&
\maxb 
&\coloneqq& \dfrac{ \lambda_{\max}\left(- \mbE \, \nabla_{\nat_B}^2 \, \ell(\truth, \bX) \right)}{\tbinom{K}{2}},
\ee
and assume that 
\beno
\sqrt{\Aavg} \; \dfrac{\sqrt{\maxw}}{\minw}
\= o\left( 
\sqrt{\dfrac{N}{p}} \; \right)
&&\mbox{and}&&
\Aavg \; 
\dfrac{\sqrt{\maxb}}{\minb}
\= o\left( 
\sqrt{\dfrac{N^2}{q}} \; \right). 
\ee
\end{assumption}

\vspace{-.5cm}

\begin{assumption}
\label{a4}
Assume the largest block size  
$A_{\max} \coloneqq \max\{|\mA_1|, \ldots, |\mA_K|\}$ 
satisfies 
\beno
A_{\max}
&\leq& \min\left\{
\left( \dfrac{N \, \maxw}{\Aavg \, p^2} \right)^{1/4}, \quad
\left( \dfrac{N^2 \, \maxb}{4 \, \Aavg^2 \, q^2} \right)^{1/4}
\,\right\}.
\ee
\end{assumption}

\vspace{-.4cm}

\begin{remark}[Discussion of Assumption \ref{a1}]
%With regards to Assumption (A.1),
We place a restriction on the scaling of the block-based
sufficient statistic vectors with respect to the sizes of the blocks.
The  need for this arises out of a need to derive concentration inequalities for gradients of the log-likelihood,
as well as a need to control third-order derivatives of the log-likelihood function
in our method of proof for deriving bounds on the error of the multivariate normal approximation.
The assumption is natural,
as it essentially requires that the values of the sufficient statistics possess an upper-bound which is proportional
to the number of edge variables in each of the respective block-based subgraphs.
An example of interest is the transitive edge count statistic of a within-block subgraph $\bX_{k,k}$,
discussed also in in Section \ref{sec:examples},
given by
\beno
\dsum_{i < j \,:\, i,j \in \mA_k} \; 
x_{i,j} \; \one\left( \, \dsum_{h \in \mA_k \setminus \{i,j\}} \,
x_{i,h} \, x_{j,h} \,\geq\, 1 \right)
&\leq& \dsum_{i < j \,:\, i,j \in \mA_k} \, x_{i,j}
&\leq& \dbinom{|\mA_k|}{2},
\ee
which can be viewed as a special case of the geometrically-weighted
edgewise shared partner statistic \citep{Hunter2006, Stewart2019}.
To further contextualize this assumption,
it is helpful to note that Assumption \ref{a1}
is related to the issue of instability of exponential-families of random graph models
\citep{Schweinberger2011}.
Maximal changes in the sufficient statistic vectors $s_{k,k}(\bx)$ $(1 \leq k \leq K$)
and $s_{k,l}(\bx)$ ($1 \leq k < l \leq K$)
due to changing the value of a single edge in $\bx$ are defining characteristics of
instability in exponential-family random graph models,
in the sense of \citet{Schweinberger2011}.
Assumption \ref{a1} implies limitations on the sensitivity of the sufficient statistic vectors to changes in the edges in the graph.
Understanding this connection
helps to explain why local dependence random graph models achieve statistical behavior and properties
not achieved in early---but flawed---statistical models of edge dependence in network data
\citep{Haggstrom99,Jonasson1999,Schweinberger2011,Chatterjee2013}.
Lastly,
it is worth noting that Assumption \ref{a1} could be relaxed further,
allowing for a larger upper bound.
The result of this,
however,
would be looser upper-bounds on the $\ell_2$-error
and slower rates of convergence.  
\end{remark}

\begin{remark}[Discussion of Assumption \ref{a2}] 
Assumption \ref{a2} places a restriction on the scaling of the smallest eigenvalue of the joint Fisher information matrix
by placing an assumption on the scaling of the smallest eigenvalue of the Fisher information matrices
$-\mbE \, \nabla_{\nat_W}^2 \, \ell(\nat, \bX)$
and
$ -\mbE \, \nabla_{\nat_B}^2 \, \ell(\nat, \bX)$
corresponding to the within-block and between-block probability distributions,
respectively,
in a neighborhood $\mB_2(\truth, \epsilon)$ of the data-generating parameter vector $\truth = (\truth_W, \truth_B)$. 
The local dependence assumption and the assumption that the parameter vector $\nat = (\nat_W, \nat_B) \in \mbR^{p+q}$
partitions the within-block and between-block parameters 
implies that 
the joint Fisher information matrix $-\mbE \, \nabla_{\nat}^2 \, \ell(\nat, \bX)$
has the form
\beno
-\mbE \, \nabla_{\nat}^2 \, \ell(\nat, \bX)
\= \left(
\begin{matrix}
-\mbE \, \nabla_{\nat_W}^2 \, \ell(\nat, \bX) & \bm{0}_{p,q} \s\\
\bm{0}_{q,p} & -\mbE \, \nabla_{\nat_B}^2 \, \ell(\nat, \bX)
\end{matrix}
\right),
\ee
where $\bm{0}_{m,n}$ is the $(m\times n)$-dimensional matrix of all zeros.
%Hence, 
%Assumption \ref{a2} can be viewed as a minimal information criterion which requires that 
%we obtain sufficient information about the parameter vector $(\nat_W, \nat_B)$ from
%an observation of the random graph. 
Assumption \ref{a2} essentially assumes that the Fisher information matrices are invertible 
in a neighborhood of the data-generating parameter vector.  
Minimum eigenvalue restrictions of Fisher information matrices
are standard in settings where the number of model parameters may tend to infinity
\citep[e.g.,][]{Portnoy1988,Ravikumar2010,Jankova2018}.
Notably,
our assumption represents a restriction on what amounts to an average minimum eigenvalue
(averaged over the block-based quantities in both the within-block and between-block cases).  
To understand why we have adopted this definition in our assumptions (relevant also to Assumption \ref{a3}),
instead of placing a restriction on the minimum eigenvalue of the Fisher information matrices corresponding
to each block-based subgraph,
observe through Weyl's inequality,
the bound 
\beno
\lambda_{\min}\left(-\mbE \, \nabla_{\nat_W}^2 \, \ell(\nat, \bX)\right)
&\geq& \dsum_{k=1}^{K} \, \lambda_{\min}\left(-\mbE \, \nabla_{\nat_W}^2 \, \ell_{k,k}(\nat_W, \bX_{k,k})\right) \s\\
&\geq& K \, \left(\min\limits_{k \in \{1, \ldots, K\}} \,
\lambda_{\min}\left(-\mbE \, \nabla_{\nat_W}^2 \, \ell_{k,k}(\nat_W, \bX_{k,k})\right) \right).
\ee
If certain subgraphs do not contain any information about
certain subsets of parameters, 
possibly due to heterogeneous parameterizations that allow different blocks to have different parameters, 
then it may be the case that $\lambda_{\min}(-\mbE \, \nabla_{\nat_W}^2 \, \ell_{k,k}(\nat_W, \bX_{k,k})) = 0$
for  some $k \in \{1, \ldots, K\}$ due to singularity.
As a result and in order to cover more general settings and heterogeneous parameterizations,
we place our minimum eigenvalue restriction on the scaling of the averaged  smallest eigenvalue of 
joint Fisher information matrices.  
\end{remark}

\begin{remark}[Discussion of Assumption \ref{a3}]
Assumption \ref{a3} places a regularity assumption on three key quantities, 
the average block size $\Aavg \coloneqq K^{-1} \sum_{k=1}^{K} |\mA_k|$,
the average minimum eigenvalues of Fisher information matrices $\minw$ and $\minb$ 
in a neighborhood of the data-generating parameter vector 
$\truth = (\truth_W, \truth_B)$ (defined in Assumption \ref{a2}),
and the average maximum eigenvalues of Fisher information matrices $\maxw$ and $\maxb$ 
at the data-generating parameter vector.
As will be seen in Theorem \ref{thm:L2},
Assumption \ref{a3} essentially outlines a scaling requirement of these three quantities (in their respective cases)
which ensures consistent estimation under Theorem \ref{thm:L2}, 
in the sense that the upper bounds on the $\ell_2$-error in Theorem \ref{thm:L2} will tend to zero 
as the size of the network $N$ tends to infinity. 
As such, 
Assumption \ref{a3} can be viewed as a minimal information criterion which requires that
we obtain sufficient information about the parameter vector $(\truth_W, \truth_B)$ from
an observation of the random graph $\bX$. 
\end{remark}

\begin{remark}[Discussion of Assumption \ref{a4}]
In our method of deriving concentration inequalities, 
we bound factors involving the influence of edge variables in the random graph 
by the size of the largest block size, 
noting that dependence is restricted to block-based subgraphs whose size is dominated by functions of the largest 
block size. 
Similar approaches have been taken in \citet{SchweinbergerStewart2020}. 
Notably, 
Assumption \ref{a4} does not assume that the sizes of blocks are fixed
and allows these quantities to grow without bound. 
However, 
this assumption places a restriction on how large  blocks can be in order to ensure that the derived 
concentration inequalities are sufficiently sharp to facilitate the development of the statistical theory of this work. 
\end{remark}

\subsubsection{Upper bounds on the $\ell_2$-error of maximum likelihood estimators} 
\label{sec:L2}

The first theoretical result we present establishes 
upper-bounds on the $\ell_2$-error of maximum likelihood estimators 
for exponential-family local dependence random graph models which hold with high probability,
presented in Theorem \ref{thm:L2}. 
This paves the way for 
establishing bounds on rates of convergence of maximum likelihood estimators with respect to the $\ell_2$-norm. 
We will address the question of optimal rates of convergence in Section \ref{sec:minimax},
where we outline a set of sufficient conditions for which we prove the upper bounds in Theorem \ref{thm:L2}
are minimax optimal,
in the sense that the upper bounds derived in Theorem \ref{thm:L2} 
match (up to an unknown constant) the minimax rate of convergence.  

\newpage

\begin{theorem}
\label{thm:L2}
Consider a minimal exponential-family local dependence random graph model satisfying  
Assumptions \ref{a1}, \ref{a2}, \ref{a3}, and \ref{a4}
and assume that $p = \dim(\truth_W) \geq \log \, N$ 
and $q = \dim(\truth_B) \geq \log \, N$.  
Then
there exist constants $C > 0$ and $N_0 \geq 3$,
independent of $N$, $p$, and $q$, 
such that,
with probability at least $1 - N^{-2}$, 
the maximum likelihood estimator $\mle = (\mle_W, \mle_B) \in \mbR^{p+q}$
exists, 
is unique,
and satisfies %the following $\ell_2$-error bounds: 
\beno
\norm{\mle_W - \truth_W}_2
&\leq& C \; \sqrt{A_{\avg}} \; \dfrac{\sqrt{\maxw}}{\minw} \, \sqrt{\dfrac{p}{N}}
\ee
\vspace{.05cm}
\beno
\norm{\mle_B - \truth_B}_2
&\leq& C \;  A_{\avg} \; \dfrac{\sqrt{\maxb}}{ \minb} \, \sqrt{\dfrac{q}{N^2}}, 
\ee
for all integers $N \geq N_0$. 
\end{theorem}

\s

Theorem \ref{thm:L2} provides the foundation for establishing 
convergence rates in the $\ell_2$-norm of 
maximum likelihood estimators of exponential-family local dependence random graph models.
The assumption that the exponential family is minimal ensures uniqueness of the maximum likelihood estimator 
when it exists 
\citep[Proposition 3.13,][]{Sundberg2019}. 
%(Proposition 3.13, \citep{Sundberg2019}),  
%The assumption that an exponential family 
%is minimal is not restrictive, 
%as any non-minimal exponential family can be reduced to a minimal exponential family 
%\citep[Proposition 1.5,][]{Brown1986}.
%(Proposition 1.5, \citep{Brown1986}). 
%Theorem \ref{thm:L2} establishes that rates of convergence will depend on
Rates of convergence will depend on 
\bi
\item  the dimensions of the parameters vectors $p = \dim(\truth_W)$ and $q = \dim(\truth_B)$; \s 
\item the ratios $\sqrt{\maxw} \;/\; \minw$
and $\sqrt{\maxb} \;/\; \minb$; and  \s
\item the average block size $\Aavg$,
\ei
with rates of convergence depending on the scaling of these quantities with respect to $N$. 
Theorem \ref{thm:L2} additionally provides a set of sufficient conditions 
for the event that the maximum likelihood estimator exists to occur with high probability. 
Related to discussions in Section \ref{sec:examples}, 
the maximum likelihood estimator exists in the event $s(\bX) \in \mbM$,
recalling the definition of $\mbM$ from Section \ref{sec:examples} 
as the mean-value parameter space of the exponential family. 
The assumptions of Theorem \ref{thm:L2} ensure that the probability of the event $s(\bX) \in \mbM$
occurs with high probability, 
provided the network size $N$ is sufficiently large.
This event essentially requires that the sufficient statistic not fall on the boundary 
of the convex hull of the image of $\mbX$ under the vector of sufficient statistics $s : \mbX \mapsto \mbR^{p+q}$,
i.e.,
$\partial \, \mbM$. 
The probabilities for any fixed network size $N$,
however, 
will depend on both properties of the network and the model specification. 
With regards to the latter, 
significant heterogeneity,
such as in Example 3 in Section \ref{sec:examples}, 
can result in a higher-dimensional parameter space and therefore sufficient statistic vector, 
which can increase the chance of the sufficient statistic vector falling on the boundary $\partial \, \mbM$,
in which event the maximum likelihood estimator will not exist. 

We permit both $\maxw$ and $\maxb$ 
to scale faster than $\minw$ and $\minb$,
respectively, 
provided consistency is still established 
(i.e., provided  Assumption \ref{a3} is met). 
Within the context of exponential families of growing dimension in classical settings 
of a random sample of independent and identically distributed random vectors, 
\citet{Portnoy1988} and \citet{Ghosal2000} obtain similar convergence rates,
in their respective settings.  
Notably,
Theorem 2.1 of \citet{Portnoy1988} arrives at a similar scaling requirement for the minimum and maximum 
eigenvalues of Fisher information matrices.  
A key difference is that both works  
place third order assumptions on the models (see the assumptions of Theorem 2.1 of \citet{Portnoy1988},
and Theorem 2.1 of \citet{Ghosal2000}). 
We avoid the need for such assumptions through the method of proof of Theorem \ref{thm:L2},
but require a smoothness condition on minimum eigenvalues of Fisher information matrices,  
as Assumptions \ref{a2} and \ref{a3} 
restrict the scaling of maximum eigenvalues of the Fisher information matrix 
at the data-generating parameter vector $\truth$ 
relative to minimum eigenvalues of the same within 
 a neighborhood $\mB_2(\truth, \epsilon)$ of $\truth$. 
If we assume additional regularity in the spectrum of the Fisher information matrices 
by assuming that  
\beno
\minwt 
&\coloneqq& \dfrac{\lambda_{\min}\left(- \mbE \, \nabla_{\nat_W}^2 \, \ell(\truth, \bX)\right)}{K}
&=& O\left(\minw\right) \s\s\\
\minbt
&\coloneqq&
\dfrac{\lambda_{\min}\left(-\mbE \, \nabla_{\nat_B}^2 \, \ell(\truth, \bX)\right)}{\binom{K}{2}}
&=& O\left(\minb\right),
\ee
then we could prove a corollary to Theorem \ref{thm:L2} which establishes the upper bounds 
\beno
\norm{\mle_W - \truth_W}_2
&\leq& C \; \sqrt{A_{\avg}} \; \dfrac{\sqrt{\maxw}}{\minwt} \, \sqrt{\dfrac{p}{N}}
\ee
\beno
\norm{\mle_B - \truth_B}_2
&\leq& C \;  A_{\avg} \; \dfrac{\sqrt{\maxb}}{ \minbt} \, \sqrt{\dfrac{q}{N^2}},
\ee 
which are more analogous to the results of \citet{Portnoy1988}. 
Related to other works within the statistical network analysis, 
our consistency results and rates of convergence have key connections to theoretical results for the $\beta$-model, 
for example those obtained in \citet{Shao2021}, 
which includes convergence rates for parameters of the $\beta$-model in the $\ell_2$-norm, 
and also the $\ell_{\infty}$- and $\ell_1$-norms.

As a final point, 
Theorem \ref{thm:L2} assumes that the block memberships are known,
i.e.,
the blocks $\mA_1, \ldots, \mA_K$ 
are observed or estimated without error. 
In many cases,
the block memberships can be observed through the observation process 
\citep[e.g.,][]{Stewart2019, SchweinbergerStewart2020}. 
However, 
in certain settings this may not be possible and the block memberships must be estimated
\citep[e.g.,][]{Babkin2020,Schweinberger2020-Bernoulli}. 
In both cases, 
the results of Theorem \ref{thm:L2} can be regarded as the estimation error of an oracle estimate
with  perfect knowledge or estimation of the block structure of the network. 
The impact of imperfect block membership knowledge on theoretical guarantees (whether through a noisy observation 
or error in the estimation of block memberships of nodes) is an open question for future research.

\subsubsection{Minimax risk in the $\ell_2$-norm and optimal rates of convergence} 
\label{sec:minimax}

We next turn to the question of whether the upper bounds on the $\ell_2$-error 
established in Theorem \ref{thm:L2} are optimal, 
in the sense that they match (up to an unknown constant) 
the rates of convergence of the minimax risk in the $\ell_2$-norm. 

We define the minimax risk with respect to the $\ell_2$-norm to be  
\be
\label{eq:minimax_risks}
\mcR_{W,N}
&\coloneqq& 
\inf\limits_{\mle_W} \;  
\sup\limits_{\nat \in \mbR^{p+q}} \; 
\mbE_{\nat} \, \norm{\mle_W - \nat_W}_2 \s\\
\mcR_{B,N}
&\coloneqq& \inf\limits_{\mle_B} \; 
\sup\limits_{\nat \in \mbR^{p+q}} \; 
\mbE_{\nat} \, \norm{\mle_B - \nat_B}_2. 
\ee
The method by which we establish lower bounds to the 
minimax risk in the $\ell_2$-norm 
requires placing an assumptions on the average value of the largest eigenvalues of Fisher information matrices, 
similar to the roles of $\maxw$ and $\maxb$ in Theorem \ref{thm:L2},
extended now to a neighborhood $\mB_2(\truth, \epsilon)$ of $\truth$. 
Fix $\epsilon > 0$,
independent of $N$, $p$, and $q$,
and define
\be
\label{eq:max_e}
\maxwe &\coloneqq& \sup\limits_{\nat \in \mB_2(\truth, \epsilon)} \;
\dfrac{\lambda_{\max}\left(- \mbE \, \nabla_{\nat_W}^2 \, \ell(\nat, \bX)\right)}{K} \s\\
\maxbe &\coloneqq& \sup\limits_{\nat \in \mB_2(\truth, \epsilon)} \;
\dfrac{\lambda_{\max}\left(- \mbE \, \nabla_{\nat_B}^2 \, \ell(\nat, \bX)\right)}{\tbinom{K}{2}}.
\ee
We first establish lower bounds to the minimax risks $\mcR_{W,N}$ and $\mcR_{B,N}$ in Theorem \ref{thm:optimal},
which enable us to outline 
sufficient conditions for 
the upper bounds on the $\ell_2$-error presented in Theorem \ref{thm:L2} to achieve (up to an unknown constant)
the minimax rates of convergence; 
see Corollary \ref{cor:optimal}. 
In the following results,  
it is helpful to recall that  
$\minw$ and $\minb$ are defined in Assumption \ref{a2},
$\maxw$ and $\maxb$ are defined in Assumption \ref{a3}, and 
$\maxwe$ and $\maxbe$ are defined in \eqref{eq:max_e}.

\begin{theorem}
\label{thm:optimal}
(Lower bound to the minimax risk) 
Consider an exponential-family local dependence random graph model
satisfying Assumption \ref{a2}.  
Then there exist constants $C_1 > 0$ and $C_2 > 0$, 
independent of $N$, $p$, and $q$, 
such that the minimax risks $\mcR_{W,N}$ and $\mcR_{B,N}$
defined in \eqref{eq:minimax_risks} satisfy 
%have the following lower bounds: 
\beno
\mcR_{W,N}
&\geq& C_1  \sqrt{\dfrac{\Aavg}{\maxwe}} \; \sqrt{\dfrac{p}{N}} 
&\geq& C_1  \left(\dfrac{\minw}{\maxwe} \right) 
\, \dfrac{\sqrt{\maxw}}{\minw} \; \sqrt{\Aavg} \; \sqrt{\dfrac{p}{N}}
\ee
\vspace{.05cm}
\beno
\mcR_{B,N}
&\geq& C_2  \dfrac{\Aavg}{\sqrt{\maxbe}} \; \sqrt{\dfrac{q}{N^2}}
&\geq& C_2  \left(\dfrac{\minb}{\maxbe} \right) 
\, \dfrac{\sqrt{\maxb}}{\minb} \; \Aavg \; \sqrt{\dfrac{q}{N^2}}, 
\ee
provided $p = \dim(\truth_W) = O(N \, \maxwe)$ 
and $q = \dim(\truth_B) = O(N^2 \, \maxbe)$.  
\end{theorem}

\s

The role of Assumption \ref{a2} in Theorem \ref{thm:optimal} is to ensure that both  
$\minw$ and $\minb$ are bounded away from $0$,
ensuring all of the lower bounds are well defined,
whereas we assume 
\be
\label{pq_cond}
p 
\;=\; \dim(\truth_W) 
\;=\; O(N \, \maxwe)
&&\mbox{and}&& 
q \;=\; \dim(\truth_B) 
\;=\; O(N^2 \, \maxbe)
\ee
in order to satisfy a technical condition in the proof of Theorem \ref{thm:optimal}. 
Under the assumption that the maximum eigenvalues of Fisher information matrices are bounded away from $0$,
the condition in \eqref{pq_cond}
requires that $p = O(N)$ and $q = O(N^2)$,
which places a much less stringent restriction on the dimensions of parameters vectors 
when compared with Assumption \ref{a3}.  
Two sets of lower bounds are presented in Theorem \ref{thm:optimal}, 
with the first being the most sharp, 
but unhelpful in our pursuit of studying whether the rates of convergence implied 
in Theorem \ref{thm:L2} are minimax optimal. 
The second, though looser, set of bounds approximately match the upper bounds 
on the $\ell_2$-error established in Theorem \ref{thm:L2}.
Indeed, 
this second set of bounds allows us to establish conditions for such minimax optimality in Corollary \ref{cor:optimal}
which is presented below.  

Note that the lower bound to the minimax risk presented in Theorem \ref{thm:optimal} 
considers $\truth \in \mbR^{p+q}$. 
The fact that the parameter space is unbounded introduces no complications when deriving lower bounds; 
however, 
when  turning to the problem of deriving an upper bound to the minimax risk, 
an unbounded parameter space presents new challenges.
The following theorem obtains upper bounds on the minimax risk with respect to the $\ell_2$-norm 
in a neighborhood of the data-generating parameter vector. 
%The upper bounds to the minimax risk help to establish the scaling of the minimax risk with respect to the $\ell_2$-norm 
%in a neighborhood $\mB_2(\truth, \epsilon)$ of a pre-specified 
%parameter vector $\truth$. 
%This allows us to confirm the scaling of the minimax risk with respect to key quantities at least locally 
%in the neighborhood of a data-generating parameter vector. 
%where $\epsilon > 0$ 
%is assumed the be the same as in Assumptions \ref{a2} and \ref{a3} and in \eqref{eq:max_e}. 

\begin{theorem}
\label{thm:upper}
(Upper bound to the minimax risk) 
Under the assumptions of Theorem \ref{thm:L2}, 
there exist constants $C_1 > 0$, $C_2 > 0$, and $N_0 \geq 3$,  
independent of $N$, $p$, and $q$, 
such that the minimax risks restricted to a local neighborhood 
$\mB_2(\truth, \epsilon)$ of a point $\truth \in \mbR^{p+q}$ 
satisfy,
for all integers $N \geq N_0$, 
\beno
\inf\limits_{\mle_W} \;
\sup\limits_{\nat \in \mB_{2}(\truth, \epsilon)} \;
\mbE_{\nat} \, \norm{\mle_W - \nat_W}_2
&\leq& C_1 \, \dfrac{\sqrt{\maxw}}{\minw} \; \sqrt{\Aavg} \; \sqrt{\dfrac{p}{N}}
\ee
\vspace{.05cm}
\beno
\inf\limits_{\mle_B} \;
\sup\limits_{\nat \in \mB_2(\truth, \epsilon)} \;
\mbE_{\nat} \, \norm{\mle_B - \nat_B}_2
&\leq& C_2 
\, \dfrac{\sqrt{\maxb}}{\minb} \; \Aavg \; \sqrt{\dfrac{q}{N^2}}, 
\ee
where $\epsilon > 0$ is the same as in Assumptions \ref{a2} and \ref{a3} and in \eqref{eq:max_e}. 
\end{theorem}

The final result of this section is concerned with outlining a set of sufficient conditions 
which allow us to establish the minimax optimality of Theorem \ref{thm:L2}.

\begin{corollary}
\label{cor:optimal}
%Consider an exponential-family local dependence random graph model
Under the assumptions of Theorems \ref{thm:L2} and \ref{thm:optimal},
and the assumption that  
\be
\label{eq:opt_cond}
\maxwe \;=\; O\left( \minw \right)
&&\mbox{and}&& 
\maxbe \;=\; O\left( \minb \right), 
\ee
the maximum likelihood estimators $\mle_W$ and $\mle_B$ 
achieve the minimax rate of convergence,
in the sense that the upper bounds on the $\ell_2$-error of $\mle_W$ and $\mle_B$
presented in Theorem \ref{thm:L2} match (up to an unknown constant which is independent of $N$, $p$, and $q$)
the lower bounds to the minimax risks in Theorem \ref{thm:optimal}.  
\end{corollary}

If the exponential-family local dependence random graph model satisfies 
\eqref{eq:opt_cond},
then Corollary \ref{cor:optimal} establishes the minimax optimality of the rates of convergence 
for maximum likelihood estimators implied via Theorem \ref{thm:L2}. 
Such an assumption is common in the high-dimensional statistics literature
\citep[e.g.,][]{Ravikumar2010,Jankova2018},
where it is common to assume that minimum and maximum eigenvalues of Fisher information matrices corresponding 
to the sampling distribution are bounded away from $0$ and from above, respectively. 
We can interpret condition \eqref{eq:opt_cond} similarly, 
however applied to the joint Fisher information for the entire collection of random variables in the random graph 
(in contrast to the sampling distribution from which a random sample is generated) 
and in a neighborhood $\mB_2(\truth, \epsilon)$ of the data-generating parameter vector $\truth$.

\subsection{Convergence rates of the multivariate normal approximation}
\label{sec:main_res_norm}

A key challenge to any statistical analysis of network data 
is finding rigorous justification for statistical inference methodology. 
The main contributing factor to this challenge lies in the fact that statistical analyses of network data
are typically in the setting of a single collection of dependent random variables 
without the benefit of replication. 
In other words, 
any statistical inference will be based on a single observation of a collection of dependent binary random variables. 
It is common for inference of model parameters in exponential-family random graph models 
to utilize the normal approximation for carrying out inference about estimated coefficients 
\citep[e.g.,][]{Krivitsky2023,Lusher2012,Stewart2019}. 
Except in select cases, 
these inferences are performed without rigorous theoretical justification, 
owing to the difficulty of obtaining theoretical results 
establishing the validity of 
the normal approximation 
in scenarios with a single observations of a collection of dependent binary random variables.

The dependence structure of local dependence random graph models facilitates 
proof of rigorous theoretical results justifying the normal approximation for estimators,
and in this section, 
we obtain rates of convergence of the multivariate normal approximation in scenarios of increasing model dimension. 
It is worth noting that our results imply the univariate normal approximation,
as multiple univariate tests are frequently utilized in applications
\citep[e.g.,][]{Stewart2019}.  
%In practice, 
%multivariate tests are not utilized as frequently as univariate tests, 
%as the most common application of the normal approximation in 
%exponential-family random graph models is to perform significance tests for elements of the estimated parameter vector. 
%Nonetheless, 
%the multivariate normal approximation establishes the univariate normal approximation. 
Similarly to our consistency results presented in Theorem \ref{thm:L2},
the quality of the multivariate normal approximation will depend on key quantities 
related to the block structure, graph, and model specification. 

Throughout,
$\bZ_d$ will denote a $d$-dimensional multivariate normal random vector 
with mean vector $\bm{0}_d$ (the $d$-dimensional vector of all zeros) 
and covariance matrix $\bI_{d}$ (the $d$-dimensional identity matrix).  
The probability distribution of $\bZ_d$ is denoted by $\Phi_d$. 

In order to establish our multivariate normal approximation theory, 
we leverage a multivariate Berry-Esseen theorem provided in \citet{Raic2019}, 
together with a Taylor expansion of the log-likelihood equation. 
Utilizing properties of exponential families, 
we are able to derive non-asymptotic bounds on the error of the multivariate normal approximation 
for a standardization of the maximum likelihood estimator,
providing the first results which elaborate conditions under which  the normal approximation is expected 
to produce valid inferences in local dependence random graph models.

\begin{theorem}
\label{thm:normal}
Consider a minimal exponential-family local dependence random graph model satisfying
Assumptions \ref{a1}, \ref{a2}, \ref{a3}, and \ref{a4}
and assume that $p = \dim(\truth_W) \geq \log \, N$  and 
$q = \dim(\truth_B) \geq \log \, N$. 
Then there exist constants $C_1 > 0$, $C_2 > 0$, and $N_0 \geq 3$,
independent of $N$, $p$, and $q$,
and a random vector $\bDelta \in \mbR^{p+q}$ 
such that,
for all integers $N \geq N_0$ and measurable convex sets $\mcC \subset \mbR^{p+q}$,
\beno
&& \left|\mbP(\mcI(\truth)^{1/2} \, (\mle - \truth) + \bDelta \,\in\,\mcC) - \Phi_d(\bZ_d \in \mcC)\right| \s\s\\
&\leq& C_1 \, (p+q)^{1/4} \, A_{\max}^{7} \,
\left[ \sqrt{\dfrac{p^3}{(\minwt)^{3} \, N}}
+ \sqrt{\dfrac{q^3}{(\minbt)^{3} \, N^2}} \; \right], 
\ee
where the random vector $\bDelta$ satisfies 
\beno
\mbP\left(\norm{\bDelta}_2
\,\leq\,
C_2 \,  A_{\max}^{6}
\sqrt{\Aavg \, \dfrac{(\maxw)^2}{(\minw)^5} \, \dfrac{p^5}{N}
+ \Aavg^2 \, \dfrac{(\maxb)^2}{(\minb)^5} \, \dfrac{q^5}{N^2}}
\;\, \right)
&\geq& 1 - \dfrac{1}{N^2}.
\ee
\end{theorem}

The standardization 
$\mcI(\truth)^{1/2}\,  (\mle - \truth)$
is of a familiar form in multivariate normal approximation settings.
The quantity $\bDelta$ can be interpreted as an error term or a random perturbation, 
arising due to a Taylor approximation.
While our result is stated 
for $\mcI(\truth)^{1/2}\,  (\mle - \truth) + \bDelta$,
an important aspect of Theorem \ref{thm:normal} lies in establishing that 
the random perturbation $\bDelta$ to $\mcI(\truth)^{1/2}\,  (\mle - \truth)$ 
is small (in the $\ell_2$-norm) with high probability,
justifying basing inferences and derivations of confidence regions on 
$\mcI(\truth)^{1/2}\,  (\mle - \truth)$
in applications.  
Indeed, 
under mild assumptions (which we state below), 
it is straightforward to establish that $\norm{\bDelta}_2$ converges almost surely to $0$ 
as $N \to \infty$.  

A remark is in order regarding the term $(p+q)^{1/4}$ in the upper bound on the error of the multivariate normal approximation
in Theorem \ref{thm:normal}. 
Current results on multivariate Berry-Esseen bounds involve terms which are functions of the dimension of the random vector
\citep{Raic2019}. 
Here, 
the total dimension of the random vector is $p+q$,
as we are proving the joint multivariate normality of a standardization of 
the entire vector of maximum likelihood estimators $(\mle_W, \mle_B)$
which has dimension $p+q$. 
In other words,
we are unable to separate the error into two terms which are functions of only 
quantities based on within-block and between-block quantities, 
as was done in our consistency theory in Section \ref{sec:main_res_estimators}.

Typically, 
both $\mcI_W(\truth_W)^{1/2}$ 
and $\mcI_B(\truth_B)^{1/2}$ will be unknown,
but can be approximated in practice.  
We can approximate both $\mcI_W(\truth_W)$
and $\mcI_B(\truth_B)$ 
through Monte-Carlo methods, 
as Fisher information matrices of canonical exponential families are the covariance matrices of the sufficient statistics. 
This is a common approach to estimating the Fisher information matrix
in the exponential-family random graph model literature,
owing to the fact that models are frequently estimated via Monte-Carlo maximum likelihood estimation,
which already requires simulating sufficient statistic vectors 
\citep[e.g.,][]{Hunter2006, Krivitsky2023},
and discussed in Section \ref{sec:exp_fam_prelims}.

Under an additional regularity assumption, 
we can simplify the bounds presented in Theorem \ref{thm:normal}.
%We state this additional assumption and then demonstrate how to simplify our bounds on the 
%error of the multivariate normal approximation. 
%We end the section with a corollary that elaborates an asymptotic normality result.  

\s

\begin{assumption}
\label{a5}
Assume that there exist constants $L > 0$ and $U > 0$ such that 
\be
\label{reg_norm}
0 &<& L
&\leq& \min\left\{\minw, \; \minb \right\}
&\leq& \max\left\{\maxw, \; \maxb \right\}
&\leq& U, 
\ee
for all values of $N$, $p$, and $q$. 
\end{assumption} 

\vspace{-.05cm}

Assumption \ref{a5} is reminiscent of minimum and maximum eigenvalue restrictions in the high-dimensional statistics literature,
where it is common to assume the minimum and maximum eigenvalues of Fisher information matrices are bounded 
away from $0$ and from above,
respectively  
\citep[e.g.,][]{Ravikumar2010,Jankova2018}. 
Assumption \ref{a5} can be interpreted similarly,  
though applied to the averaged minimum and maximum eigenvalues of the joint Fisher information matrices; 
see also the discussions following Corollary \ref{cor:optimal}.  

Under Assumptions \ref{a1}, \ref{a2}, \ref{a3}, \ref{a4}, and \ref{a5},
we may leverage Theorem \ref{thm:normal} to establish, 
for all measurable convex sets $\mcC \subset \mbR^{p+q}$,
the new bound of  
\beno
|\mbP(\mcI(\truth)^{1/2} \, (\mle - \truth) + \bDelta \,\in\,\mcC) - \Phi_d(\bZ_d \in \mcC)|
&\leq& C_1  (p+q)^{1/4}  \, A_{\max}^{7} 
\left[ \sqrt{\dfrac{p^3}{N}}
+ \sqrt{\dfrac{q^3}{N^2}} \, \right],
\ee
where $\bDelta$ now satisfies
\beno
\mbP\left(\norm{\bDelta}_2
\,\leq\,
C_2 \,  A_{\max}^{6}
\sqrt{\Aavg \,
\dfrac{p^5}{N}
+ \Aavg^2 \,
\dfrac{q^5}{N^2}}
\; \right)
&\geq& 1 - \dfrac{1}{N^2}.
\ee
In certain settings, 
it may be the case that properties of the network limit the sizes of the blocks, 
in which the size of the largest block $A_{\max}$ may be bounded for all network sizes. 
Under the additional assumption that the sizes of the blocks are bounded above, 
%and that $\Aavg$ possess a positive lower bound, 
we can absorb the quantities involving $A_{\max}$ and $\Aavg$ into the constants $C_1 > 0$ 
and $C_2 > 0$ in the above bounds. 
This results in the following simple bounds on the error of the multivariate normal approximation: 
\beno
\left|\mbP(\mcI(\truth)^{1/2} \, (\mle - \truth) + \bDelta \,\in\,\mcC) - \Phi_d(\bZ_d \in \mcC)\right|
&\leq& C_1 \, (p+q)^{1/4} \, 
\left[ \sqrt{\dfrac{p^3}{N}}
+ \sqrt{\dfrac{q^3}{N^2}} \; \right], 
\ee
where $\bDelta$ will then satisfy 
\beno
\mbP\left(\norm{\bDelta}_2
\,\leq\,
C_2 \,  
\sqrt{
\dfrac{p^5}{N}
+ 
\dfrac{q^5}{N^2}}
\; \right)
&\geq& 1 - \dfrac{1}{N^2},
\ee
for all measurable convex sets $\mcC \subset \mbR^{p+q}$. 
Note,
in the above results, 
that the probability bounds approach $1$ sufficiently fast,
allowing us to establish,  
through the Borel–Cantelli lemma,
that $\norm{\bDelta}_2$ converges 
$\mbP$-almost surely to $0$ as $N \to \infty$,
provided the upper bounds on $\norm{\bDelta}_2$ tend to $0$ as $N \to \infty$. 

Finally, 
to deliver a simple and easily interpretable result for statistical inference,  
we prove a corollary to Theorem \ref{thm:normal} establishing the asymptotic multivariate normality 
of maximum likelihood estimators.

\begin{corollary}
\label{cor:normal2}
Under the assumptions of Theorem \ref{thm:normal}, 
Assumption \ref{a5}, 
and assuming  
\beno
\lim\limits_{N \to \infty} \;
\max\left\{
A_{\max}^{6}
\sqrt{\Aavg \,
\dfrac{p^5}{N}
+ \Aavg^2 \,
\dfrac{q^5}{N^2}}, \;\;
(p+q)^{1/4}  \, A_{\max}^{7}
\left[ \sqrt{\dfrac{p^3}{N}}
+ \sqrt{\dfrac{q^3}{N^2}} \, \right]
\right\} 
\= 0,
\ee
we have the distributional limit $\mcI(\truth)^{1/2} \, (\mle - \truth) \overset{D}{\to} \bZ_{p+q}$ 
as $N \to \infty$. 
\end{corollary}

\s

Corollary \ref{cor:normal2} can be proved directly by observing that the assumptions of the corollary 
ensure the error bounds in Theorem \ref{thm:normal} converge to $0$ in the limit as $N \to \infty$. 
As a result of Corollary \ref{cor:normal2}, 
standard procedures for constructing confidence regions, 
univariate confidence intervals, 
and 
performing statistical hypothesis tests for significance of parameters 
are justified using the asymptotic approximation 
of the variance-covariance matrix $\mcI(\truth) = \var \, s(\bX)$,
which we discuss above. 
When the sizes of the blocks are bounded as above, 
the essential condition for asymptotic multivariate normality becomes  
\beno
\lim\limits_{N \to \infty} \, 
\sqrt{
\dfrac{p^5}{N}
+
\dfrac{q^5}{N^2}}
\= 0,
\ee
restricting the maximum growth with $N$ 
of the dimensions of the parameters vectors
$p = \dim(\nat_W)$ and $q = \dim(\nat_B)$,  
suggesting  
that both $p = \dim(\nat_W) = o(N^{1/5})$
and $q = \dim(\nat_B) = o(N^{2/5})$  
must hold in our theory for the error of the multivariate normal approximation to vanish in the limit as $N \to \infty$.  

Up to now, 
we required knowledge of the Fisher information matrix $\mcI(\truth) = \var \, s(\bX)$. 
We end the section with a result concerning the estimation of this term for practical implementation.  
We define  
\beno
\widetilde{\mcI}_W^\star
&\coloneqq& \dfrac{\mbE[-\nabla_{\btheta_W}^2 \, \ell(\truth, \bX)]}{K}
&&&\mbox{and}&&&
\widetilde{\mcI}_B^\star
&\coloneqq& \dfrac{\mbE[-\nabla_{\btheta_B}^2 \, \ell(\truth, \bX)]}{\tbinom{K}{2}}.
\ee
Natural estimators for each are given by 
\beno
\widehat{\mcI}_W
&\coloneqq& \dfrac{1}{K} \, \dsum_{k=1}^{K} \, 
\left(s_{k,k}(\bX_{k,k}) - \bar{s}_W(\bX_W) \right) \, 
\left(s_{k,k}(\bX_{k,k}) - \bar{s}_W(\bX_W) \right)^{\top} \s\s\\
\widehat{\mcI}_B
&\coloneqq& \dfrac{1}{\tbinom{K}{2}} \, \dsum_{1 \leq k < l \leq K} \,
\left(s_{k,l}(\bX_{k,l}) - \bar{s}_B(\bX_B) \right) \, 
\left(s_{k,l}(\bX_{k,l}) - \bar{s}_B(\bX_B) \right)^{\top}
\ee
with the definition 
\beno
\bar{s}_W(\bX_W) &\coloneqq& \dfrac{1}{K} \, \dsum_{k=1}^{K} \, s_{k,k}(\bX_{k,k}) 
&&&\mbox{and}&&&
\bar{s}_B(\bX_B) &\coloneqq& \dfrac{1}{\tbinom{K}{2}} \, \dsum_{1 \leq k < l \leq K} \, s_{k,l}(\bX_{k,l}). 
\ee
The following theorem establishes bounds on the error $\mnorm{\widehat{\mcI}_W - \widetilde{\mcI}_W^\star}_2$ 
and 
$\mnorm{\widehat{\mcI}_B - \widetilde{\mcI}_B^\star}_2$
which hold with high probability,
where $\mnorm{\,\cdot\,}_2$ denotes the spectral matrix norm.  

\begin{theorem}
\label{thm:I_est}
Under the assumptions of Theorem \ref{thm:normal},
the events 
\beno
\mnorm{\widehat{\mcI}_W - \widetilde{\mcI}_W^\star}_2
&\leq& C \, A_{\max}^{2} \, \sqrt{\Aavg \, \maxw} \;\; \left( 
\,  \sqrt{\dfrac{p \, \log(p)}{N}}
+ \sqrt{\dfrac{p^2}{N}} \right) \s\\ 
\mnorm{\widehat{\mcI}_B - \widetilde{\mcI}_B^\star}_2
&\leq& C \,  A_{\max}^{2} \, \Aavg \, \sqrt{\maxb} \;\; \left( 
\,  \sqrt{\dfrac{q \, \log(q)}{N^2}}
+ \sqrt{\dfrac{q^2}{N^2}} \right)
\ee
jointly occur with probability at least $1 - 4 \, N^{-2}$. 
\end{theorem}

The conclusions of Theorem \ref{thm:I_est} 
reiterate the conclusions of our previous theoretical results, 
that if certain quantities related to the sizes of blocks and properties of models 
through the spectral properties of Fisher information matrices are sufficiently well-behaved,
and the dimensions of the parameter vectors do not grow too quickly with $N$, 
then accurate estimation and valid inferences of parameter vectors 
of local dependence random graph models will be obtained with high probability.

\section{Simulation results} 
\label{sec:simulation}

%We conduct simulation studies to study how the realized empirical performance
%of models matches the theoretical results of Section \ref{sec:theory}.  

\subsection{Simulation study 1: Convergence rates of maximum likelihood estimators} 

Simulation study 1 demonstrates that the rate of growth of the dimension of parameter vectors
plays a key role in the finite sample performance. 
We consider three cases in a setting which controls certain aspects of the graph. 
Throughout this study, 
we assume that the sizes of the blocks are all fixed at $50$, 
i.e., 
$|\mA_k| = 50$ for all $k \in \{1, \ldots, K\}$. 
In order to vary the size of the network $N$, 
we vary the number of blocks $K \in \{1, 5, 10, 15, 20\}$,
which results in networks of size $N \in \{50, 250, 500, 750, 1000\}$. 
We focus on a special case of Example 3 from Section \ref{sec:examples},
by assuming that each node $i \in \mN$ is assigned to a group $\mG_1, \ldots, \mG_M$ ($M \geq 2$). 
The specific form of this model is then given by 
\beno
\mbP_{\nat}(\bX = \bx)
&\propto& \exp\left( 
\dsum_{m=1}^{M} \, \theta_m \, s_m(\bx)
+ \theta_{m+1} \, s_{m+1}(\bx)
\right),
\ee
where 
\beno
s_{m}(\bx) 
\= \dsum_{k=1}^{K} \, \dsum_{i \in \mA_k \cap \mG_m} \, \dsum_{j \in \mA_k \setminus \{i\}} \, x_{i,j},
&& m \in \{1, \ldots, M\}, 
\ee
and 
\beno
s_{m+1}(\bx)
\= \dsum_{k=1}^{K} \, \dsum_{i < j \,:\, i \in \mA_k, j \in \mA_k} \, x_{i,j} \, 
\one\left( \dsum_{h \in \mA_k \setminus \{i,j\}} \, x_{i,h} \, x_{j,h} \,\geq\, 1 \right). 
\ee
For this simulation study we will focus on the within-block parameter vector
in order to easily compare the trade-off between the dimension of the parameter vector $p$ 
and the size of the network $N$. 
We can then assume that $X_{i,j} = 0$ with probability one for all $\{i,j\} \subset \mN$ 
belonging to distinct blocks,
i.e.,
the between-block subgraphs $\bX_{k,l}$ ($1 \leq k < l \leq K$) are empty subgraphs with probability one.

\begin{figure}[t]
\centering 
\includegraphics[width = .70 \linewidth]{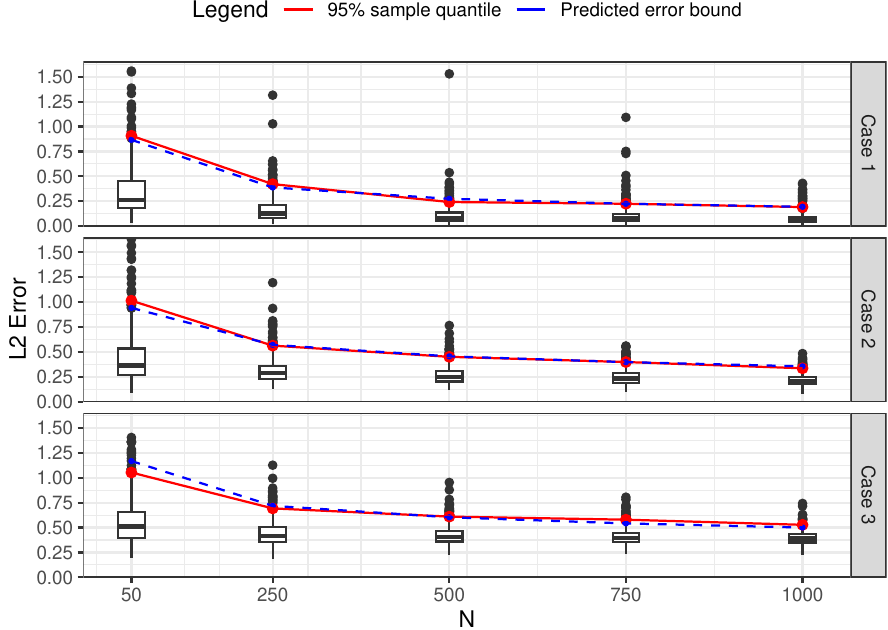}
\caption{\label{fig:sim1} 
The results of Simulation study 1,
which demonstrates the trade-off in finite sample performance of maximum likelihood estimators 
based on the number of model parameters and size of the network. 
Each boxplot for each combination of case and network size is based on $500$ replications. 
Boxplots display the empirical distribution of the $\ell_2$-error, 
whereas the red lines track the $95\%$ sample quantiles 
and the blue dashed lines track the error bounds predicted by Theorem \ref{thm:L2}.  
}
\end{figure}

We consider three cases:
\begin{itemize}[topsep=4pt,itemsep=3.5pt]
\item Case 1: $M = 3$, 
in which case $p = 4$ for all $N \in \{50, 250, 500, 750, 1000\}$.  
\item Case 2: $M = \lceil N^{2/5} \rceil$,
in which case $p \in \{6, 11, 14, 16, 17\}$ depending on the size of the network. 
\item Case 3: $M = \lceil \sqrt{N} \rceil$, 
in which case $ p \in \{9,17, 24, 29, 33\}$ depending on the size of the network.  
\ei
For each case and network size $N \in \{50, 250, 500, 750, 1000\}$, 
we simulate $500$ networks from $\mbP_{\nat}$ 
where $\theta_{M+1} = .5$ 
and $(\theta_1, \ldots, \theta_M) \overset{iid}{\sim} \text{Unif}(-1.5, -.5)$. 
The value of $\theta_{M+1}$ ensures there is a reasonably strong tendency towards transitivity in the network,
and the values of $(\theta_1, \ldots, \theta_M)$ result in networks with plausible densities.  
The results of the Simulation study 1 are summarized in Figure \ref{fig:sim1}. 

The finite sample performance of this study suggests, 
as would be expected based on the results of Theorem \ref{thm:L2},
that the rate at which the $\ell_2$-error converges to $0$ 
is fastest in Case 1 for which the model dimension is fixed,
and slowest in Case 3 for which the model dimension is on the order of $\sqrt{N}$. 
We compute a predicted error bound based on Theorem \ref{thm:L2} by estimating the constant terms, 
which in this simulation study include the average block sizes $\Aavg$ and the largest block size $A_{\max}$,
as well as the terms quantifying averaged eigenvalues of the Fisher information matrices. 
This can be accomplished by estimating constants  
for each network size by 
\beno
\widehat{C}_N
&\coloneqq& Q_{N,.95} \;\Big/\; \sqrt{\dfrac{p}{N}},
&&&& N \in \{50, 250, 500, 750, 1000\}, 
\ee
where $Q_{N,.95}$ is the $95\%$ sample quantile of the $\ell_2$-errors of the maximum likelihood estimators 
based on the $500$ replications, 
and then using the estimate 
\beno
\widehat{C} &\coloneqq& \dfrac{1}{5} \, \dsum_{N \in \{50, 250, 500, 750, 1000\}} \, \widehat{C}_N
\ee
to obtain an overall estimate of the constant term.  
The predicted error bounds are then defined as 
\beno
\widetilde{E}_N &\coloneqq& \widehat{C} \; \sqrt{\dfrac{p}{N}},
&&&& N \in \{50, 250, 500, 750, 1000\}.  
\ee 
The dashed blue lines track the values of $\widetilde{E}_N$ in Figure \ref{fig:sim1},
whereas the red lines track $Q_{N,.95}$.  

Notably, 
the predicted error bound closely matches the $95\%$ sample quantile of the simulated $\ell_2$-errors. 
Theorem \ref{thm:L2} establishes a bound which should hold with high probability,
provided $N$ is sufficiently large. 
Figure \ref{fig:sim1} demonstrates that the predicted error bounds most closely match the realized $95\%$ sample 
quantile of the simulated $\ell_2$-errors for larger network sizes. 
It is also worth noting that an additional source of variation here may be due to the fact that 
the constant term is not actually constant in the network size, 
as the quantities $\maxw$ and $\minw$ may depend on $N$.
With that said, 
though,
the simulation reveals close agreement with the predicted error bounds.

\subsection{Simulation study 2: Error of the normal approximation} 

The second simulation study we conduct explores the error of the normal approximation,
leveraging results in Theorem \ref{thm:normal}. 
We consider the same probability distribution as in Simulation study 1,
in the following two cases:
\begin{itemize}[topsep=4pt,itemsep=3.5pt]
\item Case 1: Fixed parameter dimension $p = 5$ with $M = 4$
categories of each node group and networks of size $N \in \{250, 500, 750, 1000\}$.  
\item Case 2: Growing parameter dimension $p = 2 \, K$ 
with $M = 2 \, K - 1$ categories of each node group,
where there are $50$ nodes per block and the number of blocks 
vary over $K \in \{5, 10, 15, 20\}$,
resulting in networks of size $N \in \{250, 500, 750, 1000\}$.  
\ei
We generate $500$ replications in each case,
simulating networks from the same probability distributions as in Simulation study 1
and in the same manner.  

We study the quality of the normal approximation by 
constructing confidence intervals for the transitive edge parameter 
and Quantile-Quantile plots for the standardized maximum likelihood estimator of the transitive edge parameter.
Our results demonstrate the empirical Type I error in the former matches the theoretical Type I error, 
with the Quantile-Quantile plots not revealing significant departure from normality.  
For each case, 
we constructed $95\%$ confidence intervals and computed the empirical Type I error control. 
Letting $\theta_{m+1}$ and $\widehat{\theta}_{m+1}$ denote the transitive edge parameter 
and the maximum likelihood estimator of the transitive edge parameter, 
we leverage Theorem \ref{thm:normal} to construct confidence intervals: 
\beno
\mbP\left( \theta_{m+1}^{\star} \in 
\left[ \widehat{\theta}_{m+1} - q_{1-\alpha/2} \; \sqrt{[\bS^{-1}]_{m+1,m+1}}, \;\;
\widehat{\theta}_{m+1} + q_{1-\alpha/2} \; \sqrt{[\bS^{-1}]_{m+1,m+1}} \right] \right)
\,\approx\, 1 - \alpha,
& \alpha \in (0, 1), 
\ee
where 
$q_{1-\alpha/2}$ denotes the $(1-\alpha/2)\%$-quantile of the univariate standard normal distribution and  
$\bS$ denotes the sample variance-covariance matrix obtained by sampling sufficient statistics through MCMC methods;
see the discussions in Section \ref{sec:exp_fam_prelims}. 
%estimating the standard error of $\widehat{\theta}_{m+1}$ via $\sqrt{[\bS^{-1}]_{m+1,m_1}}$.  
For Case 1, 
the empirical coverage was 
$(.96, .95, .95, .96)$ corresponding to network sizes of $(250, 500, 750, 1000)$,
and for Case 2,
the same was  
$(.96, .95, .95, .96)$ corresponding to network sizes of $(250, 500, 750, 1000)$. 
The Quantile-Quantile plots for each case across the different network sizes are presented in Figure \ref{fig:qq}. 

\begin{figure}
\centering 
\includegraphics[width = .70\linewidth]{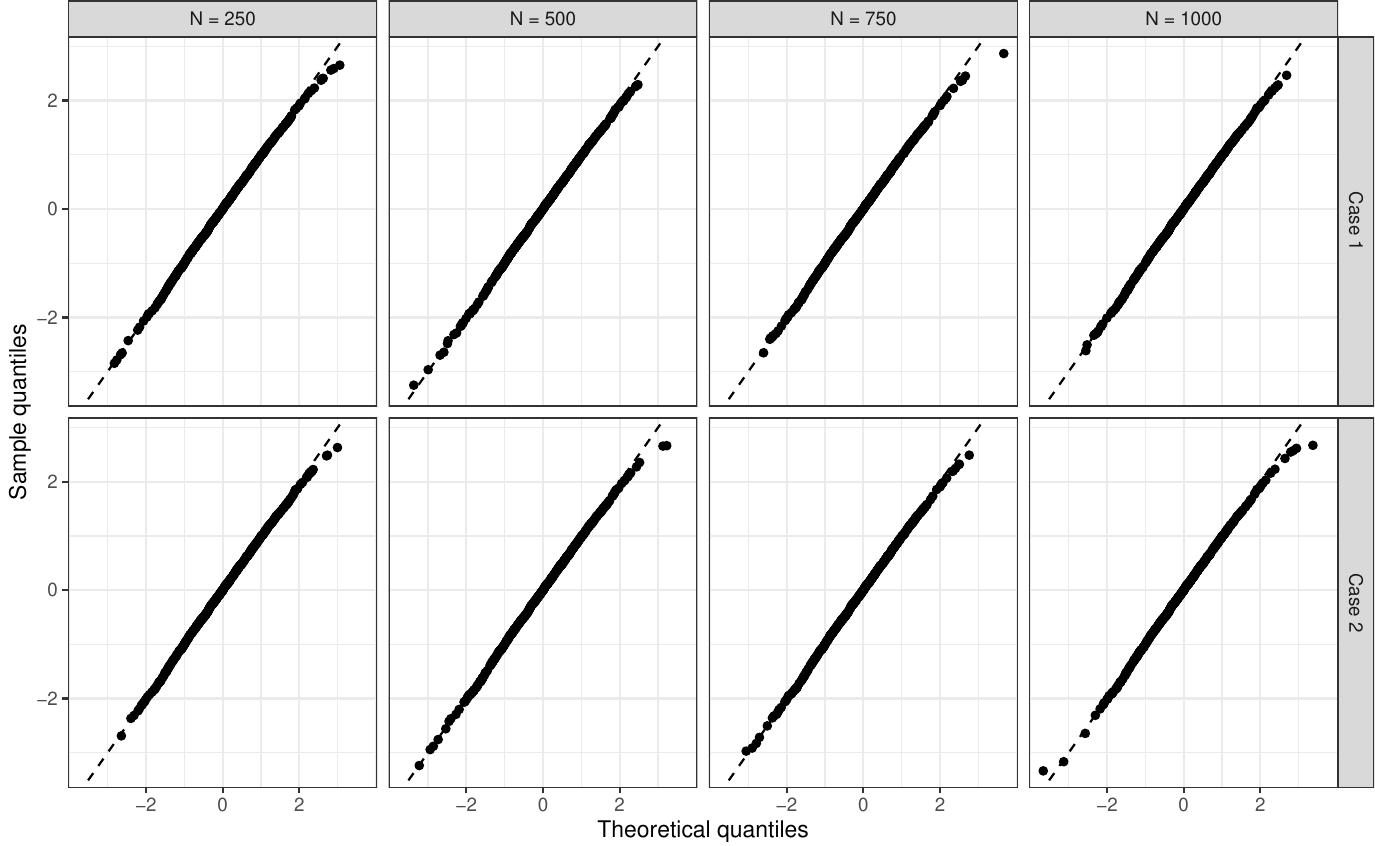}
\caption{\label{fig:qq} Quantile-Quantile plots showing the results of Simulation Study 2. 
The sample quantiles of the standardized  maximum likelihood estimates of the transitive edge parameter 
are plotted against the theoretical quantiles based on the standard normal approximation in each 
of the two cases studied in Simulation study 2 across networks of size $N \in \{250, 500, 750, 1000\}$.}
\end{figure}

\section{Conclusions}

In this work, 
we have proved the first rigorous theory for both estimation and statistical inference of local dependence random graph models.
We have established minimax optimal rates of convergence in the $\ell_2$-norm 
of maximum likelihood estimators of exponential-family local dependence random graph models,
accompanying these results with finite-sample error bounds on the multivariate normal approximation 
of a standardization of maximum likelihood estimators.
Notably, 
our results allow for both the number of parameters and the sizes of blocks to grow unbounded 
with the size of the network.  

Our consistency and normal approximation theory are non-asymptotic, 
although we have stated helpful asymptotic results along the way, 
which enable us to understand how key aspects of the model (through the spectrum of Fisher information matrices
and the dimension of parameter vectors) 
and properties of the network (through the number and sizes of blocks and nodes)
impact rates of convergence for both the statistical error (in the $\ell_2$-norm)
and the multivariate normal approximation. 
Our results cover general settings and heterogeneous parameterizations, 
as exemplified in the examples in Section \ref{sec:examples} and our simulation studies in Section \ref{sec:simulation},
which allow our results to cover a broad scope. 

Results were derived under the assumption that we have perfect knowledge of the block memberships of nodes in the network. 
This may be reasonable in certain settings where we can observe the block memberships of nodes, 
but might be violated in other settings where we obtain imperfect observations of the block memberships of nodes, 
whether through a noisy observation process or error in the estimates of the block memberships. 
The effect of imperfect knowledge of the block memberships of nodes on the aforementioned errors and convergence rates 
is an open question.

\begin{acks}[Acknowledgments]
The author is grateful to two anonymous reviewers,
an associate editor, the editor, 
and to Dr. Michael Schweinberger,
for constructive comments and suggestions
that have greatly improved the paper.
\end{acks}

\begin{funding}
The author was supported by NSF award SES-2345043 and 
the Department of Defense Test Resource Management Center under contracts FA8075-18-D-0002 and FA8075-21-F-0074. 
\end{funding}

\begin{supplement}
\stitle{Supplement to ``Rates of convergence and normal approximations for estimators of local dependence random graph models}
\sdescription{The supplementary material contains the proofs for all results in Section \ref{sec:theory},
as well as additional technical results and the proofs of these technical results used to prove the main results in Section \ref{sec:theory}} 
\end{supplement}

\bibliographystyle{imsart-nameyear.bst}
\bibliography{library}

%\newpage
%\begin{appendices}
%\input{supplement.tex}
%\end{appendices}
%\label{last.page}

\newpage 

\begin{frontmatter}
\title{Supplement to ``Rates of convergence and normal approximations for estimators of local dependence random graph models''}
\runtitle{Supplement to ``Rates of convergence for local dependence random graph models''}

\begin{aug}
\author[A]{\inits{F.}\fnms{Jonathan R.}~\snm{Stewart}\ead[label=e1]{jrstewart@fsu.edu}}
\address[A]{Department of Statistics,
Florida State University,
Tallahassee, FL, USA\printead[presep={,\ }]{e1}}
\end{aug}

\end{frontmatter}

\makeatletter

\setcounter{page}{1}

\setcounter{section}{0}

\setcounter{com}{0}

\s\s

\noindent
{Section \ref{L2}: Proof of Theorem 2.1 \dotfill\pageref{L2}}\s\\
\mbox{\hspace{.15in}} {Section \ref{cons_aux}: Auxiliary results for Theorem 2.1 \dotfill\pageref{cons_aux}}\s\\
{Section \ref{minimax_lower}: Proof of Theorem 2.2 \dotfill\pageref{minimax_lower}}\s\\
{Section \ref{minimax_upper}: Proof of Theorem 2.3 \dotfill\pageref{minimax_upper}}\s\\
{Section \ref{cor_minimax}: Proof of Corollary 2.4 \dotfill\pageref{cor_minimax}}\s\\
{Section \ref{norm}: Proof of Theorem 2.5 \dotfill\pageref{norm}}\s\\
\mbox{\hspace{.15in}} {Section \ref{norm_aux}: Auxiliary results for Theorem 2.5 \dotfill\pageref{norm_aux}}\s\\
{Section \ref{sup-sec:prop_I}: Proof of Theorem 2.7 \dotfill\pageref{sup-sec:prop_I}}\s\\
\mbox{\hspace{.15in}} {Section \ref{sec:aux_2.7}: Auxiliary results for Theorem 2.7 \dotfill\pageref{sec:aux_2.7}}\s\\
{Section \ref{sup-sec:exp_fam}: Auxiliary results for exponential families \dotfill\pageref{sup-sec:exp_fam}}\s\\

\s

\hide{
We establish a number of  auxiliary results in the supplementary materials which are used in the proofs of 
Theorems \ref{thm:L2} and \ref{thm:normal},
as well as various basic properties of exponential families which are used throughout the main manuscript 
and this supplement. 
}

\section{Proof of Theorem 2.1}
\label{L2}

Our method of proof utilizes a general M-estimation argument. 
For ease of presentation, 
we first introduce the general argument and then apply 
the general argument to maximum likelihood estimators 
of exponential-family local dependence random graph models 
to obtain the results of Theorem 2.1.

\s

{\bf General M-estimation framework for rates of convergence.}
Consider a random estimating function $m : \mbR^d \times \mbX \mapsto \mbR$ 
and define $M(\nat) \coloneqq \mbE \, m(\nat, \bX)$ for $\nat \in \mbR^d$.
We make the following assumptions concerning $m(\nat, \bx)$ and $M(\nat)$:
\ben 
\item Assume that $m(\nat, \bx)$ is concave in $\nat \in \mbR^d$ and continuously differentiable 
at all $\nat \in \mbR^d$ and for all $\bx \in \mbX$.

\vspace{.1cm}
  
\item Assume that $M(\nat)$ is strictly concave in $\nat \in \mbR^d$
and that $\truth \in \mbR^d$ is the unique global maximizer of $M(\nat)$.

\vspace{.1cm}

\item Assume that $M(\nat)$ is twice continuously differentiable
and that there exists an $\epsilon > 0$ (fixed) such that 
the negative Hessian $\bH(\nat) \coloneqq -\nabla_{\nat}^2 \, M(\nat)$ 
of $M(\nat)$ is positive definite for all $\nat \in \mB_2(\truth, \epsilon)$.  
 
\een
When $m(\nat, \bx)$ is the log-likelihood corresponding to a minimal exponential family,
standard exponential family theory establishes that the above conditions (1) and (2) hold 
(e.g., Proposition 3.10 of \citep{Sundberg2019}). 
As a result, 
$\nabla_{\nat} \, M(\truth) = \bm{0}_d$,
where $\bm{0}_d$ is the $d$-dimensional zero vector.
By Theorem 6.3.4 of \citet{Ortega2000}, 
if the event (for $\delta > 0$) 
\be
\label{eq:key_event}
\inf\limits_{\nat \in \partial\mB_2(\truth, \, \delta)} \, 
\langle \nat - \truth, \; \nabla_{\nat} \, m(\nat, \bX) \rangle 
&\geq& 0
\ee
occurs,  
where $\partial\mB_2(\truth, \, \delta)$ denotes the boundary of 
$\mB_2(\truth, \, \delta) \coloneqq \{ \nat \in \mbR^d : \norm{\nat - \truth}_2 < \delta\}$,  
then a root of $\nabla_{\nat} \, m(\nat, \bX)$ exists in 
$\bar{\mB_2}(\truth, \delta) \coloneqq \mB_2(\truth, \, \delta) \cup \partial \mB_2(\truth, \, \delta)$,
in which case a global maximizer $\nat_0 = \argmax_{\nat \in \mbR^d} \, m(\nat, \bX)$ exists and satisfies 
$\norm{\nat_0 - \truth}_2 \leq \delta$. 

The key to our approach lies in demonstrating that condition \eqref{eq:key_event} holds with high probability 
for a chosen $\delta \in (0, \epsilon)$ 
($\epsilon > 0$ fixed) 
which helps to establish rates of convergence of estimators. 
In order to do so,
we leverage the multivariate mean value theorem to establish that there exists,
for each parameter vector $\nat \in \partial \mB_2(\truth, \, \delta)$, 
a parameter vector  
\beno 
\dot\nat 
\= t \, \nat + (1-t) \, \truth 
&\in& \mB_2(\truth, \, \delta)
&\subset& \mB_2(\truth, \, \epsilon),
&& \mbox{for some } t \in (0, 1),
\ee
such that   
\beno
\langle \nat - \truth, \; \nabla_{\nat} \, M(\nat) \rangle 
\= \langle \nat - \truth, \; \nabla_{\nat} \, M(\truth) \rangle 
+ \langle \nat - \truth, \; \bH(\dot\nat) \, (\nat - \truth) \rangle \s\\ 
\= \langle \nat - \truth, \; \bH(\dot\nat) \, (\nat - \truth) \rangle,
\ee
recalling that $\nabla_{\nat} \, M(\truth) = \bm{0}_d$. 
Observe that 
\beno
\langle \nat - \truth, \; \bH(\dot\nat) \, (\nat - \truth) \rangle
\=
\dfrac{\langle \nat - \truth, \; \bH(\dot\nat) \, (\nat - \truth) \rangle}{\langle \nat - \truth, \, \nat - \truth \rangle} \, \norm{\nat - \truth}_2^2 
&\geq& \lambda_{\min}(\bH(\dot\nat)) \, \delta^2, 
\ee
noting that the Rayleigh quotient of $\bH(\dot\nat)$ is bounded below by the smallest eigenvalue 
$\lambda_{\min}(\bH(\dot\nat))$ of $\bH(\dot\nat)$
and that $\norm{\nat - \truth}_2 = \delta$ for all $\nat \in \partial\mB_2(\truth, \delta)$.
As a result,
\be
\label{eq:key_bound} 
\inf\limits_{\nat \in \partial\mB_2(\truth, \delta)} \, 
\langle \nat - \truth, \; \nabla_{\nat} \, M(\nat) \rangle
&\geq&  \inf\limits_{\nat \in \mB_2(\truth, \delta)} \, \lambda_{\min}(\bH(\nat)) \, \delta^2.
\ee
As $\partial\mB_2(\truth, \delta) \subset \mB_2(\truth, \epsilon)$,  
\be
\label{eq:positive}
\inf\limits_{\nat \in \mB_2(\truth, \delta)} \, \lambda_{\min}(\bH(\nat)) 
&\geq&
\inf\limits_{\nat \in \mB_2(\truth, \epsilon)} \, \lambda_{\min}(\bH(\nat)) 
&>& 0, 
\ee
by the assumption that $\bH(\nat)$ is positive definite, 
and thus non-singular,
on $\mB_2(\truth, \epsilon)$;
Assumption 2 ensures this condition for  maximum likelihood estimators 
of exponential-family local dependence random graph models. 
As a result of \eqref{eq:key_bound} and \eqref{eq:positive}, 
the event 
\be
\label{eq:key_condition}
\sup\limits_{\nat \in \partial \mB_2(\truth, \delta)}
\left| \langle \nat - \truth, \; \nabla_{\nat} \, m(\nat, \bX) - \nabla_{\nat} \, M(\nat) \rangle \right|
\,\leq\, \inf\limits_{\nat \in \mB_2(\truth, \epsilon)} \, \lambda_{\min}(\bH(\nat)) \, \delta^2
\ee
implies the event \eqref{eq:key_event}. 
Thus, 
demonstrating that event \eqref{eq:key_condition} occurs with probability 
at least $1 - N^{-2}$ 
demonstrates that event \eqref{eq:key_event} occurs with probability at least $1 - N^{-2}$. 

\s\s

{\bf Rates of convergence for maximum likelihood estimators.} 
The log-likelihood equation of an exponential-family local dependence random graph model has the form
\beno
\ell(\btheta, \bx)
\= 
\dsum_{k=1}^{K} \ell_{k,k}(\nat_W, \bx_{k,k}) + 
\dsum_{1 \leq k < l \leq K} \ell_{k,l}(\nat_B, \bx_{k,l}), 
\ee
which implies that the maximum likelihood estimator $\mle = (\mle_{W}, \mle_{B})$
is given by
\begin{equation}\begin{array}{lrlcrl}
\label{eq:optimizers}
\mle_{W} \= \argmax\limits_{\nat_W \in \mbR^p} \,
& \dsum_{k=1}^{K} & \ell_{k,k}(\nat_W, \bx_{k,k}) \s\\
\mle_B
\= \argmax\limits_{\nat_B \in \mbR^q} \,
& \dsum_{1 \leq k < l \leq K} & \ell_{k,l}(\nat_B, \bx_{k,l}),
\end{array}
\end{equation} 
owing to the fact that the subgraphs $\bX_{k, l}$ ($1 \leq k \leq l \leq K$) are independent
and that the parameter vectors $\nat_W$ and $\nat_B$ partition the parameters in $\nat$. 
Hence,
each optimizer in \eqref{eq:optimizers} can be found separately and independently.
Define 
\begin{equation}\begin{array}{lrccrl} \nonumber 
m_W(\nat_W, \bx_W) 
&\coloneqq& \dsum_{k=1}^{K} &\, \ell_{k,k}(\nat_W, \bx_{k,k}) \s\\
m_B(\nat_B, \bx_B) 
&\coloneqq& \dsum_{1 \leq k < l \leq K} & \, \ell_{k,l}(\nat_B, \bx_{k,l}), 
\end{array}
\end{equation}
$M_W(\nat_W) \coloneqq \mbE \, m_W(\nat_W, \bX_W)$,  
and $M_B(\nat_B) \coloneqq \mbE \, m_B(\nat_B, \bX_B)$,
where 
\beno
\bX_W \coloneqq (\bX_{1,1}, \ldots, \bX_{K,K}) 
&\mbox{and}& 
\bX_B \coloneqq (\bX_{1,2}, \ldots, \bX_{1,K}, \bX_{2,3}, \bX_{2,4}, \ldots, \bX_{K-1,K}). 
\ee 
Due to the above considerations, 
\beno
\bH(\nat) 
&\coloneqq& -\mbE \, \nabla_{\nat}^2 \, \ell(\nat, \bX)
\= \left( \begin{matrix}
\bH_W(\nat_W) & \bm{0}_{p,q} \s\\
\bm{0}_{q,p} & \bH_B(\nat_B)
\end{matrix} \right), \s
\ee
where $\bm{0}_{d,r}$ is the $(d \times r)$-dimensional matrix of all zeros, and where  
\beno
\bH_W(\nat_W) \,\coloneqq\, \dsum_{k=1}^{K} \, \bH_{k,k}(\nat_W)
&&\mbox{and}&&
\bH_B(\nat_B) \,\coloneqq\, \dsum_{1 \leq k < l \leq K} \, \bH_{k,l}(\nat_B),
\ee
with the definitions 
\beno
\bH_{k,k}(\nat_W) &\coloneqq& - \mbE \, \nabla_{\nat_W}^2 \ell_{k,k}(\nat_W, \bX_{k,k}),
&& \mbox{for all } 1 \leq k \leq K \s\\ 
\bH_{k,l}(\nat_B) &\coloneqq& - \mbE \, \nabla_{\nat_B}^2 \, \ell_{k,l}(\nat_B, \bX_{k,l}),
&& \mbox{for all } 1 \leq k < l \leq K. 
\ee
Note that the interchange of differentiation and integration in this setting is trivial 
as the expectations are finite sums. 

We demonstrate that event \eqref{eq:key_condition} occurs 
with probability at least $1 - N^{-2}$
for the within-block and between-block cases separately.
Assumption 2 ensures that  
\be
\label{eq:lmin_lower}
\inf\limits_{\nat \in \mB_2(\truth_W, \epsilon)} \, 
\lambda_{\min}(\bH_W(\nat_W))
&=& K \, \minw 
&>& 0 \s\\
\inf\limits_{\nat \in \mB_2(\truth_B, \epsilon)} \, 
\lambda_{\min}(\bH_B(\nat_B))
&=& \dbinom{K}{2} \, \minb 
&>& 0.  
\ee 
\hide{
an application of Weyl's inequality and Assumption (A.1) shows that 
\begin{equation}\begin{array}{llclll} 
\label{eq:lmin_lower} 
\lambda_{\min}(\bH_W(\nat_W))
&\geq& \dsum_{k=1}^{K} & \, \lambda_{\min}(\bH_{k,k}(\nat_W))
&\geq& \minw \, K \s\\ 
\lambda_{\min}(\bH_B(\nat_B))
&\geq& \dsum_{1 \leq k < l \leq K} & \, \lambda_{\min}(\bH_{k,l}(\nat_B))
&\geq& \minb \, \dbinom{K}{2}.
\end{array}
\end{equation}
}
Let $\delta_W \in (0, \, \epsilon \,/\, \sqrt{2})$ and 
$\delta_B \in (0, \, \epsilon \,/\, \sqrt{2})$, 
and assume $\nat_W \in \partial\mB_2(\truth_W, \delta_W)$ 
and $\nat_B \in \partial\mB_2(\truth_B, \delta_B)$.
By assumption, 
$(\nat_W,\nat_B) \in \mB_2(\truth,\epsilon)$. 
Thus,  
using \eqref{eq:lmin_lower},
we can rewrite the events in \eqref{eq:key_condition} as events  
\beno
& |\langle \nat_W - \truth_W, \, \nabla_{\nat_W} \,  m_W(\nat_W, \bX_W) - \nabla_{\nat_W} \, M_W(\nat_W) \rangle|
&\leq& \delta_W^2 \, K \, \minw \s\\
& |\langle \nat_B - \truth_B, \, \nabla_{\nat_B} \, m_B(\nat_B, \bX_B) - \nabla_{\nat_B} \, M_B(\nat_B) \rangle |
&\leq& \delta_B^2 \,  \dbinom{K}{2} \, \minb. 
\ee
%using \eqref{eq:lmin_lower}.
By the Cauchy-Schwarz inequality, 
\beno
&& |\langle \nat_W - \truth_W, \, \nabla_{\nat_W} \, m_W(\nat_W, \bX_W) - \nabla_{\nat_W} \, M_W(\nat_W) \rangle| \s\\
&\leq& \norm{\nat_W - \truth_W}_2 \, \norm{\nabla_{\nat_W}\, m_W(\nat_W, \bX_W) - \nabla_{\nat_W} \, M_W(\nat_W)}_2 \s\\
\= \delta_W \, \norm{\nabla_{\nat_W}\, m_W(\nat_W, \bX_W) - \nabla_{\nat_W} \, M_W(\nat_W)}_2.  
\ee
Similarly, 
\beno
&& |\langle \nat_B - \truth_B, \, \nabla_{\nat_B} \, m_B(\nat_B, \bX_B) - \nabla_{\nat_B}\, M_B(\nat_B) \rangle | \s\\
&\leq& \delta_B \, \norm{\nabla_{\nat_B} \, m_B(\nat_B, \bX_B) - \nabla_{\nat_B} \, M_B(\nat_B)}_{2}. 
\ee
It therefore suffices to demonstrate,
for all $\nat_W \in \partial\mB_2(\truth_W, \delta_W)$
and $\nat_B \in \partial\mB_2(\truth_B, \delta_B)$, 
that events  
\be
\label{eq:mle_key}
\norm{\nabla_{\nat_W} \, m_W(\nat_W, \bX_W) - \nabla_{\nat_W} \, M_W(\nat_W)}_{2} 
&\leq& \delta_W \, K \, \minw  \s\\
\norm{\nabla_{\nat_B} \, m_B(\nat_B, \bX_B) - \nabla_{\nat_B}\, M_B(\nat_B)}_{2} 
&\leq& \delta_B \, \tbinom{K}{2} \, \minb
\ee
occur with probability at least $1 - N^{-2}$. 
Define,
for all $t > 0$,
 the events
\beno
\mD_W(t) &\coloneqq& \left\{\bx \in \mbX \,:\, 
\sup\limits_{\nat_W \in \partial\mB_2(\truth_W, \delta_W)}
\norm{\nabla_{\nat_W} \, m_W(\nat_W, \bX_W) - \nabla_{\nat_W} \, M_W(\nat_W)}_{2} \,\geq\, t \right\} \s\\
\mD_B(t) &\coloneqq& \left\{\bx \in \mbX \,:\, 
\sup\limits_{\nat_B \in \partial\mB_2(\truth_B, \delta_B)}
\norm{\nabla_{\nat_B} \, m_B(\nat_B, \bX_B) -  \nabla_{\nat_B} \, M_B(\nat_B)}_{2} \,\geq\, t \right\}.
\ee
By Lemma \ref{lem:sup_mle}, 
\beno
\sup\limits_{\nat_W \in \partial\mB_2(\truth_W, \delta_W)}
\norm{\nabla_{\nat_W} \, m_W(\nat_W, \bX_W) - \nabla_{\nat_W} \, M_W(\nat_W)}_{2}
\= \norm{\nabla_{\nat_W} \, \ell_W(\truth_W, \bX_W)}_2 \s\\
\sup\limits_{\nat_B \in \partial\mB_2(\truth_B, \delta_B)}
\norm{\nabla_{\nat_B} \, m_B(\nat_B, \bX_B) - \nabla_{\nat_B} \, M_B(\nat_B)}_{2}
\= \norm{\nabla_{\nat_B} \, \ell_B(\truth_B, \bX_B)}_2, 
\ee
and 
applying Lemma \ref{lem:conc_L2},
we obtain the bounds  
\beno
\scalebox{1}{$\mbP\left(\mD_W\left(\delta_W \, K \, \minw  \right) \right)$}
&\leq& \scalebox{1}{$ 
\exp\left(- \dfrac{\delta_W^2 \, (\minw)^2 \, K^2}
{5 \, K \, \maxw + C_W \, A_{\max}^2 \, \sqrt{p} \, \delta_W \, \minw \, K} + \log(5) \, p \right)
$}
\ee
and  
\beno 
\scalebox{1}{$ 
\mbP\left(\mD_B\left(\delta_B \, \tbinom{K}{2} \, \minb \right) \right)$}
&\leq& 
\scalebox{1}{$
\exp\left(- \dfrac{\delta_B^2 \, (\minb)^2 \, \tbinom{K}{2}^2}
{5 \, \tbinom{K}{2} \, \maxb + 2 \, C_B \, A_{\max}^2 \, \sqrt{q} \, \delta_B \, \minb \, \tbinom{K}{2}}
+ \log(5) \, q \right).$} 
\ee
Choosing 
\beno
\delta_W 
\= \beta_W  \, \dfrac{\sqrt{\maxw}}{\minw} \, \sqrt{\dfrac{p}{K}}
&>& 0, 
\ee
for a value of $\beta_W \in (0, \infty)$ to be given, 
establishes 
\beno
\scalebox{1}{$\mbP\left(\mD_W\left(\delta_W \, \minw \, K \right) \right)$}
&\leq& \scalebox{1}{$
\exp\left(- \dfrac{\beta_W^2 \, p \, K \, \maxw}
{5 \, K \, \maxw + C_W \, \beta_W \, A_{\max}^2 \, p \, \sqrt{K \, \maxw}} + \log(5) \, p \right).
$} 
\ee
Using Assumption 4,
the assumption that  
\beno
A_{\max}
\;\leq\; 
\left(\dfrac{N \; \maxw}{\Aavg \; p^2} \right)^{1/4}
%\left(\dfrac{N}{\Aavg} \right)^{1/4} \; \sqrt{\dfrac{\maxw}{p}}
\;&\mbox{implies}&\;
A_{\max}^2  \,
p \, \sqrt{K \, \maxw}
\;\leq\;
K \, \maxw,
\ee
defining $\Aavg \coloneqq K^{-1} \sum_{k=1}^{K} |\mA_k|$ and
using the identity
\be
\label{eq:NK}
N
\= \dsum_{k=1}^{K} |\mA_k|
\= K \, \dfrac{1}{K} \, \dsum_{k=1}^{K} |\mA_k|
\= K \, \Aavg.
\ee
Hence,
\beno
&& \exp\left(- \dfrac{\beta_W^2 \, p \, K \, \maxw}
{5 \, K \, \maxw + C_W \, \beta_W \, A_{\max}^2 \, p \, \sqrt{K \, \maxw}} + \log(5) \, p \right) \s\s\\
&\leq&  \exp\left(- \dfrac{\beta_W^2 \, p \, K \, \maxw}{(5 + \beta_W \, C_W) \, K \, \maxw} + \log(5) \, p \right) \s\\
\= \exp\left(\left( -\dfrac{\beta_W^2}{5 + \beta_W \, C_W} + \log(5) \right) \, p \right).  
\ee
To obtain the desired probability guarantee, 
we require 
\beno
-\dfrac{\beta_W^2}{5 + \beta_W \, C_W} + \log(5) 
\=  -2, 
\ee
which in turn requires a solution $\beta_W \in (0, \infty)$ to the quadratic equation 
\beno
\beta_W^2 - C_W \, (2 + \log(5)) \, \beta_W - 5 \, (2 + \log(5))
\= 0.
\ee 
Using the quadratic formula, 
such a root is given by 
\beno
\beta_W 
\= \dfrac{C_W \, (2 + \log(5)) + \sqrt{C_W^2 \, (2 + \log(5))^2 \, + 20 \, (2 + \log(5))}}{2} 
&>& 0, 
\ee
which in turn establishes that 
\beno
\mbP\left(\mD_W\left(\delta_W \, \minw \, K \right) \right)
&\leq& \exp\left( - 2 \, p\right). 
\ee
Under the assumption that $p \geq \log(N)$,
\beno
\mbP\left(\mD_W\left(\delta_W \, \minw \, K \right) \right)
&\leq& \dfrac{1}{N^2}. 
\ee
Similarly, 
choosing 
\beno
\delta_B 
\= 
\beta_B \, \dfrac{\sqrt{\maxb}}{\minb} \, \sqrt{\dfrac{q}{\tbinom{K}{2}}}
&>& 0, 
\ee
for a value of $\beta_B \in (0, \infty)$ to be given, 
establishes 
\beno
\scalebox{1}{$\mbP\left(\mD_B\left(\delta_B \, \minb \, \tbinom{K}{2} \right) \right)$}
&\leq& \scalebox{1}{$
\exp\left(- \dfrac{\beta_B^2 \, q \, \tbinom{K}{2} \, \maxb}
{5 \, \tbinom{K}{2} \, \maxb + 2 \, C_B \, \beta_B \,  A_{\max}^2 \, q \, \sqrt{\tbinom{K}{2} \, \maxb}} + \log(5) \, q \right).
$}
\ee
Using Assumption 4, 
the assumption that
\beno
A_{\max}
\;\leq\;
\left( \dfrac{N^2 \, \maxb}{4 \, \Aavg^2 \, q^2} \right)^{1/4}
&&\mbox{implies}&&
A_{\max}^2 \,
q \, \sqrt{\tbinom{K}{2} \, \maxb}
\;\leq\;
\tbinom{K}{2} \, \maxb,
\ee
once more using the identity in \eqref{eq:NK}.  
Hence,
\beno
&& \exp\left(- \dfrac{\beta_B^2\, q \, \tbinom{K}{2} \, \maxb}
{5 \, \tbinom{K}{2} \, \maxb + 2 \, C_B \, \beta_B \, A_{\max}^2 \, q \, \sqrt{\tbinom{K}{2} \, \maxb}} + \log(5) \, q \right) \s\s\\
&\leq&  \exp\left(- \dfrac{\beta_B^2 \, q \, \tbinom{K}{2} \, \maxb}{(5 + 2 \, C_B \, \beta_B) \, \tbinom{K}{2} \, \maxb} + \log(5) \, q \right) \s\\
\= \exp\left( \left(-\dfrac{\beta_B^2}{5 + 2 \, C_B \, \beta_B} + \log(5) \right) \, q \right).  
\ee
To obtain the desired probability guarantee, 
we require 
\beno
-\dfrac{\beta_B^2}{5 + 2 \, C_B \, \beta_B } + \log(5) 
\=  -2, 
\ee
which in turn requires a solution $\beta_B \in (0, \infty)$ to the quadratic equation 
\beno
\beta_B^2 - 2 \, C_B \, (2 + \log(5)) \, \beta_B - 5 \, (2 + \log(5))
\= 0.
\ee 
Using the quadratic formula, 
such a root is given by 
\beno
\beta_B 
\= C_B \, (2 + \log(5)) + \sqrt{C_B^2 \, (2 + \log(5))^2 \, + 5 \, (2 + \log(5))} 
&>& 0, 
\ee
which in turn establishes that 
\beno
\mbP\left(\mD_B\left(\delta_W \, \minb \, \tbinom{K}{2} \right) \right)
&\leq& \exp\left( - 2 \, q\right). 
\ee
Under the assumption that $q \geq \log(N)$,
\beno
\mbP\left(\mD_W\left(\delta_B \, \minb \, \tbinom{K}{2} \right) \right)
&\leq& N^{-2}. 
\ee
As a result, 
event \eqref{eq:mle_key} occurs with probability at least $1 - N^{-2}$,
implying that,
with probability at least $1 - N^{-2}$,
the maximum likelihood estimator $\mle = (\mle_W, \mle_B)$ exists uniquely and satisfies 
\beno
\norm{\mle_W - \truth_W}_2 &\leq& \beta_W \, \dfrac{\sqrt{\maxw}}{\minw} \, \sqrt{\dfrac{p}{K}} \s\\ 
\norm{\mle_B - \truth_B}_2 &\leq& \beta_B \, \dfrac{\sqrt{\maxb}}{\minb} \, \sqrt{\dfrac{q}{\tbinom{K}{2}}}. 
\ee 
Uniqueness of $(\mle_W, \mle_B)$ 
follows from the assumption that the exponential-family local dependence random graph model
is minimal  
(Proposition 3.13 of \citet{Sundberg2019}).
We convert the bounds in terms of $K$ to $N$ by utilizing \eqref{eq:NK} again to show that 
\beno
\norm{\mle_W - \truth_W}_2  
&\leq& \beta_W \, \dfrac{\sqrt{\maxw}}{\minw} \, \sqrt{\dfrac{p}{K}}
\= C_1 \, \sqrt{\Aavg} \, \dfrac{\sqrt{\maxw}}{\minw} \, \sqrt{\dfrac{p}{N}} \s\\
\norm{\mle_B - \truth_B}_2 
&\leq& \beta_B \, \dfrac{\sqrt{\maxb}}{\minb} \, \sqrt{\dfrac{q}{\tbinom{K}{2}}}
&\leq& C_2 \, \Aavg \, \dfrac{\sqrt{\maxb}}{\minb} \, \sqrt{\dfrac{q}{N^2}},
\ee
using $\tbinom{K}{2} \geq K^2 \,/\, 4$ in the second case, 
and defining $C_1 \coloneqq \beta_W > 0$ and 
$C_2 \coloneqq 2 \, \beta_B > 0$.  
Both are independent of $N$, $p$, and $q$.  

Finally, 
we show the restriction to $\mB_2(\truth, \epsilon)$ to be legitimate. 
Assumption 3  ensures  
\beno
\lim\limits_{N \to \infty} \, 
\sqrt{\Aavg} \, \dfrac{\sqrt{\maxw}}{\minw} \, \sqrt{\dfrac{p}{N}}
\;=\; 0 
&&\mbox{and}&&
\lim\limits_{N \to \infty} \, 
\Aavg \, \dfrac{\sqrt{\maxb}}{\minb} \, \sqrt{\dfrac{q}{N^2}}
\;=\; 0.
\ee
As such, 
there exists an $N_0 \geq 3$ such that,
for all integers $N \geq N_0$, 
we have $\max\{\delta_W, \delta_B\} < \epsilon / \sqrt{2}$. 
Thus,
for all integers $N \geq N_0$ and  
with probability at least $1 - N^{-2}$, 
the unique vector $\mle \in \Mle$ satisfies  
\beno
\norm{\mle - \truth}_2
&=& \sqrt{\norm{\mle_W - \truth_W}_2^2 + \norm{\mle_B -\truth_B}_2^2}
&\leq& \sqrt{\delta_W^2 + \delta_B^2}
&<& \epsilon,  
\ee
which implies,
for all integers $N \geq N_0$, that  
\beno
\mbP(\norm{\mle - \truth}_2 \leq \epsilon)
&\geq& 1 - N^{-2}, 
\ee  
justifying the restriction to the subset of the parameter space $\mB_2(\truth, \epsilon) \subset \mbR^{p+q}$. 
\qed

\subsection{Auxiliary results for Theorem 2.1}
\label{cons_aux}

We prove a concentration inequality for gradients of the log-likelihood which is utilized in the proof 
of Theorem 2.1. 

\s

\begin{lemma}
\label{lem:conc_L2}
Under the assumptions of Theorem 2.1,
\beno
\scalebox{1}{$
\mbP\left( \norm{\nabla_{\btheta_W} \, \ell_W(\truth_W, \bX_W)}_2 \,\leq\, \delta \right)$}
&\geq& \scalebox{.9}{$1 - \exp\left( - \dfrac{\delta^2}{5 \, K \, \maxw
+ C_W A_{\max}^2 \, \sqrt{p} \, \delta} + \log(5) \, p \right)$} \s\s\\
\scalebox{1}{$\mbP\left( \norm{\nabla_{\btheta_B} \, \ell_B(\truth_B, \bX_B)}_2 \,\leq\, \delta \right)$}
&\geq& \scalebox{.9}{$1 - \exp\left( - \dfrac{\delta^2}{5 \, \tbinom{K}{2} \, \maxb
+ 2 \, C_B A_{\max}^2 \, \sqrt{q} \, \delta} + \log(5) \, q \right)$}, 
\ee
for all $\delta > 0$,
where $C_W > 0$ and $C_B > 0$ are the same constants as in Assumption 1.  
\end{lemma}

\llproof \ref{lem:conc_L2}. 
We first prove the result for the within-block case, 
and then discuss extensions to prove the result for the between-block case, 
noting that the two proofs are essentially the same with only a couple of notational changes. 

Following the method utilized in the proof of Lemma 8.4 of \citet{Chen2022}, 
define  
$\mcU \coloneqq \left\{  \bu \in \mbR^p : \norm{\bu}_2 \leq 1 \right\}$
to be the closed unit ball in $\mbR^p$. 
By Corollary 4.2.13 (p. 78) of \citet{Vershynin2018}, 
there exists a subset $\mcV_{\epsilon} \subset \mcU$ (for $\epsilon \in (0, 1)$)
which is an $\epsilon$-net of $\mcU \subset \mbR^p$ such that 
the cardinality of the set $\mcV_{\epsilon}$ satisfies  
$\log  |\mcV_{\epsilon}| \leq p \, \log(2 \, \epsilon^{-1} + 1)$. 
Taking $\epsilon = 1/2$,
there exists,  
for each $\bu \in \mcU$, 
a $\bv \in \mcV_{1/2}$ satisfying $\norm{\bu - \bv}_2 \leq 1/2$,
where $\log |\mcV_{1/2}| \leq p \, \log(5)$. 
For ease of presentation, 
define 
\beno
\bG &\coloneqq& \nabla_{\nat_W} \, \ell_W(\truth_W, \bX_W).  
\ee
For any $\bu \in \mcU$,
with the corresponding $\bv \in \mcV_{1/2}$, 
the Cauchy-Schwarz inequality implies  
\be
\label{eq:ChenGaoZhangBridge} 
\langle \bu, \, \bG \rangle 
\= \langle \bv, \, \bG \rangle + \langle \bu - \bv, \, \bG \rangle  \s\\
&\leq& \langle \bv, \, \bG \rangle + \norm{\bu - \bv}_2 \, \norm{\bG}_2 \s\\ 
&\leq& \langle \bv, \, \bG \rangle + \dfrac{1}{2} \, \norm{\bG}_2,
\ee
using the fact that $\norm{\bu - \bv}_2 \leq 1/2$ in the last line. 
Next, 
choosing 
\beno
u_i 
\= \dfrac{G_i}{\norm{\bG}_2},
&& i \in \{1, \ldots, p\}, 
\ee
ensures that $\norm{\bu}_2 \leq 1$ so that the chosen $\bu$ exists in $\mcU$.  
By writing 
\beno
\langle \bu, \, \bG \rangle 
\= \dfrac{1}{\norm{\bG}_2} \, \dsum_{i=1}^{p} \, G_i^2 
\= \dfrac{\norm{\bG}_2^2}{\norm{\bG}_2}
\= \norm{\bG}_2, 
\ee
we revisit \eqref{eq:ChenGaoZhangBridge} to obtain,
using the above identity and re-arrangement, 
the inequality 
\be
\label{eq:ChenGaoZhangKey} 
\norm{\bG}_2 
&\leq& 2 \, \max\limits_{\bv \in \mcV_{1/2}} \, \langle \bv, \, \bG \rangle,
\ee
where we take the maximum over $\bv \in \mcV_{1/2}$
since we cannot be sure which $\bv \in \mcV_{1/2}$ would correspond to our choice of $\bu \in \mcU$ above. 
A quick remark is in order regarding the case when $\norm{\bG}_2 = 0$. 
Note that the above implicitly assumed $\norm{\bG}_2 \neq 0$. 
In the event where $\norm{\bG}_2 = 0$,
the inequality \eqref{eq:ChenGaoZhangKey} remains true trivially, 
because $\bG = \bm{0}_p$ and $\langle \bv, \, \bG \rangle = 0$ for all $\bv \in \mcV_{1/2}$. 
As a result of \eqref{eq:ChenGaoZhangKey} and for $\delta > 0$,  
\beno
\mbP\left(\norm{\nabla_{\btheta_W} \, \ell_W(\truth, \bX)}_2 \leq \,\delta\, \right) 
&\geq& \mbP\left( 2 \, \max\limits_{\bv \in \mcV_{1/2}} \, \langle \bv, \, \bG \rangle \,\leq\, \delta \right) \s\\
&=& 1 - \mbP\left( 2 \, \max\limits_{\bv \in \mcV_{1/2}} \, \langle \bv, \, \bG \rangle \,>\, \delta \right). 
\ee
We next focus on bounding the probability 
\beno
\mbP\left( 2 \, \max\limits_{\bv \in \mcV_{1/2}} \, \langle \bv, \, \bG \rangle > \delta \right) 
&\leq& \exp(p \, \log(5)) \, \max\limits_{\bv \in \mcV_{1/2}} \, 
\mbP\left( \langle \bv, \, \bG \rangle > \dfrac{\delta}{2} \right), 
\ee
where the inequality follows from a union bound over the set of $\bv \in \mcV_{1/2}$ 
and using the fact that $\log |\mcV_{1/2}| \leq p \, \log(5)$. 
For a given $\bv \in \mcV_{1/2}$, 
Lemmas \ref{lem:exp.fam.deriv} and \ref{lem:sup_mle} allow us to write 
\beno
\langle \bv, \, \bG \rangle
\= \dsum_{i=1}^{p} \, v_i \, \left[ \nabla_{\nat_W} \, \ell_W(\truth_W, \bX_W) \right]_i 
\= \dsum_{i=1}^{p} \, v_i \, \dsum_{k=1}^{K} \, \left[s_{k,k,i}(\bX_{k,k}) - \mbE \, s_{k,k,i}(\bX_{k,k}) \right]. 
\ee
\hide{
By the local dependence assumption, 
the collection of random variables 
%(for fixed $\bv \in \mcV_{1/2}$)
\be
\label{eq:collection}
\dsum_{i=1}^{p} \, v_i \, \left[s_{k,k,i}(\bX_{k,k}) - \mbE \, s_{k,k,i}(\bX_{k,k}) \right], 
&& k \in \{1, \ldots, K\},  
\ee
is a collection of independent random variables for fixed $\bv \in \mcV_{1/2}$.
}
Observe the following two key facts: 
%This yields two key facts about $\langle \bv, \, \bG \rangle$: 
\ben
\item (Mean zero) The sum of random variables $\langle \bv, \, \bG \rangle$ satisfies $\mbE \, \langle \bv, \, \bG \rangle = 0$.\s 
\item (Sum of independent random variables) 
The sum of random variables 
\beno
\langle \bv, \, \bG \rangle
\= \dsum_{k=1}^{K} \, \left[ 
\dsum_{i=1}^{p} \, v_i \, \left[s_{k,k,i}(\bX_{k,k}) - \mbE \, s_{k,k,i}(\bX_{k,k}) \right]
\right]
\ee
is a sum of mean zero  independent random variables 
for fixed $\bv \in \mcV_{1/2}$,
because,  
by the local dependence assumption, 
the collection of random variables 
\be
\label{eq:collection}
\dsum_{i=1}^{p} \, v_i \, \left[s_{k,k,i}(\bX_{k,k}) - \mbE \, s_{k,k,i}(\bX_{k,k}) \right], 
&& k \in \{1, \ldots, K\},  
\ee
is a collection of independent random variables for fixed $\bv \in \mcV_{1/2}$.

\een
Together, 
these two points ensure the assumptions of Bernstein's inequality are met
\citep[e.g., Theorem 2.8.4, p. 35,][]{Vershynin2018}.  
Along this path, 
we first evaluate the variance term by writing    
\beno
\var \, \langle \bv, \, \bG \rangle
\= \dsum_{k=1}^{K} \, \var\left(\dsum_{i=1}^{p} \, v_i \, 
\left[s_{k,k,i}(\bX_{k,k}) - \mbE \, s_{k,k,i}(\bX_{k,k}) \right] \right) \s\\
\= \dsum_{k=1}^{K} \, \dsum_{i=1}^{p} \, \dsum_{j=1}^{p} \, \cov(v_i \, s_{k,k,i}(\bX_{k,k}), \; v_j \, s_{k,k,j}(\bX_{k,k})) \s\\
\= \dsum_{k=1}^{K} \, \dsum_{i=1}^{p} \, \dsum_{j=1}^{p} \, v_i \, v_j \, \cov(s_{k,k,i}(\bX_{k,k}), \, s_{k,k,j}(\bX_{k,k})) \s\\
\= \dsum_{k=1}^{K} \, \left\langle\bv, \, \mcI_{k,k}(\truth) \, \bv \right\rangle \s\\
\= \left\langle\bv, \, \mcI_{W}(\truth) \, \bv \right\rangle \s\\ 
&\leq& \lambda_{\max}(\mcI_W(\truth)) \, \norm{\bv}_2^2 \s\\
&\leq& \dfrac{9}{4} \, \lambda_{\max}(\mcI_W(\truth)), 
\ee
where $\lambda_{\max}(\mcI_W(\truth))$ is the largest eigenvalue of $\mcI_W(\truth)$,  
and using the inequality 
\beno
\norm{\bv}_2 
&\leq& \norm{\bu}_2 + \norm{\bu - \bv}_2
&\leq& 1 + \dfrac{1}{2}
&\leq& \dfrac{3}{2}, 
\ee
where the construction of the $\epsilon$-net $\mcV_{1/2}$ of $\mcU$ with $\epsilon = 1/2$ ensures the existence of such a $\bu \in \mcU$ 
to make the above inequality valid. 
This yields the final inequality 
\beno
\max\limits_{\bv \in \mcV_{1/2}} \, \var \, \langle \bv, \, \bG \rangle
&\leq& \dfrac{9}{4} \,  \lambda_{\max}(\mcI_W(\truth)) 
\= \dfrac{9}{4} \, K \, \maxw, 
%&\leq& \dfrac{9}{4} \, N \, \maxw, 
\ee 
defining 
\beno
\maxw
&\coloneqq& \dfrac{\lambda_{\max}(-\mbE \, \nabla_{\btheta_W}^2 \, \ell(\truth, \bX))}{K} 
\= \dfrac{\lambda_{\max}(\mcI_W(\truth))}{K}. 
\ee
We next bound the absolute value of each random variable in \eqref{eq:collection} $\mbP$-almost surely. 
By Assumption 1, 
there exists a constant $C_W > 0$,
independent of $N$, $p$, and $q$, 
such that 
\be
\label{eq:bound_on_sk_22}
\max\limits_{\bx_{k,k} \in \mbX_{k,k}} \, \norm{s_{k,k}(\bx_{k,k})}_{\infty}
&\leq& C_W \, \dbinom{|\mA_k|}{2}, 
&& k \in \{1, \ldots, K\}.  
\ee
Hence,
by the Cauchy-Schwarz inequality and using \eqref{eq:bound_on_sk_22},  
\beno
&& \max\limits_{k \in \{1, \ldots, K\}} \, 
\sup\limits_{\bx_{k,k} \in \mbX_{k,k}} \, 
\left|\dsum_{i=1}^{p} \, v_i \, \left(s_{k,k,i}(\bx_{k,k}) - \mbE \, s_{k,k,i}(\bX_{k,k}) \right) \right| \s\\
&\leq& 
\max\limits_{k \in \{1, \ldots, K\}} \,
\sup\limits_{\bx_{k,k} \in \mbX_{k,k}} \,
\norm{\bv}_2 \, 
\norm{s_{k,k}(\bx_{k,k}) - \mbE \, s_{k,k}(\bX_{k,k})}_{2} \s\\
&\leq& \dfrac{3}{2} \, 
\left( \max\limits_{k \in \{1, \ldots, K\}} \,
\sup\limits_{\bx_{k,k} \in \mbX_{k,k}} \,
\norm{s_{k,k}(\bx_{k,k}) - \mbE \, s_{k,k}(\bX_{k,k})}_{2} \right) \s\\
&\leq& \dfrac{3 \, \sqrt{p}}{2} \, 
\left( \max\limits_{k \in \{1, \ldots, K\}} \,
\sup\limits_{\bx_{k,k} \in \mbX_{k,k}} \, 
\norm{s_{k,k}(\bx_{k,k}) - \mbE \, s_{k,k}(\bX_{k,k})}_{\infty} \right) \s\\
&\leq& 3 \, \sqrt{p} \, 
\left( \max\limits_{k \in \{1, \ldots, K\}} \,
\sup\limits_{\bx_{k,k} \in \mbX_{k,k}} \, 
\norm{s_{k,k}(\bx_{k,k})}_{\infty} \right) \s\\
&\leq& 3 \, \sqrt{p} \; C_W \, \dbinom{A_{\max}}{2} \s\\
&\leq& \dfrac{3 \, C_W}{2} \, A_{\max}^2 \, \sqrt{p},   
\ee
noting the bound $\norm{\bv}_2 \leq 3 \,/\, 2$ demonstrated above.  
With these bounds, 
we apply Bernstein's inequality (for just the upper-tail) \citep[e.g., Theorem 2.8.4, p. 35,][]{Vershynin2018}
to obtain 
\beno
\mbP\left(\langle \bv, \, \bG \rangle \,>\, \delta \right)
&\leq& \exp\left( -\dfrac{\delta^2 \,/\, 2}{(9 \,/\, 4) \, K \, \maxw  
+ (3 \, C_W \,/\, 2) \, A_{\max}^2 \, \sqrt{p} \, \delta \,/\, 3} \right) \s\\
&\leq& \exp\left( - \dfrac{\delta^2}{5 \, K \, \maxw 
+ C_W A_{\max}^2 \, \sqrt{p} \, \delta} \right).  
\ee
Collecting results, 
we have shown, 
for $\delta > 0$,
that 
\beno
\mbP\left( \norm{\nabla_{\btheta_W} \, \ell_W(\truth_W, \bX_W)}_2 \,\leq\, \delta \right)
&\geq& 1 - \exp\left( - \dfrac{\delta^2}{5 \, K \, \maxw
+ C_W A_{\max}^2 \, \sqrt{p} \, \delta} + \log(5) \, p \right). 
\ee

\s

{\bf Changes for the between-bock case.} 
By a similar argument, 
\beno
\mbP\left( \norm{\nabla_{\btheta_B} \, \ell_B(\truth_B, \bX_B)}_2 \,\leq\, \delta \right)
&\geq& 1 - \exp\left( - \dfrac{\delta^2}{5 \, \tbinom{K}{2} \, \maxb
+ 2 \, C_B A_{\max}^2 \, \sqrt{q} \, \delta} + \log(5) \, q \right),
\ee
for $\delta > 0$. 
We highlight the main changes to the above argument. 
First,
we now take  
\beno
\bG &\coloneqq& \nabla_{\btheta_B} \, \ell_B(\truth_B, \bX_B).
\ee
Second, 
the dimension $p$ of $\truth_W \in \mbR^p$ is replaced by the dimension $q$ of $\truth_B \in \mbR^q$. 
This implies that $\log |\mcV_{1/2}| \leq q \log(5)$. 
Third, 
\beno
\langle \bv, \, \bG \rangle 
\= \dsum_{1 \leq k < l \leq K} \, \left[ \dsum_{i=1}^{q} \, v_i \, 
\left[ s_{k,l,i}(\bX_{k,l}) - \mbE \, s_{k,l,i}(\bX_{k,l}) \right] \right], 
\ee
which implies 
\beno
\var \, \langle \bv, \, \bG \rangle
&\leq& \dfrac{9}{4} \, \lambda_{\max}(\mcI_{B}(\truth)) 
&=& \dfrac{9}{4} \, \dbinom{K}{2} \, \maxb,
%&\leq& \dfrac{9}{4} \, N^2 \, \maxb, 
\ee
defining 
\beno
\maxb
&\coloneqq& \dfrac{\lambda_{\max}(- \mbE \, \nabla_{B}^2 \, \ell(\truth, \bX))}{\tbinom{K}{2}} 
\= \dfrac{ \lambda_{\max}(\mcI_{B}(\truth))}{\tbinom{K}{2}}. 
\ee
Fourth, 
and finally,  
we have the bound 
\beno
\max\limits_{\{k,l\} \subseteq \{1, \ldots, K\}} \, 
\sup\limits_{\bx_{k,l} \in \mbX_{k,l}} \,
\left|\dsum_{i=1}^{q} \, v_i \, \left[s_{k,l,i}(\bx_{k,l}) - \mbE \, s_{k,l,i}(\bX_{k,l}) \right] \right| 
&\leq& 3 \, C_B \, A_{\max}^2 \, \sqrt{q}. 
\ee
Together, 
these changes will yield the inequality
\beno
\mbP\left( \norm{\nabla_{\btheta_B} \, \ell_B(\truth_B, \bX_B)}_2 \,\leq\, \delta \right)
\;\geq\; 1 - \exp\left( - \dfrac{\delta^2}{5 \, \tbinom{K}{2} \, \maxb 
+ 2 \, C_B A_{\max}^2 \, \sqrt{q} \, \delta} + \log(5) \, q \right),
\ee
for $\delta > 0$. 
\qed

\s\s

\s

\section{Proof of Theorem 2.2}
\label{minimax_lower}

%\ttproof \ref{thm:optimal}.
Our method of proof utilizes Fano's method for lower bounding the minimax risk.
We present a general argument for lower bounding the minimax risk for exponential families utilizing 
Fano's method and then apply the obtained general argument to our specific cases. 

\s

{\bf General argument.} 
Consider an exponential family of densities $\{f_{\nat} : \nat \in \mbR^m\}$ 
for a finite support $\mbX$
given by 
\beno
f_{\nat}(\bx) 
\= h(\bx) \, \exp\left(\langle \nat, \, s(\bx) \rangle - \psi(\nat) \right)
&>& 0, 
&& \bx \in \mbX, 
\ee
data-generating parameter vector $\truth$, 
and define the minimax risk in the $\ell_2$-norm to be 
\beno
\mcR &\coloneqq& 
\inf\limits_{\mle} \; 
\sup\limits_{\nat \in \mbR^m} \; 
\mbE_{\nat} \, \norm{\mle - \nat}_2. 
\ee
In the case that $f_{\nat}(\bx) = 0$ for some $\bx \in \mbX$, 
we would simply reduce the support to 
\beno 
\mbX_0 
&\coloneqq& \{\bx \in \mbX \,:\, f_{\nat}(\bx) > 0 \},
\ee
obtaining a family of strictly positive densities on $\mbX_0$. 
For ease of presentation, 
we therefore proceed without loss of generality assuming that 
$\{f_{\nat} : \nat \in \mbR^m\}$ are strictly positive on the support $\mbX$. 

Let $\epsilon > 0$ be fixed and consider  
$\gamma \in (0, \epsilon)$. 
Assume that $\{\nat_1, \ldots, \nat_M\} \subset \mB_2(\truth, \gamma)$ ($M \geq 2$)
is a $2 \, \delta$-separated set in the metric 
$d(\bv, \bw) = \norm{\bv - \bw}_2$,
i.e.,
$d(\nat_i, \nat_j) \geq 2 \, \delta$ for all pairs $\{i,j\} \subseteq \{1, \ldots, M\}$. 
Then, 
by Proposition 15.12 (p. 502) of \citet{Wainwright2019} and the discussions following, 
the minimax risk $\mcR$ has the lower bound 
\beno
\mcR
&\geq& 
\delta \, \left[1 - \dfrac{\mI + \log(2)}{\log \, M} \right],
\ee
where 
\be
\label{eq:inf_crit}
\mI
&\coloneqq& \max\limits_{\{i,j\} \subseteq \{1, \ldots, M\}} \; 
\text{KL}(\nat_i, \nat_j),
\ee
defining the Kullback–Leibler divergences 
\beno
\text{KL}(\nat_i, \nat_j) &\coloneqq& \dsum_{\bx \in \mbX}  \, f_{\nat_i}(\bx) \, 
\log \, \dfrac{f_{\nat_i}(\bx)}{f_{\nat_j}(\bx)},
&& \{i,j\} \subseteq \{1, \ldots, M\}.  
\ee
For an exponential family, 
we can express $\text{KL}(\nat_i, \nat_j)$ as 
\be
\label{eq:express}
\text{KL}(\nat_i, \nat_j) 
\= \dsum_{\bx \in \mbX} \, 
f_{\nat_i}(\bx) 
\left[ \log h(\bx) - \log h (\bx) + \langle \nat_i - \nat_j, \, s(\bx) \rangle - \psi(\nat_i) + \psi(\nat_j) \right] \s\\
\= \mbE_{\nat_i}\langle \nat_i - \nat_j, \, s(\bX) \rangle 
- \psi(\nat_i) + \psi(\nat_j) \s\\
\= \langle \nat_i - \nat_j, \, \bmu(\nat_i) \rangle  - \psi(\nat_i) + \psi(\nat_j), 
\ee
defining $\bmu(\nat) \coloneqq \mbE_{\nat} \, s(\bX)$ to be the mean-value parameter map
of the exponential family.  
Using Corollary 2.3 of \citet{Brown1986}, 
we perform the expansion 
\be
\label{eq:expansion_psi}
\psi(\nat_j) 
\= \psi(\nat_i) + \langle \nat_j - \nat_i, \, \bmu(\nat_i) \rangle 
+ \dfrac{1}{2} \, \langle \nat_j - \nat_i, \, \mcI(\dot\nat) \, (\nat_j - \nat_i) \rangle \s\\
\=  \psi(\nat_i) - \langle \nat_i - \nat_j,  \, \bmu(\nat_i) \rangle
+ \dfrac{1}{2} \, \langle \nat_i - \nat_j, \, \mcI(\dot\nat) \, (\nat_i - \nat_j) \rangle,  
\ee
where $\dot\nat = t \, \nat_i + (1 - t) \, \nat_j$ (for some $t \in (0, 1)$) 
and 
$\mcI(\nat) \coloneqq -\mbE \,\nabla_{\nat}^2 \, \log f_{\nat}(\bX)$ is the Fisher information matrix
corresponding to $f_{\nat}$.  
Combining \eqref{eq:express} and \eqref{eq:expansion_psi}, 
\beno
\text{KL}(\nat_i, \nat_j)
\= \dfrac{1}{2} \, \langle \nat_i - \nat_j, \, \mcI(\dot\nat) \, (\nat_i - \nat_j) \rangle \s\\
&\leq& \dfrac{1}{2} \, n \, \widetilde{\lambda}_{\max}^{\epsilon} \, \norm{\nat_i - \nat_j}_2^2 \s\\
&\leq& \dfrac{1}{2} \, n \, \widetilde{\lambda}_{\max}^{\epsilon} \, 
\left( \norm{\nat_i - \truth}_2 + \norm{\nat_j - \truth}_2\right)^2 \s\\
&\leq& 2 \, n \, \epsilon^2 \, \widetilde{\lambda}_{\max}^{\epsilon}, 
\ee
defining for a fixed $\epsilon > 0$ the quantity 
\beno
\widetilde{\lambda}_{\max}^{\epsilon}
&\coloneqq& \sup\limits_{\nat \in \mB_2(\truth, \epsilon)} \, \dfrac{\lambda_{\max}(\mcI(\nat))}{n}, 
\ee
noting that 
$\{\nat_1, \ldots, \nat_M\} \subset \mB_2(\truth, \gamma) \subset \mB_2(\truth, \epsilon)$ by assumption.  
The size $M$ of the largest possible $2 \, \delta$-separated set 
$\{\nat_1, \ldots, \nat_M\} \subset \mB_2(\truth, \gamma) \subset \mbR^m$ 
is the packing number of $\mB_2(\truth, \gamma)$,
which by Lemma 4.2.8 (equivalence of covering and packing numbers) 
and Corollary 4.2.13 (covering numbers of the Euclidean ball) of \citet{Vershynin2018}, 
satisfies 
\beno
M &\geq& \left(\dfrac{\gamma}{2 \, \delta} \right)^m,
\ee
taking the $2 \, \delta$-separated set $\{\nat_1, \ldots, \nat_M\} \subset \mB_2(\truth, \gamma)$ to be as large as possible
and applying the results to a Euclidean ball of arbitrary radius $\gamma > 0$.  
As a result, 
\beno
\log  M
&\geq& m \, \log(\gamma \,/\, 2 \, \delta). 
\ee
Altogether, 
we have demonstrated the bound 
\beno
\mcR
&\geq&
\delta \, \left[1 - \dfrac{2 \, n \, \gamma^2 \, \widetilde{\lambda}_{\max}^{\epsilon} + \log(2)}{m \, \log(\gamma \,/\, 2 \, \delta)} \right]. 
\ee
We desire that 
\beno
\dfrac{2 \, n \, \gamma^2 \, \widetilde{\lambda}_{\max}^{\epsilon} + \log(2)}{m \, \log(\gamma \,/\, 2 \, \delta)}
&\leq& \dfrac{1}{2}, 
\ee
in order to show that $\mcR \geq \delta \,/\, 2$. 
Re-arranging this inequality, 
we have 
\beno
\dfrac{4 \, n \, \gamma^2 \, \widetilde{\lambda}_{\max}^{\epsilon}}{m}
&\leq& 
\dfrac{4 \, n \, \gamma^2 \, \widetilde{\lambda}_{\max}^{\epsilon}}{m} + \dfrac{2 \, \log(2)}{m}
&\leq& \log(\gamma \,/\, 2) - \log(\delta), 
\ee 
and exponentiating we obtain 
\beno
\exp\left(\dfrac{4 \, n \, \gamma^2 \, \widetilde{\lambda}_{\max}^{\epsilon}}{m} \right) \; 
%\exp\left( \dfrac{2 \, \log(2)}{m} \right)
&\leq& \dfrac{\gamma \,/\, 2}{\delta}
&\leq& \dfrac{\gamma}{\delta}.
\ee
This leads us to the following inequality 
\beno
\delta &\leq& 
\gamma \; \exp\left(-\dfrac{4 \, n \, \gamma^2 \, \widetilde{\lambda}_{\max}^{\epsilon}}{m} \right).  
\ee
Choosing 
\beno
\gamma 
\=  C \; 
\sqrt{\dfrac{m}{n \, \widetilde{\lambda}_{\max}^{\epsilon}}},
\ee
for some $C > 0$ which is presumed to be fixed, but freely chosen, 
yields the bound 
\be
\label{eq:big_bound}
\delta
&\leq& 
C \; \exp(-4 \, C^2) \, 
\sqrt{\dfrac{m}{n \, \widetilde{\lambda}_{\max}^{\epsilon}}}.  
\ee
As long as $m = O(n \, \widetilde{\lambda}_{\max}^{\epsilon})$,
we can choose $C > 0$ to ensure that $\gamma \in (0, \epsilon)$.  
Thus, 
for all $\delta > 0$ satisfying \eqref{eq:big_bound},
we can lower bound the minimax risk $\mcR$ by 
\beno
\mcR 
&\geq& \dfrac{\delta}{2}. 
\ee
The remainder of the proof will utilize this general argument 
to lower bound the minimax risk in the $\ell_2$-norm for exponential-family local dependence random graphs.

\s\s

{\bf Lower bounds to the minimax risk in the $\ell_2$-norm for exponential-family local dependence random graphs.} 
We will first handle the within-block case by considering 
\beno
\mcR_{W,N}
&\coloneqq& \inf\limits_{\mle_W} \; \sup\limits_{\nat \in \mbR^{p+q}} \; 
\mbE_{\nat} \, \norm{\mle_W - \nat_W}_2. 
\ee 
Fix $\epsilon > 0$, 
independent of $N$, $p$, and $q$, 
and define 
\beno
\maxwe 
&\coloneqq& \sup\limits_{\nat \in \mB_2(\truth, \epsilon)} \, 
\dfrac{\lambda_{\max}\left(- \mbE \, \nabla_{\nat_W}^2 \, \ell(\nat, \bX)\right)}{K}. 
\ee
With this definition,
we revisit \eqref{eq:big_bound} 
taking $m = p$ and $n = K$ to obtain 
\beno
\mcR_{W,N}
&\geq& C_1 \; \exp(-4 \, C_1^2)  \; \sqrt{\dfrac{p}{K \, \maxwe}}, 
\ee
for some $C_1 > 0$ assumed to be fixed, but freely chosen.
Using the relationship 
\beno
N 
\= \dsum_{k=1}^{K} \, |\mA_k|
\= K \, \dfrac{1}{K} \dsum_{k=1}^{K} \, |\mA_k|
\= K \, \Aavg,
\ee
defining $\Aavg \coloneqq K^{-1} \, \sum_{k=1}^{K} |\mA_k|$, 
we obtain 
\beno
\mcR_{W,N}
&\geq& C_1 \; \exp(-4 \, C_1^2) \; \sqrt{\dfrac{\Aavg}{\maxwe}} \; \sqrt{\dfrac{p}{N}}.  
\ee
Then there exists $B_1 \coloneqq C_1\exp(-4 \, C_1^2) > 0$,
independent of $N$, $p$, and $q$, 
such that 
\beno
\mcR_{W,N}
&\geq& B_1 \; \sqrt{\dfrac{\Aavg}{\maxwe}} \; \sqrt{\dfrac{p}{N}} \s\\
\= B_1 \; \dfrac{1}{\sqrt{\maxwe}} \; \sqrt{\dfrac{\maxw}{\maxw}} \; \dfrac{\minw}{\minw} \; 
\sqrt{\Aavg} \; \sqrt{\dfrac{p}{N}} \s\\ 
\= B_1 \; \left(\dfrac{\minw}{ \sqrt{\maxw \; \maxwe}} \right) \; \dfrac{\sqrt{\maxw}}{\minw} \; \sqrt{\Aavg} \; \sqrt{\dfrac{p}{N}} \s\\
&\geq& B_1 \; \left(\dfrac{\minw}{\maxwe} \right) \; \dfrac{\sqrt{\maxw}}{\minw} \; \sqrt{\Aavg} \; \sqrt{\dfrac{p}{N}}. 
\ee
%One can show $B_1 \leq 1 \,/\, 4 >0$. 
The above inequality establishes 
\beno
\mcR_{W,N}
&\geq& B_1 \; \sqrt{\dfrac{\Aavg}{\maxwe}} \; \sqrt{\dfrac{p}{N}} 
&\geq& B_1 \; \left(\dfrac{\minw}{\maxwe} \right) \; \dfrac{\sqrt{\maxw}}{\minw} \; \sqrt{\Aavg} \; \sqrt{\dfrac{p}{N}}.
\ee
Next,
we prove the between-block case and  
consider 
\beno
\mcR_{B,N}
&\coloneqq& \inf\limits_{\mle_B} \; \sup\limits_{\nat \in \mbR^{p+q}} \; 
\mbE_{\nat} \, \norm{\mle_B - \nat_B}_2. 
\ee 
Fix $\epsilon > 0$, 
independent of $N$, $p$, and $q$, 
and define 
\beno
\maxbe 
&\coloneqq& \sup\limits_{\nat \in \mB_2(\truth, \epsilon)} \, 
\dfrac{\lambda_{\max}\left(- \mbE \, \nabla_{\nat_B}^2 \, \ell(\nat, \bX)\right)}{\tbinom{K}{2}}. 
\ee
With this definition,
we revisit \eqref{eq:big_bound} 
taking $m = q$ and $n = \tbinom{K}{2}$ to obtain 
\beno
\mcR_{B,N}
&\geq& C_2 \; \exp(-4 \, C_2^2) \; 
\sqrt{\dfrac{q}{\tbinom{K}{2} \, \maxbe}}, 
\ee
for some $C_2 > 0$ assumed to be fixed, but freely chosen.
Using the relationship 
\beno
N 
\= \dsum_{k=1}^{K} \, |\mA_k|
\= K \, \dfrac{1}{K} \dsum_{k=1}^{K} \, |\mA_k|
\= K \, \Aavg,
\ee
defining $\Aavg \coloneqq K^{-1} \, \sum_{k=1}^{K} |\mA_k|$, 
we obtain 
\beno
\mcR_{B,N}
&\geq& C_2 \; \exp(-4 \, C_2^2) \; \dfrac{1}{\sqrt{\maxbe}} \; \sqrt{\dfrac{q}{\binom{N \,/\, \Aavg}{2}}} \s\\
&\geq& C_2 \; \exp(-4 \, C_2^2) \; \dfrac{1}{\sqrt{\maxbe}} \; \sqrt{\dfrac{2 \, q \, \Aavg^2}{N^2}} 
\= B_2 \; \dfrac{\Aavg}{\sqrt{\maxbe}} \; \sqrt{\dfrac{q}{N^2}}, 
\ee
where $B_2 \coloneqq  \sqrt{2} \, C_2 \; \exp(-4 \, C_2^2) > 0$ is independent of $N$, $p$, and $q$. 
Hence, 
\beno
\mcR_{B,N}
&\geq& B_2 \; \dfrac{\Aavg}{\sqrt{\maxwe}} \; \sqrt{\dfrac{q}{N^2}} \s\\
\= B_2 \; \dfrac{1}{\sqrt{\maxwe}} \; \sqrt{\dfrac{\maxw}{\maxw}} \; \dfrac{\minw}{\minw} \; 
\Aavg \; \sqrt{\dfrac{q}{N^2}} \s\\ 
\= B_2 \; \left(\dfrac{\minw}{ \sqrt{\maxw \; \maxwe}} \right) \; \dfrac{\sqrt{\maxw}}{\minw} \; \Aavg \; \sqrt{\dfrac{q}{N^2}} \s\\
&\geq& B_2 \; \left(\dfrac{\minw}{\maxwe} \right) \; \dfrac{\sqrt{\maxw}}{\minw} \; \Aavg \; \sqrt{\dfrac{q}{N^2}}, 
\ee
showing the claimed result of 
\beno
\mcR_{B,N}
&\geq& B_2 \; \dfrac{\Aavg}{\sqrt{\maxwe}} \; \sqrt{\dfrac{q}{N^2}}
&\geq& B_2 \; \left(\dfrac{\minw}{\maxwe} \right) \; \dfrac{\sqrt{\maxw}}{\minw} \; \Aavg \; \sqrt{\dfrac{q}{N^2}}. 
\ee

Lastly, 
the assumption in the general argument that $m = O(n \, \widetilde{\lambda}_{\max}^{\epsilon})$
requires:
\bi \item $p = \dim(\truth_W)$ 
satisfies $p = O(N \, \maxwe)$, and  
\item $q = \dim(\truth_B)$ satisfies $q = O(N^2 \, \maxbe)$,
\ei 
substituting the relevant quantities into $m = O(n \, \widetilde{\lambda}_{\max}^{\epsilon})$
for each case. 
\qed

\s 

\section{Proof of Theorem 2.3}
\label{minimax_upper}

We start by considering the restricted minimax risk 
\beno
\inf\limits_{\mle_W} \, 
\sup\limits_{\nat \in \mB_2(\truth, \epsilon)} \, \mbE_{\nat} \, \norm{\mle_W - \nat_W}_2,  
%&\leq& \sup\limits_{\nat \in \mB_2(\truth, \epsilon)} \, \mbE_{\nat} \, \norm{\widetilde\btheta_W - \nat_W}_2,
\ee
where $\epsilon > 0$ is the same as in Assumptions 2 and 3. 
Let $\nat \in \mB_2(\truth, \epsilon) \subset \mbR^{p+q}$ be arbitrary. 
We partition the support $\mbX$ of $\bX$ as follows: 
\beno
\mbX_1 &\coloneqq&
\left\{ 
\bx \in \mbX \,:\, \norm{\mle_W - \truth_W}_2 \,\leq\, C_1 \; \sqrt{A_{\avg}} \; \dfrac{\sqrt{\maxw}}{\minw} \, \sqrt{\dfrac{p}{N}}
\right\} \s\\
\mbX_2 &\coloneqq& \mbX \;\setminus\; \mbX_1, 
\ee
where the constant $C_1 > 0$ is the same as the one guaranteed in Theorem 2.1.
Then, 
\be
\label{eq:upper_inq1}
\mbE_{\nat} \, \norm{\mle_W - \nat_W}_2
\= \mbE_{\nat}\left[\norm{\mle_W - \nat_W}_2 \,|\, \bX \in \mbX_1\right] \, \mbP_{\nat}(\bX \in \mbX_1) \s\\ 
&+& \mbE_{\nat}\left[\norm{\mle_W - \nat_W}_2 \,|\, \bX \in \mbX_2\right] \, \mbP_{\nat}(\bX \in \mbX_2),
\ee
by the law of total expectation,
and  
\beno
\inf\limits_{\mle_W} \,
\sup\limits_{\nat \in \mB_2(\truth, \epsilon)} \, \mbE_{\nat} \, \norm{\mle_W - \nat_W}_2 
&\leq& 
%\inf\limits_{\mle_W \in \mbR^p} \,
\sup\limits_{\nat \in \mB_2(\truth, \epsilon)} \, 
\mbE_{\nat}\left[\norm{\mle_W - \nat_W}_2 \,|\, \bX \in \mbX_1\right] \, \mbP_{\nat}(\bX \in \mbX_1) \s\\
&+& 
%\inf\limits_{\mle_W \in \mbR^p} \,
\sup\limits_{\nat \in \mB_2(\truth, \epsilon)} \, \mbE_{\nat}\left[\norm{\widetilde\btheta_W - \nat_W}_2 \,|\, \bX \in \mbX_2\right] \, \mbP_{\nat}(\bX \in \mbX_2), 
\ee
where $\mle_W$ in the first term in the upper bound is the maximum likelihood estimator
and 
\beno
\widetilde\btheta_W
&\coloneqq& \argmax\limits_{\btheta_W \in \mB_2(\bm{0}_p, 100)} \, \ell_W(\nat_W, \bX_W)
\ee
is the maximum likelihood estimator restricted to the subset $\mB_2(\bm{0}_p, 100)$. 
Note that Theorem 2.1 establishes that there exists $N_0 \geq 1$,
independent of $N$, $p$, or $q$, 
such that 
$\mbP(\bX \in \mbX_1) \geq 1 - N^{-2}$,
in which case we have the bound $\mbP(\bX \in \mbX_2) \leq N^{-2}$, 
for all integers $N \geq N_0$. 
Observe that the upper bound becomes trivial in the case when $\mbP(\bX \in \mbX_2) = 0$. 
We then obtain the upper bound 
\beno
\sup\limits_{\nat \in \mB_2(\truth, \epsilon)} \, 
\mbE_{\nat}\left[\norm{\widetilde\btheta_W - \nat_W}_2 \,|\, \bX \in \mbX_2\right] 
&\leq& 
\sup\limits_{\widetilde\btheta_W \in \mB_2(\bm{0}_p, 100)} \,
\sup\limits_{\nat \in \mB_2(\truth, \epsilon)} \, 
\norm{\widetilde\btheta_W - \nat_W}_2.
\ee
Define 
\beno
M &\coloneqq& \sup\limits_{\widetilde\btheta_W \in \mB_2(\bm{0}_p, 100)} \,
\sup\limits_{\nat \in \mB_2(\truth, \epsilon)} \,
\norm{\widetilde\btheta_W - \nat_W}_2
&\in& (0, \infty). 
\ee
Continuing from \eqref{eq:upper_inq1} leads us to the upper bound   
\beno
\mbE_{\nat} \, \norm{\mle_W - \nat_W}_2
&\leq& \mbE_{\nat}\left[\norm{\mle_W - \nat_W}_2 \,|\, \bX \in \mbX_1\right] \, \mbP_{\nat}(\bX \in \mbX_1) 
+ M \;  \mbP_{\nat}(\bX \in \mbX_2) \s\\ 
&\leq& \mbE_{\nat}\left[\norm{\mle_W - \nat_W}_2 \,|\, \bX \in \mbX_1\right] 
\left(1 - \dfrac{1}{N^2} \right)
+ M\left( \dfrac{1}{N^2} \right) \s\\ 
&\leq& \mbE_{\nat}\left[\norm{\mle_W - \nat_W}_2 \,|\, \bX \in \mbX_1\right] 
+ \dfrac{M}{N^2} \s\\ 
&\leq& C_1 \; \sqrt{A_{\avg}} \; \dfrac{\sqrt{\maxw}}{\minw} \, \sqrt{\dfrac{p}{N}}
+ \dfrac{M}{N^2}. 
\ee
Note that Assumption 2 implies that 
\beno
\lim\limits_{N \to \infty} \; 
\sqrt{A_{\avg}} \; \dfrac{\sqrt{\maxw}}{\minw} \, \sqrt{\dfrac{p}{N}}
\= 0, 
\ee 
which further implies 
\beno
\sqrt{A_{\avg}} \; \dfrac{\sqrt{\maxw}}{\minw} \, \sqrt{p} 
\= o\left(\sqrt{N}\,\right),
\ee
in turn ultimately implying 
\beno
\sqrt{A_{\avg}} \; \dfrac{\sqrt{\maxw}}{\minw} \, \sqrt{\dfrac{p}{N}}
\= o\left(N^2\right). 
\ee
As a result, 
there exists an $N_1 \geq N_0$,
independent of $N$, $p$, and $q$, 
such that 
\beno
\sqrt{A_{\avg}} \; \dfrac{\sqrt{\maxw}}{\minw} \, \sqrt{\dfrac{p}{N}} 
&\geq& \dfrac{1}{N^2},
&& N \geq N_1.  
\ee
Defining $C_2 \coloneqq C_1 + M > 0$,
also independent of $N$, $p$, and $q$, 
we have 
\beno 
\inf\limits_{\mle_W} \, 
\sup\limits_{\nat \in \mB_2(\truth, \epsilon)} \,
\mbE_{\nat} \, \norm{\mle_W - \nat_W}_2
&\leq& C_2 \;  \sqrt{A_{\avg}} \; \dfrac{\sqrt{\maxw}}{\minw} \, \sqrt{\dfrac{p}{N}}.  
\ee
Repeating the same essential argument for 
\beno
 \inf\limits_{\mle_B} \,
\sup\limits_{\nat \in \mB_2(\truth, \epsilon)} \, \mbE_{\nat} \, \norm{\mle_B - \nat_B}_2
\ee
we obtain, 
for $C_3 > 0$ defined similarly to $C_2 > 0$, 
\beno
 \inf\limits_{\mle_B} \,
\sup\limits_{\nat \in \mB_2(\truth, \epsilon)} \, \mbE_{\nat} \, \norm{\mle_B - \nat_B}_2
&\leq& C_3 \; \Aavg \; \dfrac{\sqrt{\maxb}}{\minb} \, \sqrt{\dfrac{q}{N^2}}.  
\ee
\qed

\s

\section{Proof of Corollary 2.4} 
\label{cor_minimax}

%\ccproof \ref{cor:optimal}.
The assumptions of both Theorems 2.1 and 2.2 are met. 
Theorem 2.1 supplies the following upper bounds to the $\ell_2$-error of maximum likelihood estimators: 
\beno
\norm{\mle_W - \truth_W}_2
&\leq& C \; \sqrt{A_{\avg}} \; \dfrac{\sqrt{\maxw}}{\minw} \, \sqrt{\dfrac{p}{N}}
\ee
\vspace{.05cm}
\beno
\norm{\mle_B - \truth_B}_2
&\leq& C \;  A_{\avg} \; \dfrac{\sqrt{\maxb}}{ \minb} \, \sqrt{\dfrac{q}{N^2}}.
\ee
Theorem 2.2 provides the lower bounds 
on the minimax risks $\mcR_{W,N}$ and $\mcR_{B,N}$:  
\beno
\mcR_{W,N}
&\geq& B_1  \left(\dfrac{\minw}{\maxwe} \right)
\, \dfrac{\sqrt{\maxw}}{\minw} \; \sqrt{\Aavg} \; \sqrt{\dfrac{p}{N}}
\ee
\vspace{.05cm}
\beno
\mcR_{B,N}
&\geq& B_2  \left(\dfrac{\minw}{\maxwe} \right)
\, \dfrac{\sqrt{\maxw}}{\minw} \; \Aavg \; \sqrt{\dfrac{q}{N^2}}.
\ee
Inspecting the two sets of bounds, 
the assumption that 
\beno
\maxwe \;=\; O\left( \minw \right)
&&\mbox{and}&&
\maxbe \;=\; O\left( \minb \right)
\ee
ensures that the two sets of bounds match (up to an unknown, but fixed, constant). 
As a result, 
the upper bounds on the $\ell_2$-error presented in Theorem 2.1 achieve (up to an unknown, but fixed, constant)
the minimax rate of convergence. 
\qed

\s

\section{Proof of Theorem 2.5}
\label{norm}

For ease of presentation, 
we first present a general argument for bounding the error of the  
multivariate normal approximation,
and then show how it can be applied to
maximum likelihood estimators of exponential-families of local dependence random graph models
in order to establish the desired result. 

\s

\noindent {\bf Bounding the error of the multivariate normal approximation.} 
Consider a general estimating function $m : \mbR^d \times \mbX \mapsto \mbR$
which admits the following form:
\beno
m(\nat, \bx)
\= \dsum_{k=1}^{K} \, m_{k,k}(\nat_W, \bx_{k,k})
+ \dsum_{1 \leq k < l \leq  K} \, m_{k,l}(\nat_B, \bx_{k,l}),
\ee
and assume $m(\nat,\bx)$ is thrice continuously differentiable in the elements of $\nat \in \mbR^{p+q}$.  
By assumption,
the subgraphs $\bX_{k,l}$ ($1 \leq k \leq l \leq K$)
are mutually independent,
implying, for a fixed $\nat \in \mbR^{p+q}$,
that 
the collection of random variables $m_{k,k}(\nat_W, \bX_{k,k})$ ($1 \leq k \leq K$) 
and $m_{k,l}(\nat_B, \bX_{k,l})$ ($1 \leq k < l \leq K$)
are likewise mutually independent.
As such, 
\be
\label{eq:gen_grad_m}
\nabla_{\nat} \, m(\nat, \bx)
\= \dsum_{k=1}^{K} \, \nabla_{\nat} \,  m_{k,k}(\nat_W, \bx_{k,k})
+ \dsum_{1 \leq k < l \leq  K} \, \nabla_{\nat} \,  m_{k,l}(\nat_B, \bx_{k,l})
\ee
is a sum of mutually independent random vectors.
Assume that
\be
\label{eq:zeros}
\mbE \, \nabla_{\nat} \,  m_{k,k}(\truth_W, \bX_{k,k})
\= \bm{0}_{p+q},
&& 1 \leq k \leq K \s\\ 
\mbE \, \nabla_{\nat} \,  m_{k,l}(\truth_B, \bX_{k,l})
\= \bm{0}_{p+q},
&& 1 \leq k < l \leq K,
\ee
implying that 
$\mbE \, \nabla_{\nat} \, m(\truth, \bX) = \bm{0}_{p+q}$.
 
Let $\nat \in \mbR^{p+q}$ and $\bx \in \mbX$ be fixed.  
By a multivariate Taylor expansion,
\be
\label{eq:gen_mvte}
\nabla_{\nat} \, m(\nat, \bX)
\= \nabla_{\nat} \, m(\truth, \bX)
+ \nabla_{\nat}^2 \, m(\truth, \bX) \, (\nat - \truth) + \bR,
\ee
where $\bR \in \mbR^{p+q}$ is a vector of remainders given in the Lagrange form,
where each of the 
remainder terms $R_i$ ($i \in \{1, \ldots, p+q\}$) is given by
\be
\label{eq:remainder}
R_i
\=  \dsum_{j=1}^{p+q} \, \dfrac{1}{2} \,
\left[ \, \dfrac{\partial^2}{\partial \, \theta_j^2} \, 
\left[\nabla_{\nat} \, m(\dot\nat^{(i)}, \bX) \right]_i \right] \,
(\theta_j - \theta_j^\star)^2 \s \\
&& +  \dsum_{1 \le j < r \le p+q} \,\dfrac{1}{2} \,  
\left[ \, \dfrac{\partial^2}{\partial \, \theta_j \; \partial \, \theta_r} \, 
\left[ \nabla_{\nat} \, m(\dot\nat^{(i)}, \bX) \right]_i \right] \,
(\theta_j - \theta^\star_j) \, (\theta_r - \theta^\star_r),
\ee
where $\dot\nat^{(i)} = t_i \, \nat + (1 - t_i) \, \truth$ (for some $t_i \in [0, 1]$, $i \in \{1, \ldots, p+q\}$).
Assume that 
$\bC \coloneqq \var \, \nabla_{\nat} \, m(\truth, \bX)$ is non-singular and that 
$\nabla_{\nat} \,  m(\nat, \bx)$ has a root given by $\nat_0 \in \mbR^{p+q}$.
Taking $\nat = \nat_0$, 
we re-arrange \eqref{eq:gen_mvte} with the observation $\bX = \bx$ in order to obtain
\be
\label{norm_key}
\nabla_{\nat} \, m(\truth, \bx)
\= \nabla_{\nat}^{2} \, m(\truth, \bx) \, (\truth - \nat_0) - \bR. 
\ee
Bear in mind, from the form of \eqref{eq:gen_grad_m},
that $\nabla_{\nat} \, m(\truth, \bx)$ is a sum of independent random vectors.
Define
\beno
\bY_{k,l} &\coloneqq& \begin{cases}
\bC^{-1/2} \, \; \nabla_{\nat} \, m_{k,k}(\truth_W, \bX_{k,k}),
& \mbox{if } k = l \s\\ 
\bC^{-1/2} \, \; \nabla_{\nat} \, m_{k,l}(\truth_B, \bX_{k,l}), 
& \mbox{if } k \neq l
\end{cases}, 
&&& 1 \leq k \leq l \leq K,
\ee
and 
\beno 
\bS 
&\coloneqq& \dsum_{1 \leq k \leq l \leq K} \, \bY_{k,l}.
\ee
Observe that,
by \eqref{eq:zeros},
$\mbE \, \bS = \bm{0}_{p+q}$ 
and that,
by the definition of $\bC$, 
$\var \, \bS = \bI_{p+q}$,
where $\bI_{p+q}$ is the $(p+q)$-dimensional identity matrix. 
Applying Lemma \ref{lem:raic}, 
for all measurable convex sets $\mcC \subset \mbR^{p+q}$,
\beno
|\mbP(\bS \in \mcC) - \Phi_{p+q}(\bZ_{p+q} \in \mcC)|
&\leq& (42 \, (p+q)^{1/4} + 16) \, \dsum_{1 \leq k \leq l \leq K} \, \mbE \, \norm{\bY_{k,l}}_2^3 \s\\
&\leq& 58 \, (p+q)^{1/4} \,  \dsum_{1 \leq k \leq l \leq K} \, \mbE \, \norm{\bY_{k,l}}_2^3. 
\ee
Normality results for $\bS$ can be extended to a standardization of $(\truth - \nat_0)$ via \eqref{norm_key}:
\be
\label{eq:norm_key}
\bS
&\overset{D}{=}&
\bC^{-1/2} \, \left[ \nabla_{\nat}^{2} \, m(\truth, \bx) \, (\truth - \nat_0) - \bR \right], 
\ee
where $\overset{D}{=}$ indicates equality in distribution. 
%In what follows,
%we will utilize the above argument,
%together with additional auxiliary results in the supplementary materials to establish the multivariate normality
%of maximum likelihood estimators of exponential-family local dependence random graph models. 

\s\s

\noindent {\bf Multivariate normal approximation for maximum likelihood estimators.}
Define
\beno
m_{k,k}(\nat_W, \bx_{k,k})
&\coloneqq& \ell_{k,k}(\nat_W, \bx_{k,k}),
&& 1 \leq k \leq K \s\\ 
m_{k,l}(\nat_B, \bx_{k,l})
&\coloneqq& \ell_{k,l}(\nat_B, \bx_{k,l}),
&& 1 \leq k < l \leq K.
\ee
We next verify that the assumptions placed on $m(\nat, \bx)$ in the general argument presented above
are met in the case of maximum likelihood estimation
for exponential-family local dependence random graphs. 

By Lemma \ref{lem:exp.fam.deriv},
\beno
\nabla_{\nat} \, m_{k,k}(\nat_W, \bx_{k,k})
\= (s_{k,k}(\bx_{k,k}) - \mbE_{\nat} \, s_{k,k}(\bX_{k,k}), \, \bm{0}_q), \s\\
\nabla_{\nat} \, m_{k,l}(\nat_B, \bx_{k,l})
\= (\bm{0}_p, \, s_{k,l}(\bx_{k,l}) - \mbE_{\nat} \, s_{k,l}(\bX_{k,l})), 
\ee
noting that $\nat = (\nat_W, \nat_B) \in \mbR^{p+q}$,
implying 
$\nabla_{\nat} \, m(\nat, \bx) = s(\bx) - \mbE_{\nat} \, s(\bx)$.
Observe that
\beno
\mbE\left[ \nabla_{\nat} \, m_{k,k}(\truth_W, \bX_{k,k}) \right] \= \bm{0}_{p+q},
&& 1 \leq k \leq K \s\\ 
\mbE\left[ \nabla_{\nat} \, m_{k,l}(\truth_B, \bX_{k,l}) \right] \= \bm{0}_{p+q},
&& 1 \leq k < l \leq K,
\ee
implying $\mbE \nabla_{\nat} \, m(\truth, \bX) = \bm{0}_{p+q}$.
Lemma \ref{lem:exp.fam.deriv} additionally establishes that 
\beno
\nabla_{\nat}^{2} \, m(\truth, \bx)
\= \var \, s(\bX)
\= \var \, (s(\bX) - \mbE \, s(\bX))
\= \var \, \nabla_{\nat} \, m(\truth, \bX),
\ee
implying $\bC = \var \, s(\bX) = \nabla_{\nat}^{2} \, m(\nat, \bx) = - \mbE \, \nabla_{\nat}^2 \, \ell(\nat, \bX)$,
which is non-singular for all $\nat \in \mB_2(\truth, \epsilon)$ by Assumption 2.  
Restricting to $\nat \in \mB_2(\truth, \epsilon)$,
we have verified all conditions placed on $m(\nat, \bx)$ in the general argument outlined above.
In this case, 
\eqref{eq:norm_key} can be expressed as 
\beno
\bS
&\overset{D}{=}&
\mcI(\truth)^{1/2}\,  (\truth - \nat) -  \mcI(\truth)^{-1/2} \, \bR,
\ee
where $\mcI(\truth) = \var \, s(\bX)$ is the Fisher information matrix 
evaluated 
at the data-generating parameter vector $\truth \in \mbR^{p+q}$. 
The local dependence assumption and the partitioning of 
the sufficient statistics vector $s(\bX) = (s_W(\bX_W), s_B(\bX_B))$ imply that  
\beno
\mcI(\truth)
\= \var \, s(\bX)
\= \left( \begin{matrix}
\mcI_W(\truth_W) & \bm{0}_{p,q} \s\\
\bm{0}_{q,p} & \mcI_{B}(\truth_B)
\end{matrix} \right),
\ee
where $\bm{0}_{d,r}$ is the $(d\times r)$-dimensional matrix consisting of all zeros,
with the definitions 
\beno
\mcI_W(\truth_W) 
&\coloneqq& \dsum_{k=1}^{K} \, - \mbE \, \nabla_{\nat_W}^2 \, \ell_{k,k}(\truth_W, \bX_{k,k}) \s\\
\mcI_{B}(\truth_B)
&\coloneqq& \dsum_{1 \leq k < l \leq K} \,  - \mbE \, \nabla_{\nat_B}^2 \, \ell_{k,l}(\truth_B, \bX_{k,l}).  
\ee

The proof is completed by establishing the following two additional results. 
\ben
\item[I.] {\bf Convergence rate of the multivariate normal approximation} 

\s

We establish the convergence rate of the multivariate normal approximation 
by bounding $\sum_{1 \leq k \leq l \leq K} \, \mbE \, \norm{\bY_{k,l}}_2^3$.  
In order to do so, 
we bound each term:  
\beno
\norm{\bY_{k,k}}_2
\= \norm{\mcI_W(\truth_W)^{-1/2} \, \left[ s_{k,k}(\bx_{k,k}) - \mbE_{\nat} \, s_{k,k}(\bX_{k,k}) \right]}_2 \s\\
&\leq& \mnorm{\mcI_W(\truth_W)^{-1/2}}_2 \, \norm{s_{k,k}(\bx_{k,k}) - \mbE_{\nat} \, s_{k,k}(\bX_{k,k})}_2 \s\\
&\leq& \dfrac{\norm{s_{k,k}(\bx_{k,k}) - \mbE_{\nat} \, s_{k,k}(\bX_{k,k})}_2}{\sqrt{K \, \minwt}}, 
\ee
using the bound on the spectral norm 
$\mnorm{\mcI_W(\truth_W)^{-1/2}}_2$ of the matrix $\mcI_W(\truth_W)^{-1/2}$: 
\beno
\mnorm{\mcI_W(\truth_W)^{-1/2}}_2
&\leq& \dfrac{1}{\sqrt{\lambda_{\min}(- \mbE \, \nabla_{\nat_W}^2 \ell(\truth, \bX))}}  
\= \dfrac{1}{\sqrt{K \, \minwt}}, 
\ee
which follows from Assumption 2,
because if $\lambda_{\min}(- \mbE \, \nabla_{\nat_W}^2 \, \ell(\truth, \bX))$ 
is the smallest eigenvalue of $\mcI_W(\truth_W) \coloneqq - \mbE \, \nabla_{\nat_W}^2 \ell(\truth, \bX)$,
then $1 \,/\, \lambda_{\min}(- \mbE \, \nabla_{\nat_W}^2 \, \ell(\truth, \bX))$ 
is the largest eigenvalue of $\mcI_W(\truth_W)^{-1}$. 
By Assumption 1, 
there exists $C_W > 0$,
independent of $N$, $p$, and $q$, such that  
\beno
\sup\limits_{\bx_{k,k} \in \mbX_{k,k}} \, 
\norm{s_{k,k}(\bx_{k,k})}_{\infty}
&\leq& C_W \, \dbinom{|\mA_k|}{2}, 
&&& 1 \leq k \leq K,  
\ee
which in turn implies the inequality
\beno
\norm{s_{k,k}(\bx_{k,k}) - \mbE_{\nat} \, s_{k,k}(\bX_{k,k})}_2
&\leq& \sqrt{p} \, \norm{s_{k,k}(\bx_{k,k}) - \mbE_{\nat} \, s_{k,k}(\bX_{k,k})}_{\infty} \s\\
&\leq& 2 \,  \sqrt{p} \, \norm{s_{k,k}(\bx_{k,k})}_{\infty} \s\\
&\leq& 2 \, \sqrt{p} \; C_W \, \dbinom{|\mA_k|}{2} \s\\
&\leq& \sqrt{p} \; C_W \, A_{\max}^2,
\ee
using the inequality $\tbinom{|\mA_k|}{2} \leq |\mA_k|^2 \,/\, 2 \leq A_{\max}^2 \,/\, 2$. 
Collecting the above bounds,
\beno
\norm{\bY_{k,k}}_2^3
&\leq& \left( \dfrac{\sqrt{p} \, C_W \, A_{\max}^2}{\sqrt{K \, \minwt}} \right)^3 
\= \dfrac{p^{3/2} \, C_W^3 \, A_{\max}^{6}}{K^{3/2} \, (\minwt)^{3/2}}, 
\ee
which implies
\beno
\dsum_{k=1}^{K} \, \mbE \, \norm{\bY_{k,k}}_2^3
&\leq& \dfrac{p^{3/2} \; C_W^3 \; A_{\max}^{6}}{K^{1/2} \, (\minwt)^{3/2}}. 
\ee
A similar argument will reveal,
using Assumption 1 once more, that 
\beno
\sup\limits_{\bx_{k,l} \in \mbX_{k,l}} \, 
\norm{s_{k,l}(\bx_{k,l}) - \mbE_{\nat} \, s_{k,l}(\bX_{k,l})}_{\infty}
&\leq& 2 \, C_B \; |\mA_k| \, |\mA_l|
&\leq& 2 \, C_B \; A_{\max}^2, 
\ee
for $1 \leq k < l \leq K$,
which will instead yield the bound 
\beno
\dsum_{1 \leq k < l \leq K} \, \mbE \, \norm{\bY_{k,l}}_2^3
&\leq& \dfrac{2 \, q^{3/2} \, C_B^3 \, A_{\max}^{6}}{\tbinom{K}{2}^{1/2} \, (\minbt)^{3/2}}
&\leq& \dfrac{2 \sqrt{2} \, q^{3/2} \, C_B^3 \, A_{\max}^{6}}{K \, (\minbt)^{3/2}}, 
\ee
noting there are $\tbinom{K}{2}$ between-block subgraphs $\bX_{k,l}$ ($1 \leq k < l \leq K$),
as opposed to $K$ within-block subgraphs $\bX_{k,k}$ ($1 \leq k \leq K$)
and using the bound $\tbinom{K}{2} \geq K^2 \,/\, 2$. 
Collecting terms and using the bound 
$N = \sum_{k=1}^{K} \, |\mA_k| \leq K \, A_{\max}$, 
\hide{
\beno
N 
\= \dsum_{k=1}^{K} \, |\mA_k|
&\leq& K \, A_{\max},  
\ee 
%defining $\Aavg \coloneqq K^{-1} \, \sum_{k=1}^{K} |\mA_k|$,
}
we obtain the bound  
\beno
\dsum_{1 \leq k \leq l \leq K} \, \mbE \, \norm{\bY_{k,l}}_2^3
&\leq& A_{\max}^{6} \, \left[ \dfrac{p^{3/2} \, C_W^3 }{K^{1/2} \, (\minwt)^{3/2} } +
\dfrac{2 \, \sqrt{2} \, q^{3/2} \, C_B^3}{K \, (\minbt)^{3/2}} \right] \s\\
&\leq& 2 \, \sqrt{2} \, A_{\max}^{6} \, 
\left[C_W^3 \, \sqrt{\dfrac{p^3}{K \, (\minwt)^{3}}}
+ C_B^3 \, \sqrt{\dfrac{q^3}{K^2 \, (\minbt)^{3}}} \; \right] \s\\
&\leq& 2 \, \sqrt{2} \, A_{\max}^{7}  \,
\left[C_W^3 \, \sqrt{\dfrac{p^3}{N \, (\minwt)^{3}}}
+ C_B^3 \, \sqrt{\dfrac{q^3}{N^2 \, (\minbt)^{3}}} \; \right].  
\ee
 
Thus, 
there exists a constant 
$C \coloneqq (2) (58) \, \sqrt{2} \, \max\{C_W^3, C_B^3\} > 0$, 
independent of $N$, $p$, and $q$,
and a random vector $\bDelta \coloneqq \mcI(\truth)^{-1/2} \, \bR$, 
such that, 
for all measurable convex sets $\mcC \subset \mbR^{p+q}$,
the error of the multivariate normal approximation is bounded above by
\beno
&& |\mbP(\mcI(\truth)^{1/2} \, (\mle - \truth) + \bDelta \,\in\,\mcC) - \Phi_d(\bZ_d \in \mcC)| \s\s\\
&\leq& 58 \, (p+q)^{1/4} \,  \dsum_{1 \leq k \leq l \leq K} \, \mbE \, \norm{\bY_{k,l}}_2^3 \s\\
&\leq& C \, (p+q)^{1/4} \, A_{\max}^{7}
\left[ \sqrt{\dfrac{p^3}{(\minwt)^{3} \, N}} 
+ \sqrt{\dfrac{q^3}{(\minbt)^{3} \, N^2}} \; \right].  
\ee

\s\s

\item[II.] {\bf Demonstrating that $\norm{\mcI(\truth)^{-1/2} \, \bR}_2$ is small with high probability.} \s

Recall that 
\beno
\mcI(\truth)
\= \left( \begin{matrix}
\mcI_W(\truth_W) & \bm{0}_{p,q} \s\\
\bm{0}_{q,p} & \mcI_B(\truth_B) 
\end{matrix} \right), 
\ee
where 
\beno
\mcI_W(\truth_W)
&\coloneqq& \dsum_{k=1}^K \, -\mbE \, \nabla_{\nat_W}^2 \, \ell_{k,k}(\truth_W, \bX_{k,k}) \s\\
\mcI_B(\truth_B)
&\coloneqq& \dsum_{1 \leq k < l \leq K} \, - \mbE \, \nabla_{\nat_B}^2 \, \ell_{k,l}(\truth_B, \bX_{k,l}), 
\ee
which implies that 
\beno
\mcI(\truth)^{-1}
\= \left( \begin{matrix}
\mcI_W(\truth_W)^{-1} & \bm{0}_{p,q} \s\\
\bm{0}_{q,p} & \mcI_B(\truth_B)^{-1}
\end{matrix} \right). 
\ee
Using Assumption 2, 
\be
\label{eq:eig_bound_1}
\lambda_{\min}(\mcI_W(\truth))
&=& K \, \minwt 
&>& 0 
\ee 
and 
\be
\label{eq:eig_bound_2}
\lambda_{\min}(\mcI_B(\truth_B))
&=& \dbinom{K}{2} \, \minbt 
&>& 0. 
\ee
Using \eqref{eq:eig_bound_1} and \eqref{eq:eig_bound_2}, 
we can bound $\norm{\mcI(\truth)^{-1/2} \, \bR}_2$ by 
\beno
\norm{\mcI(\truth)^{-1/2} \, \bR}_2^2
&=& \norm{\mcI_W(\truth)^{-1/2} \, \bR_W}_2^2 + \norm{\mcI_B(\truth)^{-1/2} \, \bR_B}_2^2 \s\\
&\leq& \mnorm{\mcI_W(\truth)^{-1/2}}_2^2 \, \norm{\bR_W}_2^2 + \mnorm{\mcI_B(\truth)^{-1/2}}^2 \, \norm{\bR_B}_2^2 \s\\
&\leq& \dfrac{\norm{\bR_W}_2^2}{K \, \minwt}
+ \dfrac{\norm{\bR_B}_2^2}{\tbinom{K}{2} \, \minbt},
\ee
where $\bR_W \coloneqq (R_1, \ldots, R_p)$ and $\bR_B \coloneqq (R_{p+1}, \ldots, R_{p+q})$.
As a result, 
\beno
\norm{\mcI(\truth)^{-1/2} \, \bR}_2
&\leq& \sqrt{ \dfrac{\norm{\bR_W}_2^2}{K \, \minwt}
+ \dfrac{\norm{\bR_B}_2^2}{\tbinom{K}{2} \, \minbt}}.
\ee
Applying Lemma \ref{lem:lag_rem}, 
\beno
\norm{\bR_W}_2^2
&\leq& p \, A_{\max}^{12} \, C_W^4 \, (C_W + 2)^2 \, K^2 \, \norm{\mle_W - \truth_W}_1^4 \s\\
\norm{\bR_B}_2^2
&\leq& 4 \, q \, A_{\max}^{12} \, C_B^4 \,  (C_B + 2)^2 \, \dbinom{K}{2}^2 \, \norm{\mle_B - \truth_B}_1^4.
\ee
Using the identity 
\beno
N 
\= \dsum_{k=1}^{K} \, |\mA_k|
\= K \, \dfrac{1}{K} \, \dsum_{k=1}^{K} \, |\mA_k|
\= K \, \Aavg,
\ee
with the definition $\Aavg \coloneqq K^{-1} \sum_{k=1}^{K} |\mA_k|$, 
we have the bound 
\beno
\norm{\mcI(\truth)^{-1/2} \, \bR}_2 
&\leq& C_3 \, A_{\max}^{6} 
\sqrt{ \dfrac{
p \, K^2 \, \norm{\mle_W - \truth_W}_1^4}
{K \, \minwt}
+ \dfrac{
q \,  \tbinom{K}{2}^2 \, \norm{\mle_B - \truth_B}_1^4}
{\tbinom{K}{2} \, \minbt}} \s\\
\= C_3 \, A_{\max}^{6}
\sqrt{ \dfrac{
p \, K \, \norm{\mle_W - \truth_W}_1^4}
{\minwt}
+ \dfrac{
q \,  \tbinom{K}{2} \, \norm{\mle_B - \truth_B}_1^4}
{\minbt}} \s\\
&\leq& 
C_3 \, A_{\max}^{6},
\sqrt{\dfrac{p \, N \, \norm{\mle_W - \truth_W}_1^4}{\Aavg \, \minwt}
+ \dfrac{q \, N^2 \, \norm{\mle_B - \truth_B}_1^4}{\Aavg^2 \, \minbt}}, 
\ee
where 
$C_3 \coloneqq 2 \, \max\{C_W^2 \, (C_W+ 2), \, C_B^2 \, (C_B + 2)\} > 0$
is a constant independent of $N$, $p$, and $q$.
By Theorem 2.1,
there exist constants $C_4 > 0$, $C_5 > 0$, and $N_0 \geq 3$,
independent of $N$, $p$, and $q$,
such that,
for all $N \geq N_0$ and  
with probability at least $1 - N^{-2}$, 
the maximum likelihood estimator exists, 
is unique, 
and satisfies 
\beno
\norm{\mle_W - \truth_W}_2
&\leq& C_4 \, \sqrt{\Aavg} \, \dfrac{\sqrt{\maxw}}{\minw} \, \sqrt{ \dfrac{p}{N} } 
\ee
and 
\beno
\norm{\mle_B - \truth_B}^2 
&\leq& C_5 \, \Aavg \, \dfrac{\sqrt{\maxb}}{\minb} \, \sqrt{ \dfrac{q}{N^2} }. 
\ee
As a result, 
\beno
\norm{\mle_W - \truth_W}_1
&\leq&
\sqrt{p} \, \norm{\mle_W - \truth_W}_2
&\leq& C_4 \, \sqrt{\Aavg} \, \dfrac{\sqrt{\maxw}}{\minw} \, \dfrac{p}{\sqrt{N}} \s\\
\norm{\mle_B - \truth_B}_1 
&\leq& \sqrt{q} \, \norm{\mle_B - \truth_B}_2
&\leq& C_5 \, \Aavg \, \dfrac{\sqrt{\maxb}}{\minb} \, \dfrac{q}{N}, 
\ee
which leads to the bound 
\beno
\norm{\mcI(\truth)^{-1/2} \, \bR}_2 
&\leq& 
C_3 \, A_{\max}^{6}
\sqrt{\dfrac{p \, N \, \norm{\mle_W - \truth_W}_1^4}{\Aavg \, \minwt}
+ \dfrac{q \, N^2 \, \norm{\mle_B - \truth_B}_1^4}{\Aavg^2 \, \minbt}} \s\s\\
&\leq& C_6 \,  A_{\max}^{6}
\sqrt{\dfrac{\Aavg \, (\maxw)^2 \, p^5}{(\minw)^4 \, \minwt \, N}
+ \dfrac{\Aavg^2 \, (\maxb)^2 \, q^5}{(\minb)^4 \, \minbt \, N^2}} \s\\
&\leq& C_6 \,  A_{\max}^{6}
\sqrt{\Aavg \, \dfrac{(\maxw)^2}{(\minw)^5} \, \dfrac{p^5}{N}
+ \Aavg^2 \, \dfrac{(\maxb)^2}{(\minb)^5} \, \dfrac{q^5}{N^2}}, 
\ee
defining $C_6 \coloneqq C_3 \, \max\{C_4^4, \, C_5^4\} > 0$. 
Thus,
there exists a constant $C \coloneqq C_6 > 0$, 
independent of $N$, $p$, and $q$,
such that 
$\bDelta \coloneqq \mcI(\truth)^{-1/2} \, \bR$ satisfies 
\beno
\mbP\left(\norm{\bDelta}_2 
\,\leq\,
C \,  A_{\max}^{6}
\sqrt{\Aavg \, \dfrac{(\maxw)^2}{(\minw)^5} \, \dfrac{p^5}{N}
+ \Aavg^2 \, \dfrac{(\maxb)^2}{(\minb)^5} \, \dfrac{q^5}{N^2}}
\; \right)
&\geq& 1 - \dfrac{1}{N^2}. 
\ee

\een

\s

\noindent {\bf Conclusion of proof.}
We have thus shown---recycling notation of constants---that 
there exist $C_1 > 0$, $C_2 > 0$, and $N_0 \geq 3$,
independent of $N$, $p$, and $q$,
and a random vector $\bDelta \in \mbR^{p+q}$
such that,
for all integers $N > N_0$ and all measurable convex sets $\mcC \subset \mbR^{p+q}$,
\beno
&& |\mbP(\mcI(\truth)^{1/2}\,  (\mle - \truth) + \bDelta \,\in\, \mcC) - \Phi_d(\bZ_d \in \mcC)|  \s\s\\
&\leq& C_1 \, (p+q)^{1/4} \, A_{\max}^{7} \, 
\left[ \sqrt{\dfrac{p^3}{(\minwt)^{3} \, N}}
+ \sqrt{\dfrac{q^3}{(\minbt)^{3} \, N^2}} \; \right], 
\ee
where $\bDelta$ satisfies 
\beno
\mbP\left(\norm{\bDelta}_2
\,\leq\,
C_2 \,  A_{\max}^{6}
\sqrt{\Aavg \, \dfrac{(\maxw)^2}{(\minw)^5} \, \dfrac{p^5}{N}
+ \Aavg^2 \, \dfrac{(\maxb)^2}{(\minb)^5} \, \dfrac{q^5}{N^2}}
\; \right)
&\geq& 1 - \dfrac{1}{N^2}.
\ee
\qed

\s

\subsection{Auxiliary results for Theorem 2.5}
\label{norm_aux}

We first recall a theorem due to \citet{Raic2019},
restated in Lemma \ref{lem:raic}. 

\begin{lemma}[Theorem 1.1, \citet{Raic2019}]
\label{lem:raic}
Consider a sequence of independent random vectors given by $\bW_1, \bW_2, \ldots \in \mbR^p$
with $\mbE \, \bW_i = 0$ for all $i \in \{1, 2, \ldots\}$.
Define
\beno
\bS_n
&\coloneqq& \dsum_{i=1}^{n} \bW_i,
&& n \in \{1, 2, \ldots\},
\ee
and assume that $\var \, \bS_n = \bI_p$.
Then,
for all measurable convex sets $\mathcal{C} \subset \mbR^p$,
\beno
\left| \mbP(\bS_n \in \mathcal{C}) - \Phi_p(\bZ_p \in \mathcal{C}) \right|
&\leq& (42 \, p^{1/4} + 16) \, \dsum_{i=1}^{n} \, \mbE \, \norm{\bW_i}_2^3,
\ee
where $\bZ_p$ is a multivariate normal random vector with mean vector $\bm{0}_p$
and covariance matrix $\bI_p$ and $\Phi_p$ is the corresponding probability distribution.
\end{lemma}

\llproof \ref{lem:raic}. 
The lemma is proved as Theorem 1.1 of \citet{Raic2019}.\
\qed

\subsection{Auxiliary results for Part II in the proof of Theorem 2.5}
\label{sup-sec:p2}

\begin{lemma}
\label{lem:lag_rem}
Under the assumptions of Theorem 2.5, 
\beno
\norm{\bR_W}_2^2
&\leq& p \, A_{\max}^{12} \, C_W^4 \, (C_W + 2)^2 \, K^2 \, \norm{\mle_W - \truth_W}_1^4 \s\\
\norm{\bR_B}_2^2
&\leq& 4 \, q \, A_{\max}^{12} \, C_B^4 \,  (C_B + 2)^2 \, \dbinom{K}{2}^2 \, \norm{\mle_B - \truth_B}_1^4, 
\ee
where $\norm{\bR_W}_2^2$ and $\norm{\bR_B}_2^2$ are the normed remainder terms 
in the proof of Theorem 2.5. 
\end{lemma}

\vspace{-.15cm}

\llproof \ref{lem:lag_rem}. 
We bound the remainder terms that arose out of the multivariate Taylor approximation in 
the proof of Theorem 2.5 using derivatives. 
Recall
that each of the remainder terms $R_i$ ($i \in \{1, \ldots, p+q\}$) in the Lagrange form is given by 
\be
\label{eq:remainder}
R_i
\=  \dsum_{j=1}^{p+q} \, \dfrac{1}{2} \,
\left[ \, \dfrac{\partial^2}{\partial \, \theta_j^2} \,
\left[\nabla_{\nat} \, m(\dot\nat^{(i)}, \bX) \right]_i \right] \,
(\theta_j - \theta_j^\star)^2 \s \\
&& +  \dsum_{1 \le j < r \le p+q} \,\dfrac{1}{2} \,
\left[ \, \dfrac{\partial^2}{\partial \, \theta_j \; \partial \, \theta_r} \,
\left[ \nabla_{\nat} \, m(\dot\nat^{(i)}, \bX) \right]_i \right] \,
(\theta_j - \theta^\star_j) \, (\theta_r - \theta^\star_r),
\ee
where $\dot\nat^{(i)} = t_i \, \nat + (1 - t_i) \, \truth$ (for some $t_i \in (0, 1)$, $i \in \{1, \ldots, p+q\}$).
If 
\beno
\sup\limits_{\nat \in \mbR^{p+q} \,:\, \norm{\nat - \truth}_1 \leq \norm{\mle - \truth}_1} \;
& \left|  \dfrac{\partial^2}{\partial \, \theta_j \; \partial \, \theta_r} \,
\left[\nabla_{\nat} \, m(\nat, \bX) \right]_i \right|
&\leq& M_i,
&& 1 \leq j \leq r \leq p,
\ee 
for all $i \in \{1, \ldots , p\}$ 
and 
\beno
\sup\limits_{\nat \in \mbR^{p+q} \,:\, \norm{\nat - \truth}_1 \leq \norm{\mle - \truth}_1} \;
& \left| \dfrac{\partial^2}{\partial \, \theta_j \; \partial \, \theta_r} \,
\left[ \nabla_{\nat} \, m(\nat, \bX) \right]_i \right|
&\leq& M_i,
&& 1 +p \leq j \leq r \leq p + q, 
\ee
for all $i \in \{1+p, \ldots, p+q\}$,
then the Lagrange remainders are bounded above by
\beno
|R_i|
&\leq& \begin{cases}
\dfrac{M_i}{2} \, \norm{\mle_W - \truth_W}_1^2 & \mbox{if } i \in \{1, \ldots, p\} \s\s\\ 
\dfrac{M_i}{2} \, \norm{\mle_B - \truth_B}_1^2  & \mbox{if } i \in \{p+1, \ldots, p+q\} 
\end{cases}
\ee
on the set
\beno
\left\{ \nat_W \in \mbR^p \,:\, \norm{\nat_W - \truth_W}_1 \leq \norm{\mle_W - \truth_W}_1 \right\}
\;\times\;  
\left\{ \nat_B \in \mbR^p \,:\, \norm{\nat_B - \truth_B}_1 \leq \norm{\mle_W - \truth_B}_1 \right\}.
\ee
For the rest of the proof,
assume that $\nat$ belongs to the above set.
By Assumption 2,
there exists $C_W > 0$ and $C_B > 0$,
independent of $N$, $p$, and $q$, 
such that
\beno
\sup\limits_{\bx_{k,k} \in \mbX_{k,k}}
\norm{s_{k,k}(\bx_{k,k})}_{\infty}
&\leq& C_W \, \dbinom{|\mA_k|}{2},
&&& 1 \leq k \leq K, 
\ee
and
\beno
\sup\limits_{\bx_{k,l} \in \mbX_{k,l}}
\norm{s_{k,l}(\bx_{k,l})}_{\infty}
&\leq& C_B \, |\mA_k| \, |\mA_l|,
&&& 1 \leq k < l \leq K, 
\ee
Lemmas \ref{lem:rem_bound} and \ref{lem:rem_bound_between}
establish that
\beno
\scalebox{1}{$
\left|  \dfrac{\partial^2}{\partial \, \theta_j \; \partial \, \theta_h} \,
\left[\nabla_{\nat} \, \ell(\nat, \bX) \right]_i \right|
$}
&\leq&
\scalebox{.95}{$
\begin{cases}
A_{\max}^6 \, C_W^3 \, (C_W + 2) \, K,
& (i,j,h) \in \{1, \ldots,p\}^3 \s\\
 2 \, A_{\max}^{6} \, C_B^2 \, (C_B + 2) \, \tbinom{K}{2},
& (i,j,h) \in \{p+1, \ldots, p+q\}^3 \s\\
0 & \mbox{otherwise}. \s\\
\end{cases}
$}
\ee
As a result,
when $m(\nat, \bX) = \ell(\nat, \bX)$ in the proof of Theorem 2.5,  
\beno
|R_i|
&\leq& \begin{cases}
A_{\max}^6 \, C_W^3 \, (C_W + 2) \, K \, 
\norm{\mle_W - \truth_W}_1^2, & 1 \leq i \leq p \s\s\\
2 \, A_{\max}^{6} \, C_B^2 \, (C_B + 2) \, \tbinom{K}{2} \, 
\norm{\mle_B - \truth_B}_1^2, & p+1 \leq i \leq p+q,
\end{cases}
\ee
which implies the bounds 
\beno
\norm{\bR_W}_2^2
&\leq& \dsum_{i=1}^{p} \, R_i^2
&\leq& p \, A_{\max}^{12} \, C_W^4 \, (C_W + 2)^2 \, K^2 \, \norm{\mle_W - \truth_W}_1^4 \s\\
\norm{\bR_B}_2^2
&\leq& \dsum_{i=p+1}^{p+q} \, R_i^2
&\leq& 4 \, q \, A_{\max}^{12} \, C_B^4 \,  (C_B + 2)^2 \, \dbinom{K}{2}^2 \, \norm{\mle_B - \truth_B}_1^4. 
\ee
\qed

\s\s

\begin{lemma}
\label{lem:rem_bound}
Consider an exponential-family local dependence random graph model
which satisfies Assumption 1.
%Assume that
%there exists a constant $C_W > 0$ such that,
%for all $k \in \{1, \ldots, K\}$, 
%\beno
%\sup\limits_{\bx_{k,k} \in \mbX_{k,k}} \, \norm{s_{k,k}(\bx_{k,k})}_{\infty}
%&\leq& C_W \, \dbinom{|\mA_k|}{2}. 
%\ee
Then,
for all $(i,j,h) \in \{1, \ldots, p\}^3$, 
\beno
\left| \dfrac{\partial^2 \, [\nabla_{\nat} \, \ell(\nat, \bx)]_i}{\partial \theta_{h} \, \partial \theta_{j}} \right|
&\leq& A_{\max}^{6} \, C_W^2 \, (C_W + 2). 
\ee
\end{lemma}

\llproof \ref{lem:rem_bound}.
By Lemma \ref{lem:exp.fam.deriv}, 
the second derivatives of the log-likelihood taken with respect to the natural parameters 
are equal to the variances and covariances of the sufficient statistics of the exponential family,
implying,
for all $(i,j) \in \{1, \ldots, p\}^2$, that   
\beno
\dfrac{\partial^2 \, \ell(\nat, \bx)}{\partial \theta_j \; \partial \theta_i}
\= 
\dfrac{\partial \, [\nabla_{\nat} \, \ell(\nat, \bx)]_i}{\partial \theta_{j}} 
\= \cov_{\nat}\left(\dsum_{k=1}^{K} \, s_{k,k,i}(\bX_{k,k}), \; \dsum_{k=1}^{K} \, s_{k,k,j}(\bX_{k,k}) \right), 
\ee
where $\cov_{\nat}$ denotes the covariance operator corresponding to the probability distribution $\mbP_{\nat}$.  
By the independence of the block-based subgraphs $\bX_{k,k}$ ($1 \leq k\leq K$),
\beno
\cov_{\nat}\left(\dsum_{k=1}^{K} \, s_{k,k,i}(\bX_{k,k}), \; \dsum_{k=1}^{K} \, s_{k,k,j}(\bX_{k,k}) \right)
\= \dsum_{k=1}^{K} \, \cov_{\nat}(s_{k,k,i}(\bX_{k,k}), \; s_{k,k,j}(\bX_{k,k})).
\ee
Taking $h \in \{1, \ldots, p\}$ and using the triangle inequality, 
we obtain the bound
\be
\label{eq:bound_rm}
\left| \dfrac{\partial^2 \, [\nabla_{\nat} \, \ell(\nat, \bx)]_i}{\partial \theta_{h} \, \partial \theta_{j}} \right|
&\leq& \dsum_{k=1}^{K} \, \left| \dfrac{\partial}{\partial \, \theta_{h}} \, 
\cov_{\nat}(s_{k,k,i}(\bX_{k,k}), \; s_{k,k,j}(\bX_{k,k}))\right|.
\ee
To proceed from here, 
we apply Lemma \ref{lem:3D}.
To do so, 
we verify the assumptions of Lemma \ref{lem:3D}.
By Assumption 1, 
there exists $C_W > 0$ such that  
\beno
\sup\limits_{\bx_{k,k} \in \mbX_{k,k}} \, 
\norm{s_{k,k}(\bx_{k,k})}_{\infty}
&\leq& C_W \, \dbinom{|\mA_k|}{2}
&\leq& \dfrac{C_W}{2} \, A_{\max}^2, 
&& k \in \{1, \ldots, K\}, 
\ee
which implies,
for all $k \in \{1, \ldots, K\}$,
that   
\beno
\sup\limits_{\bx_{k,k} \in \mbX_{k,k}} \, 
\norm{s_{k,k}(\bx_{k,k}) - \mbE_{\nat} \, s_{k,k}(\bX_{k,k})}_{\infty} 
&\leq& 2 \, C_W \, \dbinom{|\mA_k|}{2}
&\leq& C_W \, A_{\max}^2.  
\ee
Taking $U_1 = (C_W \,/\, 2)\, A_{\max}^2> 0 $ and $U_2 = C_W \, A_{\max}^2 > 0$
and applying Lemma \ref{lem:3D}, 
\beno
\left| \dfrac{\partial}{\partial \, \theta_{h}} \, \cov_{\nat}(s_{k,k,i}(\bX_{k,k}), \; s_{k,k,j}(\bX_{k,k}))\right|
&\leq&  \left( C_W \, A_{\max}^2  \right) \, 
\left( C_W \, A_{\max}^2  \right) \, 
\left(C_W \, A_{\max}^2  + 2 \right) \s\s\\
&\leq& A_{\max}^{6} \, C_W^2 \, (C_W + 2).  
\ee
Hence,
for all $\{i,j,h\} \subseteq \{1, \ldots, p\}$,
\beno
\left| \dfrac{\partial^2 \, [\nabla_{\nat} \, \ell(\nat, \bx)]_i}{\partial \theta_{h} \, \partial \theta_{j}} \right|
&\leq& \dsum_{k=1}^{K} \, A_{\max}^{6} \, C_W^2 \, (C_W + 2)
&\leq&  
A_{\max}^{6} \, C_W^2 \, (C_W + 2) \, K. 
\ee
\qed

\s

\begin{lemma}
\label{lem:rem_bound_between}
Consider an exponential-family local dependence random graph model
which satisfies Assumption 1. 
%Assume that
%there exists a constant $C_B > 0$ such that,
%for all $1 \leq k < l \leq K$,
%\beno
%\sup\limits_{\bx_{k,l} \in \mbX_{k,l}} \, \norm{s_{k,l}(\bx_{k,l})}_{\infty}
%&\leq& C_B \, |\mA_k| \, |\mA_l|.
%\ee
Then,
for all $(i,j,h) \in \{p+1, \ldots, p+q\}^3$,
\beno
\left| \dfrac{\partial^2 \, [\nabla_{\nat} \, \ell(\nat; \bx)]_i}{\partial \theta_{h} \, \partial \theta_{j}} \right|
&\leq& 2 \,  A_{\max}^{6} \, C_B^2 \, (C_B + 2) \, \dbinom{K}{2}.
\ee
\end{lemma}

\llproof \ref{lem:rem_bound_between}.
Lemma \ref{lem:rem_bound_between} is proved similarly to Lemma \ref{lem:rem_bound},
with the notable exception that the sum in \eqref{eq:bound_rm} is over
the index set $1 \leq k < l \leq K$,
for the between-block subgraphs.
As a result,
the factor of $K$ in the bound in Lemma \ref{lem:rem_bound}
is replaced with $\tbinom{K}{2}$.
The bound $|\mA_k| \, |\mA_l| \leq A_{\max}^2$ 
is used in place of $\tbinom{|\mA_k|}{2} \leq A_{\max}^2$,
resulting in an extra factor of $2$. 
The rest of the proof can be repeated unchanged, 
with the appropriate adjustments to indexing 
(e.g., 
using $C_B > 0$ in place of $C_W > 0$). 
\qed

\s

\s

\begin{lemma}
\label{lem:inq}
Let $a_1, a_2, b_1, b_2 \in \mbR$.
Then
\beno
|a_1 \, b_1 - a_2 \, b_2|
&\leq& |a_1| \, |b_1 - b_2| + |b_2| \, |a_1 - a_2|.
\ee
\end{lemma}

\llproof \ref{lem:inq}.
Write
\beno
|a_1 \, b_1 - a_2 \, b_2|
\= |a_1 \, b_1 - (a_2 - a_1 + a_1) \, b_2| \s\\ 
\= |a_1 \, b_1 - b_2 \, (a_2 - a_1) - a_1 \, b_2| \s\\
\= | a_1 \, (b_1 - b_2) \, - b_2 \, (a_2 - a_1) |  \s\\
&\leq& |a_1| \, |b_1 - b_2| + |b_2| \, |a_1 - a_2|.
\ee
\qed

\newpage 

\section{Proof of Theorem 2.7}
\label{sup-sec:prop_I}

Observe that both 
\beno
\widehat{\mcI}_W
&\coloneqq& \dfrac{1}{K} \, \dsum_{k=1}^{K} \, 
\left(s_{k,k}(\bX_{k,k}) - \bar{s}_W(\bX_W)\right) \,
\left(s_{k,k}(\bX_{k,k}) - \bar{s}_W(\bX_W)\right)^{\top} \s\s\\
\widehat{\mcI}_B
&\coloneqq& \dfrac{1}{\tbinom{K}{2}} \, \dsum_{1 \leq k < l \leq K} \, 
\left(s_{k,l}(\bX_{k,l}) - \bar{s}_B(\bX_B)\right) \, 
\left(s_{k,l}(\bX_{k,l}) - \bar{s}_B(\bX_B)\right)^{\top}
\ee
are unbiased estimators of 
\beno
\widetilde{\mcI}_W^\star 
&\coloneqq& \dfrac{\mbE[-\nabla_{\btheta_W}^2 \, \ell(\btheta^\star, \bX)]}{K}
\= \dfrac{1}{K} \, \dsum_{k=1}^{K} \, -\mbE \, \nabla_{\btheta_W}^2 \, \ell_{k,k}(\btheta_W^\star, \bX_{k,k}) \s\s\\
\widetilde{\mcI}_B^\star &\coloneqq& \dfrac{\mbE[-\nabla_{\btheta_B}^2 \, \ell(\btheta^\star, \bX)]}{\tbinom{K}{2}}
\= \dfrac{1}{\tbinom{K}{2}} \, \dsum_{1 \leq k < l \leq K} \, - \mbE \, \nabla_{\btheta_B}^2 \, \ell_{k,l}(\btheta_B^\star, \bX_{k,l}),
\ee
respectively,
and defining 
\beno
\bar{s}_W(\bX_W) &\coloneqq& \dfrac{1}{K} \, \dsum_{k=1}^{K} \, s_{k,k}(\bX_{k,k}) 
&&&\mbox{and}&&& 
\bar{s}_B(\bX_B) &\coloneqq& \dfrac{1}{\tbinom{K}{2}} \, \dsum_{1 \leq k < l \leq K} \, s_{k,l}(\bX_{k,l}).  
\ee 
For ease of presentation, 
we will write $\bar{s}_W \equiv \bar{s}_W(\bX_W)$
and $\bar{s}_B \equiv \bar{s}_B(\bX_B)$. 
Note that the Fisher information matrices 
$-\mbE \, \nabla_{\btheta_W}^2 \, \ell(\btheta, \bX)$ 
and $-\mbE \, \nabla_{\btheta_B}^2 \, \ell(\btheta, \bX)$
of canonical exponential families 
are the variance-covariance matrices of the vectors of sufficient statistics 
$s_W(\bX)$ and  $s_B(\bX)$,
respectively. 

We first consider the within-block case and the term $\widehat{\mcI}_W$;
we discuss extensions to prove the result for $\widehat{\mcI}_B$ afterwards. 
We can represent $\widehat{\mcI}_W$ by 
\beno
\widehat{\mcI}_W
\= \dfrac{1}{K} \, \dsum_{k=1}^{K} \, s_{k,k}(\bX_{k,k}) \, s_{k,k}(\bX_{k,k})^{\top} 
- \left( \dfrac{1}{K} \, \dsum_{k=1}^{K} \, s_{k,k}(\bX_{k,k}) \right) \, 
\left( \dfrac{1}{K} \, \dsum_{k=1}^{K} \, s_{k,k}(\bX_{k,k}) \right)^{\top} \s\\
\= \dfrac{1}{K} \, \dsum_{k=1}^{K} \, s_{k,k}(\bX_{k,k}) \, s_{k,k}(\bX_{k,k})^{\top} 
- \bar{s}_W \, \bar{s}_W^{\top}. 
\ee 
By the triangle inequality,   
\beno
\mnorm{\widehat{\mcI}_W - \widetilde{\mcI}_W}_2
&\leq& \left|\left|\left|
\dfrac{1}{K} \, \dsum_{k=1}^{K} \, s_{k,k}(\bX_{k,k}) \, s_{k,k}(\bX_{k,k})^{\top} 
- \widetilde{\bmu}_W \, \widetilde{\bmu}_W^{\top} - \widetilde{\mcI}_W
\right|\right|\right|_2 +  \left|\left|\left|  
\bar{s}_W \, \bar{s}_W^{\top} 
- \widetilde{\bmu}_W \, \widetilde{\bmu}_W^{\top}
\right|\right|\right|_2, 
\ee
defining 
\beno
\widetilde{\bmu}_W 
&\coloneqq& \dfrac{1}{K} \, \dsum_{k=1}^{K} \, \mbE \, s_{k,k}(\bX_{k,k}),
\ee
and in the between-block case, 
defining 
\beno
\widetilde{\bmu}_B
&\coloneqq& \dfrac{1}{\tbinom{K}{2}} \, \dsum_{1\leq k < l \leq K} \, \mbE \, s_{k,l}(\bX_{k,l}). 
\ee

The advantage of this is two-fold:
\ben
\item First, 
we can express
\beno
&& \dfrac{1}{K} \, \dsum_{k=1}^{K} \, s_{k,k}(\bX_{k,k}) \, s_{k,k}(\bX_{k,k})^{\top}
- \widetilde{\bmu}_W \, \widetilde{\bmu}_W^{\top} \s\\
\= \dfrac{1}{K} \, \dsum_{k=1}^{K} \, (s_{k,k}(\bX_{k,k}) - \widetilde{\bmu}_W) \, (s_{k,k}(\bX_{k,k}) - \widetilde{\bmu}_W)^{\top},
\ee
which is a sum of independent matrices, 
and \\  
\item Second,
\beno 
\mbE \, \left[ \dfrac{1}{K} \, \dsum_{k=1}^{K} \, s_{k,k}(\bX_{k,k}) \, s_{k,k}(\bX_{k,k})^{\top}
- \widetilde{\bmu}_W \, \widetilde{\bmu}_W^{\top} \right] 
\= \widetilde{\mcI}_W,
\ee 
meaning the statistic is an unbiased estimator of $\widetilde{\mcI}_W$. \\  
\een
Defining 
\beno
\mathring{\widehat{\mcI}_W}
&\coloneqq& \dfrac{1}{K} \, \dsum_{k=1}^{K} \, s_{k,k}(\bX_{k,k}) \, s_{k,k}(\bX_{k,k})^{\top}
- \widetilde{\bmu}_W \, \widetilde{\bmu}_W^{\top},
\ee
and in the between-block case 
\beno
\mathring{\widehat{\mcI}_B}
&\coloneqq& \dfrac{1}{\tbinom{K}{2}} \, \dsum_{1 \leq k < l \leq K} \, s_{k,l}(\bX_{k,l}) \, s_{k,l}(\bX_{k,l})^{\top}
- \widetilde{\bmu}_B \, \widetilde{\bmu}_B^{\top},
\ee
we can therefore apply the matrix Bernstein's inequality to obtain 
\beno
\mbP\left( \mnorm{\mathring{\widehat{\mcI}_W} - \widetilde{\mcI}_W}_2 \,\geq\, t \right)
&\leq& 2 \, p \, \exp\left( - \dfrac{K \, t^2}{2 \, p \, C_W^2 \; A_{\max}^4
(\widetilde{\lambda}^\star_{\max,W} 
+ 2 \, t \,/\, 3)} \right) \s\s\\
&\leq& 
2 \, \exp\left( - \dfrac{K \, t^2}{2 \, p \, C_W^2 \; A_{\max}^4
(\widetilde{\lambda}^\star_{\max,W}
+ t)} + \log(p)\right),
\ee
noting that,
by Assumption 1,
there exists constant $C_W > 0$,
independent of $N$, $p$, and $q$, 
such that 
\beno
\max\limits_{k \in \{1, \ldots, K\}} \, 
\sup\limits_{\bx_{k,k} \in \mbX_{k,k}} \,
&\norm{s_{k,k}(\bx_{k,k})}_{\infty}
&\leq& C_W \; \dbinom{|\mA_k|}{2}
&\leq& C_W \, A_{\max}^2,
\ee
which implies the bound 
\beno
\max\limits_{k \in \{1, \ldots, K\}} \,
\sup\limits_{\bx_{k,k} \in \mbX_{k,k}} \; 
\norm{s_{k,k}(\bx_{k,k}) - \widetilde{\bmu}_W}_2
&\leq& 2 \, \sqrt{p} \, C_W \, A_{\max}^2,  
\ee
using the triangle inequality 
%$\norm{s_{k,k}(\bx_{k,k}) - \widetilde{\bmu}_W}_2 \leq \norm{s_{k,k}(\bx_{k,k})}_{2} 
%+ \norm{\widetilde{\bmu}_W}_{2}$
and the inequality $\norm{\bz}_{2} \leq \sqrt{p} \, \norm{\bz}_{\infty}$ (for $\bz \in \mbR^p$). 
Choosing 
\beno
t 
\= \beta_W \, \sqrt{\dfrac{2 \, p \, C_W^2 \, A_{\max}^4 \, \widetilde{\lambda}^{\star}_{\max,W} \, \log(p)}{K}}
&>& 0,
\ee
for some $\beta_W > 0$ to be specified later, 
we obtain 
\beno
\mbP\left( \mnorm{\mathring{\widehat{\mcI}_W} - \widetilde{\mcI}_W}_2 \,\geq\, t \right)
&\leq& 2 \, \exp\left( - 
\dfrac{\beta_W^2 \, \log(p) \, \widetilde{\lambda}^{\star}_{\max,W}}{\widetilde{\lambda}^{\star}_{\max,W} + 
\beta_W \, \sqrt{2 \, p \, C_W^2 \, A_{\max}^4 \, \widetilde{\lambda}^{\star}_{\max,W} \, \log(p) \,/\, K}} + \log(p) \right),
\ee
By Assumption 4, 
 the largest block size
$A_{\max} \coloneqq \max\{|\mA_1|, \ldots, |\mA_K|\}$
satisfies
\beno
A_{\max}
&\leq& \min\left\{
\left( \dfrac{N \, \maxw}{\Aavg \, p^2} \right)^{1/4}, \quad
\left( \dfrac{N^2 \, \maxb}{4 \, \Aavg^2 \, q^2} \right)^{1/4}
\,\right\},
\ee
which,
using the identity $N = \Aavg \, K$, 
implies that 
\beno
\sqrt{\dfrac{p \, C_W^2 \, A_{\max}^4 \, \widetilde{\lambda}^{\star}_{\max,W} \, \log(p)}{K}}
&\leq& \widetilde{\lambda}^{\star}_{\max,W},
\ee
resulting in the inequality
\beno
&& \mbP\left( \mnorm{\mathring{\widehat{\mcI}_W} - \widetilde{\mcI}_W}_2 \;\geq\; 
\beta_W \, \sqrt{\dfrac{2 \, p \, C_W^2 \, A_{\max}^4 \, \widetilde{\lambda}^{\star}_{\max,W} \, \log(p)}{K}}
\right) \s\\
&\leq& 2 \, \exp\left( - \dfrac{\beta_W^2 \, \log(p)}{1 + \sqrt{2} \,C_W^2 \, \beta_W} + \log(p) \right) \s\\
\= 2 \, \exp\left(-\left( \dfrac{\beta_W^2}{1 + \sqrt{2} \, C_W^2 \, \beta_W} - 1 \right) \, \log(p) \right).  
\ee 
To obtain the desired probability guarantee, 
we require 
\beno
\dfrac{\beta_W^2}{1 + \sqrt{2} \, C_W^2 \, \beta_W} - 1  
\= 2,
\ee
which in turn requires a solution $\beta_W \in (0, \infty)$ 
to the quadratic equation 
\beno
\beta_W^2 - 3 \sqrt{2} \, C_W^2 \, \beta_W - 3 
\= 0. 
\ee
Using the quadratic formula, 
such a root,
which incidentally is independent of $N$, $p$, and $q$, 
is given by 
\beno
\beta_W 
\= \dfrac{3 \sqrt{2}}{2} \, C_W^2 + \dfrac{1}{2} \, \sqrt{18 \, C_W^4 + 12}
&>& 0. 
\ee
Thus,
there exists a constant $C_1 > 0$,
independent of $N$, $p$, and $q$, 
such that 
\beno
\mbP\left( \mnorm{\mathring{\widehat{\mcI}_W} - \widetilde{\mcI}_W}_2 \,<\, C_1 \, A_{\max}^2 \, \sqrt{\maxw} \, \sqrt{\dfrac{p \, \log(p)}{K}} 
\right)
&\geq& 1 - \exp(-2 \, p)
&\geq& 1 - \dfrac{2}{N^2},  
\ee
with the last inequality following from the assumption that  $p \geq \log(N)$.
Using the identity $N = \Aavg \, K$,
\beno
\mbP\left( \mnorm{\mathring{\widehat{\mcI}_W} - \widetilde{\mcI}_W}_2 
\,<\, 
C_1 \, A_{\max}^{2} \, \sqrt{\Aavg \, \maxw} \;\; \sqrt{\dfrac{p \, \log(p)}{N}}
\right)
&\geq& 1 - \dfrac{2}{N^2}. 
\ee  
Next, 
we handle the term 
$\mnorm{\bar{s}_W \, \bar{s}_W^{\top}
- \widetilde{\bmu}_W \, \widetilde{\bmu}_W^{\top}}_2$.  
We first use the inequality 
\beno
\mnorm{
\bar{s}_W \, \bar{s}_W^{\top} 
- \widetilde{\bmu}_W \, \widetilde{\bmu}_W^{\top}
}_2
&\leq& 
\mnorm{\bar{s}_W \, \bar{s}_W^{\top}
- \widetilde{\bmu}_W \, \widetilde{\bmu}_W^{\top}}_F
\ee
where $\mnorm{\,\cdot\,}_F$ denotes the Frobenius norm. 
By Lemma \ref{lem:mat_bound},
\beno
\mnorm{\bar{s}_W \, \bar{s}_W^{\top}
- \widetilde{\bmu}_W \, \widetilde{\bmu}_W^{\top}}_F
&\leq& 2 \, C_W  \, A_{\max}^2 \sqrt{p} \, \norm{\bar{s}_W - \widetilde{\bmu}_W}_2 
\= \dfrac{2 \, C_W \, A_{\max}^2 \sqrt{p}}{K} \; \norm{s_W(\bX) - \mbE \, s_W(\bX)}_2, 
\ee
noting that $\norm{s_{k,k}(\bx_{k,k})}_{\infty} \leq C_W \, A_{\max}^2$,
by Assumption 1 as discussed above. 
By Lemma \ref{lem:exp.fam.deriv}, 
\beno
\nabla_{\btheta_W} \, \ell(\truth, \bX) 
\= s_W(\bX) - \mbE \, s_W(\bX), 
\ee
which allows us to apply 
Lemma \ref{lem:conc_L2} with $\delta_W > 0$
to obtain 
\beno
\mbP\left( \norm{
s_W(\bX) - \mbE \, s_W(\bX)
}_2 \,\leq\, \delta_W \right)
\geq  1 - \exp\left( - \dfrac{\delta_W^2}{5 \, K \, \widetilde{\lambda}_{\max,W}^\star + C_W A_{\max}^2 \, \sqrt{p} \, \delta_W} + \log(5) \, p \right). 
\ee
This is close to the same inequality 
that we arrived at in the proof of Theorem 1.  
Choosing 
\beno
\delta_W 
\= \gamma_W \, \sqrt{p \, K \, \maxw} 
&>& 0,
\ee
for some $\gamma_W > 0$ to be specified, 
the probability 
\beno
\mbP\left( \norm{
s_W(\bX) - \mbE \, s_W(\bX)
}_2 \,\leq\, \gamma_W \, \sqrt{p \, K \, \maxw} \, \right) 
\ee
is bounded below by 
\beno
&& 1 - \exp\left( - \dfrac{\gamma_W^2 \, p \, K \, \maxw}
{5 \, K \, \maxw + C_W A_{\max}^2 \, p \, \sqrt{K} \, \maxw \, \gamma_W } + \log(5) \, p \right) \s\s\\
\= 1 - \exp\left( - \dfrac{\gamma_W^2 \, p \, K}{5 \, K + C_W A_{\max}^2 \, p \, \sqrt{K} \, \gamma_W} + \log(5) \, p \right). 
\ee
Under Assumption 4,
the largest block size
$A_{\max} \coloneqq \max\{|\mA_1|, \ldots, |\mA_K|\}$
satisfies
\beno
A_{\max}
&\leq& \min\left\{
\left( \dfrac{N \, \maxw}{\Aavg \, p^2} \right)^{1/4}, \quad
\left( \dfrac{N^2 \, \maxb}{4 \, \Aavg^2 \, q^2} \right)^{1/4}
\,\right\},
\ee
which implies that 
\beno
A_{\max}^2 \, p \, \sqrt{K} &\leq& K
&&&\mbox{since}&&&
A_{\max}
&\leq& \left(\dfrac{K}{p^2} \right)^{1/4}.
\ee
This in turn implies that 
\beno
\mbP\left( \norm{
s_W(\bX) - \mbE \, s_W(\bX)
}_2 \,\leq\, \gamma_W \sqrt{p \, K \, \maxw} \right)
&\geq& 1 - \exp\left( - \dfrac{\gamma_W^2 \, p}{5 + C_W  \gamma_W} + \log(5) \, p \right) \s\\
\= 1 - \exp\left( - \left( \dfrac{\gamma_W^2}{5 + C_W \gamma_W} + \log(5) \right) \, p \right). 
\ee
To obtain the desired probability guarantee, 
we require 
\beno
\dfrac{\gamma_W^2}{5 + C_W \gamma_W} + \log(5)
\= 2,
\ee
which in turn requires a solution $\gamma_W \in (0, \infty)$ to the quadratic equation 
\beno
\gamma_W^2 - C_W \, (2 + \log(5)) \, \gamma_W - 5 \, (2 + \log(5)) 
\= 0. 
\ee
Using the quadratic formula, 
such a root,
which is incidentally independent of $N$, $p$, and $q$, 
is given by 
\beno
\gamma_W 
\= \dfrac{C_W \, (2 + \log(5)) + \sqrt{C_W^2 \, (2 + \log(5))^2 + 20 \, (2 + \log(5))}}{2} 
&>& 0,
\ee
which in turn establishes there exists a constant $C_2 > 0$,
independent of $N$, $p$, and $q$,
such that  
\beno
\mbP\left( \norm{
s_W(\bX) - \mbE \, s_W(\bX)
}_2 \,\leq\, C_2 \, \sqrt{p \, K \, \maxw} \right)
&\geq& 1 - \exp(-2p)
&\geq& 1 - \dfrac{2}{N^2}, 
\ee
where the last inequality follows from the assumption that 
$p \geq \log(N)$. 
Finally,
\beno
\mnorm{\bar{s}_W \, \bar{s}_W^{\top} - \widetilde{\bmu}_W \, \widetilde{\bmu}_W^{\top}}_2
%&\leq& \dfrac{2 \, C_W \, A_{\max}^2 \sqrt{p}}{K} \, \norm{s_W(\bX) - \mbE \, s_W(\bX)}_2 \s\\ 
&\leq& C_3 \, A_{\max}^2 \, \sqrt{\maxw} \, \sqrt{\dfrac{p^2}{K}}, 
\ee
with probability at least $1 - 2 \, N^{-2}$,
defining $C_3 \coloneqq 2 \, C_W C_2 > 0$. 
Using the identity $N = \Aavg \, K$, 
\beno
\mnorm{\bar{s}_W \, \bar{s}_W^{\top} - \widetilde{\bmu}_W \, \widetilde{\bmu}_W^{\top}}_2
&\leq& C_3 \, A_{\max}^2 \, \sqrt{\Aavg \, \maxw} \, \sqrt{\dfrac{p^2}{N}},
\ee
with probability at least $ 1 - 2 \, N^{-2}$. 

Collecting results,
we have shown that the event   
\beno
\mnorm{\widehat{\mcI}_W - \widetilde{\mcI}_W}_2
&\leq& C \, A_{\max}^2 \, \sqrt{\Aavg \, \maxw} \;\; \left(  
\,  \sqrt{\dfrac{p \, \log(p)}{N}}
+ \sqrt{\dfrac{p^2}{N}} \right),
\ee
occurs with probability at least $1 - 2 \, N^{-2}$. 

We can prove a similar bound on $\mnorm{\widehat{\mcI}_B - \widetilde{\mcI}_B}_2$
by making appropriate adjustment to indexing of certain quantities and replacing $K$ by $\tbinom{K}{2}$ in all places,
establishing that the event  
\beno
\mnorm{\widehat{\mcI}_B - \widetilde{\mcI}_B}_2
&\leq& C \, A_{\max}^2 \, \Aavg \, \sqrt{\maxb} \;\; \left( 
\,  \sqrt{\dfrac{q \, \log(q)}{N^2}}
+ \sqrt{\dfrac{q^2}{N^2}} \right),
\ee 
occurs with probability at least $1 - 2 \, N^{-2}$.
\qed

\subsection{Auxiliary results for Theorem 2.7}
\label{sec:aux_2.7}

\begin{lemma}
\label{lem:mat_bound}
Let $\bv \in \mbR^p$, 
$\bw \in \mbR^p$,
and $M > 0$ 
be such that $\max\{\norm{\bv}_{\infty}, \, \norm{\bw}_{\infty}\} \leq M$.  
Then 
\beno
\mnorm{\bv \, \bv^{\top} - \bw \, \bw^{\top}}_F
&\leq& 2 \, M \sqrt{p} \, \norm{\bv - \bw}_2. 
\ee 
\end{lemma}

\llproof \ref{lem:mat_bound}. 
We start by writing 
\beno
\mnorm{\bv \, \bv^{\top} - \bw \, \bw^{\top}}_F^2
\= \dsum_{i=1}^{p} \, \dsum_{j=1}^{p} \, (v_i \, v_j - w_i \, w_j)^2 \s\\
\= M^4 \, \dsum_{i=1}^{p} \, \dsum_{j=1}^{p} \, \left( \dfrac{v_i \, v_j}{M^2} - \dfrac{w_i \, w_j}{M^2} \right)^2 \s\\
&\leq& M^4 \, \dsum_{i=1}^{p} \, \dsum_{j=1}^{p} \, \left( \left|\dfrac{v_i - w_i}{M}\right| + 
\left| \dfrac{v_j - w_j}{M} \right| \right)^2, 
\ee
where the inequality follows from Lemma \ref{lem:prod_to_sum}. 
Next,
noting that $M > 0$, 
we have  
\beno
M^4 \, \dsum_{i=1}^{p} \, \dsum_{j=1}^{p} \, \left( \left|\dfrac{v_i - w_i}{M}\right| +
\left| \dfrac{v_j - w_j}{M} \right| \right)^2
\= M^4 \, \dsum_{i=1}^{p} \, \dsum_{j=1}^{p} \, \left( \dfrac{|v_i - w_i|}{M} + \dfrac{|v_j - w_j|}{M}\right)^2 \s\\ 
\= M^4 \, \dsum_{i=1}^{p} \, \dsum_{j=1}^{p} \, \dfrac{(|v_i - w_i| + |v_j - w_j|)^2}{M^2} \s\\
&\leq& M^4 \, \dsum_{i=1}^{p} \, \dsum_{j=1}^{p} \, \left( \dfrac{(v_i - w_i)^2}{M^2 \,/\, 2} 
+ \dfrac{(v_j - w_j)^2}{M^2 \,/\, 2} \right) \s\\
\= 2 \, M^2 \,  \dsum_{i=1}^{p} \, \dsum_{j=1}^{p} \, \left( (v_i - w_i)^2 + (v_j - w_j)^2 \right) \s\\
\= 4 \, p \, M^2 \,  \dsum_{i=1}^{p} \, (v_i - w_i)^2 \s\\
\= 4 \, p \, M^2 \, \norm{\bv - \bw}_2^2,  
\ee
where the inequality follows by Titu's lemma.  
As a result, 
\beno
\mnorm{\bv \, \bv^{\top} - \bw \, \bw^{\top}}_F
&\leq& 2 \, M \sqrt{p} \, \norm{\bv - \bw}_2, 
\ee
for all $\bv \in \mbR^p$ and $\bw \in \mbR^2$ with $\norm{\bv}_{\infty} \leq M$ 
and $\norm{\bw}_{\infty} \leq M$. 
\qed

\begin{lemma}
\label{lem:prod_to_sum}
For every $(a_1, a_2, b_1, b_2) \in [-1, 1]^4$,
we have  
$|a_1 \, a_2 - b_1 \, b_2|
\leq |a_1 - b_1| + |a_2 - b_2|$. 
\end{lemma}

\llproof \ref{lem:prod_to_sum}. 
Start by defining a function $g : [0, 1] \mapsto [-1, 1]$ by 
\beno
g(t) 
\= (t \, a_1 + (1 - t) \, b_1) \, (t \, a_2 - (1 - t) \, b_2),
&& t \in [0, 1], 
\ee
for a given $(a_1, a_2, b_1, b_2) \in [-1, 1]^4$. 
Then,
by the product rule,  
\beno
g^{\prime}(t) 
\= (a_1 - b_1) \, (t \, a_2 - (1 - t) \, b_2)
+ (a_2 - b_2) \, (t \, a_1 + (1 - t) \, b_1).  
\ee
By the triangle inequality,
for all $t \in [0, 1]$,  
\be
\label{bbb}
|g^{\prime}(t)|
&\leq& |a_1 - b_1| \, |t \, a_2 - (1 - t) \, b_2| 
+ |a_2 - b_2| \, |t \, a_1 + (1 - t) \, b_1| \s\\
&\leq& |a_1 - b_1| + |a_2 - b_2|, 
\ee
as $(a_1, a_2, b_1, b_2) \in [-1, 1]^4$, 
by assumption,
ensuring
$|t \, a_i - (1 - t) \, b_i| \leq 1$ ($i \in \{1, 2\}$).
Lastly, 
\beno
|a_1 \, a_2 - b_1 \, b_2|
\= |g(1) - g(0)|
\= |g^{\prime}(t^\star)|
&\leq& |a_1 - b_1| + |a_2 - b_2|,
\ee 
where the mean value theorem guarantees the existence of some $t^\star \in [0, 1]$
and 
the upper-bound follows from \eqref{bbb} which holds for all $t \in [0, 1]$. 
\qed

\section{Auxiliary results for exponential families} 
\label{sup-sec:exp_fam}

%We prove various auxiliary results which establish some properties of exponential families. 

\begin{lemma}
\label{lem:exp.fam.deriv}
Consider a random vector $\bY$ with finite support $\mbY$ (i.e.,  $|\mbY| < \infty$)
and assume that the probability mass function $f_{\nat} : \mbY \mapsto (0, 1)$ 
belongs to an $m$-dimensional exponential family,
i.e.,
\beno
f_{\nat}(\by) 
\= h(\by) \, \exp\left(\langle \nat, \, s(\by) \rangle - \psi(\nat) \right),
&& \by \in \mbY, \; 
\nat \in \mbR^m. 
\ee
Then 
\beno
\nabla_{\nat} \, \psi(\nat)
\= \mbE_{\nat} \, s(\bY) \s\\
\nabla_{\nat} \, \ell(\nat, \by) 
\= s(\by) - \mbE_{\nat} \, s(\bY) \s\\
\nabla_{\nat}^2 \, \psi(\nat) 
\= - \nabla_{\nat}^2 \, \ell(\nat, \by)
\;\;=\;\; \var_{\nat} \, s(\bY).  
\ee
\end{lemma}

\llproof \ref{lem:exp.fam.deriv}. 
All results follow from Propositions 3.8 and 3.10 
of \citet{Sundberg2019}. \qed

\s

\begin{lemma}
\label{lem:sup_mle}
Consider a random vector $\bY$ with finite support $\mbY$ (i.e.,  $|\mbY| < \infty$)
and assume that the probability mass function $f_{\nat} : \mbY \mapsto (0, 1)$
belongs to an $m$-dimensional exponential family,
i.e.,
\beno 
f_{\nat}(\by) 
\= h(\by) \, \exp\left(\langle \nat, \, s(\by) \rangle - \psi(\nat) \right), 
&& \by \in \mbY, \; \nat \in \mbR^m. 
\ee
Then  
\beno
\nabla_{\nat} \, \ell(\nat, \bY)
- \mbE \, \nabla_{\nat} \, \ell(\nat, \bY)
\= s(\bY) - \mbE \, s(\bY),
&& \nat \in \mbR^m, 
\ee
and 
\beno
\sup\limits_{\nat \in \mbR^m} \, 
\norm{\nabla_{\nat} \, \ell(\nat, \bY)
- \mbE \, \nabla_{\nat} \, \ell(\nat, \bY)}_{\infty}
\= \norm{s(\bY) - \mbE \, s(\bY)}_{\infty}
\= \norm{\nabla_{\nat} \, \ell(\truth, \bY)}_{\infty}. 
\ee
\end{lemma}

\llproof \ref{lem:sup_mle}.
Applying Lemma \ref{lem:exp.fam.deriv},
$\nabla_{\nat} \, \ell(\nat, \by) = s(\by) - \mbE_{\nat} \, s(\bY)$.
Hence, 
for all $\nat \in \mbR^m$,
\beno
\nabla_{\nat} \, \ell(\nat, \bY)
- \mbE \, \nabla_{\nat} \, \ell(\nat, \bY)
\= s(\bY) - \mbE_{\nat} \, s(\bY) - \mbE \, s(\bY) + \mbE \, \mbE_{\nat} \, s(\bY) 
\= s(\bY) - \mbE \, s(\bY), 
\ee
which implies the final result 
\beno
\sup\limits_{\nat \in \mbR^m} \,
\norm{\nabla_{\nat} \, \ell(\nat, \bY)
- \mbE \, \nabla_{\nat} \, \ell(\nat, \bY)}_{\infty}
\= \sup\limits_{\nat \in \mbR^m} \, 
\norm{s(\bY) - \mbE \, s(\bY)}_{\infty} 
\= \norm{s(\bY) - \mbE \, s(\bY)}_{\infty}. 
\ee
\qed

\s

\begin{lemma}
\label{lem:3D}
Let $\bY$ be a random vector with finite support $\mbY$ 
(i.e., $|\mbY| < \infty$)
 and   
assume that the distribution of $\bY$ belongs to an exponential family with with probability mass functions of the form
\beno
f_{\nat}(\by)
\= h(\by) \, \exp(\langle \nat, \, \by \rangle - \psi(\nat)),
&& \by \in \mbY, \; \nat \in \mbR^m.
\ee
Assume that there exist constants $U_1 > 0$ and $U_2 > 0$ such that,
for all $ t \in \{1, \ldots, m\}$,
\beno
|Y_t| &\leq& U_1
&&\mbox{and}&& |Y_t - \mbE_{\nat} \, Y_t| &\leq& U_2,
&&& \mbP\mbox{-almost surely}. 
\ee
Then,
for all $(i,j,h) \in \{1, \ldots, m\}^3$,
\beno
\left| \dfrac{\partial}{\partial \, \theta_h} \, \cov_{\nat}(Y_i, \, Y_j) \right|
&\leq& 2 \, U_1 \, U_2 \, (U_2 + 2).
\ee
\end{lemma}

\llproof \ref{lem:3D}.
Let $(i,j,h) \in \{1, \ldots, m\}^3$
and define $\mu_t(\nat)  \coloneqq \mbE_{\nat} \, Y_t$ ($t \in \{1, \ldots, m\}$. 
Then
\beno
\dfrac{\partial}{\partial \, \theta_h} \, \cov_{\nat}(Y_i, \, Y_j)
\= \dfrac{\partial}{\partial \, \theta_h} \, \mbE_{\nat}[Y_i \, Y_j - \mu_i(\nat) \, \mu_j(\nat)] \s\\
\= \dfrac{\partial}{\partial \, \theta_h} \, \dsum_{\by \in \mbY} \, 
\left[ y_i \, y_j - \mu_i(\nat) \, \mu_j(\nat) \right] \, f_{\nat}(\by) \s\\
\= \dsum_{\by \in \mbY} 
\left[ f_{\nat}(\by) \dfrac{\partial}{\partial \, \theta_h} \left[ y_i \, y_j - \mu_i(\nat) \, \mu_j(\nat) \right] + 
\left[ y_i \, y_j - \mu_i(\nat) \, \mu_j(\nat) \right] \, \dfrac{\partial}{\partial \, \theta_h} \,  f_{\nat}(\by) 
\right].  
\ee
Using Lemma \ref{lem:exp.fam.deriv} and applying the chain rule, 
\beno
\dfrac{\partial}{\partial \, \theta_h}
\, f_{\nat}(\by)
\= \dfrac{\partial}{\partial \, \theta_h} \, h(\by) \, \exp(\langle \nat, \, \by \rangle - \psi(\nat)) 
\= \left[y_h - \mu_h(\nat) \right] \, f_{\nat}(\by). 
\ee
Hence, 
\beno
\dfrac{\partial}{\partial \, \theta_h} \, \cov_{\nat}(Y_i, \, Y_j)
\= \mbE_{\nat}[(Y_i \, Y_j - \mu_i(\nat) \, \mu_j(\nat)) \, (Y_h - \mu_h(\nat))]
- \dfrac{\partial}{\partial \, \theta_h} \left[ \mu_i(\nat) \, \mu_j(\nat) \right].
\ee
We next compute,
using Lemma \ref{lem:exp.fam.deriv},  
\beno
\dfrac{\partial}{\partial \, \theta_h} \left[ \mu_i(\nat) \, \mu_j(\nat) \right]
\= \mu_i(\nat) \, \cov_{\nat}(Y_j,\, Y_h) + \mu_j(\nat) \, \cov_{\nat}(Y_i, \, Y_h) \s\\
\= \mu_i(\nat) \, \mbE_{\nat}[Y_j \, Y_h - \mu_j(\nat) \, \mu_h(\nat)]
+ \mu_j(\nat) \, \mbE_{\nat}[Y_i \, Y_h - \mu_i(\nat) \, \mu_h(\nat)].
\ee
By Jensen's inequality and the triangle inequality
\beno
\left| \dfrac{\partial}{\partial \, \theta_h} \, \cov_{\nat}(Y_i, \, Y_j) \right|
&\leq& \mbE_{\nat}\left[\left|(Y_i \, Y_j - \mu_i(\nat) \, \mu_j(\nat))\right| \, \left|(Y_h - \mu_h(\nat)) \right| \right] \\
&&+ \;\; |\mu_i(\nat)| \; \mbE_{\nat} \left| Y_j \, Y_h - \mu_j(\nat) \, \mu_h(\nat) \right| \s\\
&&+\;\; |\mu_j(\nat)| \; \mbE_{\nat} \left| Y_i \, Y_h - \mu_i(\nat) \, \mu_h(\nat) \right|.
\ee
The assumption there exist constants 
$U_1 > 0$ 
and 
$U_2 > 0$
such that 
$|Y_t| \leq U_1$ for all $t \in \{1, \ldots, m\}$
and 
$|Y_t - \mu_t(\nat)| \leq U_2$ ($t \in \{1, \ldots, m\}$)
hold $\mbP$-almost surely 
implies that $|\mu_t(\nat)| \leq U_1$ for all $t \in \{1, \ldots, m\}$, 
and,
through an application of Lemma \ref{lem:inq},
that 
\beno
\left| Y_j \, Y_h - \mu_j(\nat) \, \mu_h(\nat) \right|
&\leq& |Y_j| \, |Y_h - \mu_h(\nat)| + |\mu_h(\nat)| \, |Y_j - \mu_j(\nat)|
&\leq& 2 \, U_1 \, U_2. 
\ee
Hence,
\beno
\left| \dfrac{\partial}{\partial \, \theta_h} \, \cov_{\nat}(Y_i, \, Y_j) \right|
&\leq&
2 \, U_1 \, U^2_2 +
4 \, U_1 \, U_2
\= 2 \, U_1 \, U_2 \, (U_2 + 2).
\ee
\qed

\newpage

\begin{lemma}
\label{lem:exp_deriv}
Consider a random vector $\bY$ with finite support $\mbY$ (i.e., $|\mbY| < \infty$)
and assume that the  distribution of $\bY$ belongs to an exponential family
with probability mass functions of the form 
\beno
\mbP_{\nat}(\bY = \by)
\= h(\by) \, \exp\left( \langle \nat, \, s(\by) \rangle - \psi(\nat) \right),
&& \by \in \mbY, \; \nat \in \mbR^m. 
\ee
Then for all functions $f : \mbY \mapsto \mbR$,
\beno
\dfrac{\partial}{\partial \, \theta_i} \, \mbE_{\nat} \, f(\bY)
\= \mbE_{\nat}[ f(\bY) \, (s_i(\bY) - \mbE_{\nat} \, s_i(\bY))],
\ee
for all $i \in \{1, \ldots, m\}$.
\end{lemma}

\llproof \ref{lem:exp_deriv}.
Write
\beno
\dfrac{\partial}{\partial \, \theta_i} \, \mbE_{\nat} \, f(\bY)
\= \dfrac{\partial}{\partial \, \theta_i} \, \dsum_{\by \in \mbY} \, f(\by) \,
h(\by) \, \exp\left( \langle \nat, \, s(\by) \rangle - \psi(\nat) \right) \s\\
\= \dsum_{\by \in \mbY} \, f(\by) \, h(\by) \, 
\left[ \dfrac{\partial}{\partial \, \theta_i} \,
\exp\left( \langle \nat, \, s(\by) \rangle - \psi(\nat) \right) \right] \s\\
\= \dsum_{\by \in \mbY} \, f(\by) \, h(\by) \, \exp\left( \langle \nat, \, s(\by) \rangle - \psi(\nat) \right) \,
(s_i(\by) - \mbE_{\nat} \, s_i(\bY)) \s\\
\= \mbE_{\nat}[f(\bY) \, (s_i(\bY) - \mbE_{\nat} \, s_i(\bY))],
\ee
as applying Lemma \ref{lem:exp.fam.deriv} shows that 
\beno
\dfrac{\partial}{\partial \, \theta_i} \, \psi(\nat)
\= \mbE_{\nat} \, s_i(\bX),
&& i = 1, \ldots, m.
\ee
\qed

%\s

%\input{cond.exp.fam.tex}

%\s

%\input{cond.exp.fam.deriv.tex}

\bibliographystyle{imsart-nameyear.bst}
\bibliography{library}

\end{document}